\documentclass[11pt, twoside]{amsart}

\usepackage[utf8]{inputenc}
\usepackage[T1]{fontenc}
\usepackage{amssymb,amsmath,amstext}
\usepackage{hyperref, graphicx}

\newtheorem{theor}{Theorem}[section]
\newtheorem{lem}[theor]{Lemma}
\newtheorem{defin}[theor]{Definition}

\newtheorem{prop}[theor]{Proposition} 

\newtheorem{notation}[theor]{Notation}
\newtheorem{exam}[theor]{Example}

\newtheorem{cor}[theor]{Corollary}
\newtheorem{rem}[theor]{Remark}

\newtheorem{assump}[theor]{Assumption}

\numberwithin{equation}{section}

\newcommand{\mr}{\mathrm}

\newcommand{\es}{\emptyset}

\newcommand{\uhr}{\upharpoonright}

\newcommand{\nts}{\negthickspace}
\newcommand{\uhrc}{\nts \upharpoonright \nts}

\newcommand{\mcA}{\mathcal{A}}
\newcommand{\mcB}{\mathcal{B}}
\newcommand{\mcC}{\mathcal{C}}
\newcommand{\mcD}{\mathcal{D}}

\newcommand{\mcL}{\mathcal{L}}

\newcommand{\mfG}{\mathfrak{G}}

\newcommand{\mbE}{\mathbf{E}}

\newcommand{\mbP}{\mathbf{P}}

\newcommand{\mbW}{\mathbf{W}}
\newcommand{\mbX}{\mathbf{X}}
\newcommand{\mbY}{\mathbf{Y}}
\newcommand{\mbZ}{\mathbf{Z}}

\newcommand{\mbbP}{\mathbb{P}}
\newcommand{\mbbN}{\mathbb{N}}
\newcommand{\mbbR}{\mathbb{R}}

\newcommand{\msfP}{\mathsf{P}}

\newcommand{\rng}{\mathrm{rng}}

\title[Conditional probability logic]
{Conditional probability logic, lifted bayesian networks, and almost sure quantifier elimination}

\author{Vera Koponen}

\address{Vera Koponen, Department of Mathematics, Uppsala University, Box 480,
75106 Uppsala, Sweden.}

\email{vera.koponen@math.uu.se}

\date{18 August 2021}

\begin{document}

\maketitle

\begin{abstract}
We introduce a formal logical language, called {\em conditional probability logic (CPL)}, which extends first-order logic
and which can express probabilities, conditional probabilities and which can compare conditional probabilities.
Intuitively speaking, although formal details are different, CPL can express  the same kind of statements
as some languages which have been considered in the artificial intelligence community.
We also consider a way of making precise the notion of {\em lifted Bayesian network}, where this notion is a type of 
(lifted) probabilistic graphical model used in machine learning, data mining and artificial intelligence.
A lifted Bayesian network (in the sense defined here) determines, in a natural way, a probability distribution on 
the set of all structures (in the sense of first-order logic) with a common finite domain $D$.
Our main result 
(Theorem~\ref{main result on quantifier elimination})
is that for every ``noncritical'' CPL-formula $\varphi(\bar{x})$ there is a quantifier-free formula $\varphi^*(\bar{x})$
which is ``almost surely'' equivalent to $\varphi(\bar{x})$ as the cardinality of $D$ tends towards infinity.
This is relevant for the problem of making probabilistic inferences on large domains $D$, because (a)
the problem of evaluating, by ``brute force'', the probability of $\varphi(\bar{x})$ being true for some sequence $\bar{d}$ of elements from $D$
has, in general, (highly) exponential time complexity in the cardinality of $D$, and (b)
the corresponding probability for the quantifier-free
$\varphi^*(\bar{x})$ depends only on the lifted Bayesian network and not on $D$.
Some conclusions regarding the computational complexity of finding $\varphi^*$ are given in
Remark~\ref{remark about computational complexity}.
The main result has two corollaries, one of which is a convergence law (and zero-one law) for noncritial CPL-formulas.
\end{abstract}

\section{Introduction}

\noindent
We consider an extension of first-order logic which we call {\em conditional probability logic}
(Definition~\ref{definition of CPL}), abbreviated CPL,
with which it is possible to express statements about probabilities, conditional probabilities, and
to compare conditional probabilities which makes it possible to express
statements about the (conditional) independence (or dependence) of events or random variables.
Remarks~\ref{remark on abbreviations}, \ref{remark on expressing independence} and Example~\ref{example 1} below 
illustrate this.
The semantics of CPL deals only with finite structures and assumes that all elements in a structure are equally likely, so
(conditional) probabilities correspond to proportions.
Quite similar formal languages, which aim at expressing the same sort of statements, have been studied within the field of artificial intelligence by
Halpern \cite[Section 2]{Hal90} and Bacchus et. al. \cite[Definition 4.1]{Bacchus96}.
CPL is more expressive than the probability logic $\mcL_{\omega P}$ considered by Keisler and Lotfallah in~\cite{KL} 
(which cannot express {\em conditional} probabilities) 
and our first theorem (Theorem~\ref{main result on quantifier elimination})
is a generalization of their main result~\cite[Theorem~4.9]{KL}, both in the sense that the language considered here is more expressive and that
we consider a wider range of probability distributions.

A {\em graphical model} for a probability distribution and a set of random variables is a ``graphical'' way of describing
the conditional dependencies and independencies between the random variables. In such a probabilistic model the random variables 
are also viewed as the vertices of a directed or undirected graph where edges indicate conditional dependencies and independencies
\cite{BK, Pea}.
The notion of a Bayesian network is one of the most well-known graphical models.
A {\em Bayesian network} $\mfG$ for a probability space $(S, \mu)$ and binary random variables
$X_1, \ldots, X_n$
is determined by the following data where, for any distinct $i_1, \ldots, i_k \in \{1, \ldots, n\}$ and
values $x_1, \ldots, x_k$ that $X_{i_1}, \ldots, X_{i_k}$ can take, the tuple $(x_1, \ldots, x_k)$ denotes the
event that $X_{i_1} = x_1, \ldots, X_{i_k} = x_k$:
\begin{enumerate}
\item A (not necessarily connected) directed acyclic graph (DAG), also denoted $\mfG$, with vertex set $V = \{X_1, \ldots, X_n\}$ such that if there is
an arrow (directed arc) from $X_i$ to $X_j$ then $i < j$.

\item To each vertex $X_i \in V$, if $X_{j_1}, \ldots, X_{j_k}$ are the parents of $X_i$ (in the DAG)
and $x_i, x_{j_1}, \ldots, x_{j_k}$ are values that $X_i, X_{j_k}, \ldots, X_{j_k}$, respectively, can take,
then the conditional probability $\mbbP(x_i \ | \ x_{j_1}, \ldots, x_{j_k})$ is specified\footnote{
The expression $\mbbP(x_i \ | \ x_{j_1}, \ldots, x_{j_k})$ denotes the conditional probability
that $X_i = x_i$ given that $X_{j_1} = x_{j_1}, \ldots, X_{j_k} = x_{j_k}$} in such a way that
the following holds:
\begin{itemize}
\item[(a)] For each $j$, the set of parents of $X_j$, denoted $par(X_j)$, is minimal (with respect to set inclusion) with the property that
for every $i < j$, $X_i$ and $X_j$ are conditionaly independent over $par(X_j)$.
The conditional independence means that if $par(X_j) = \{X_{l_1}, \ldots, X_{l_k}\}$, then for all possible values 
$x_i, x_j, x_1, \ldots, x_k$ of $X_i, X_j, X_{l_1}, \ldots, X_{l_k}$, respectively,
we have 
\[
\mbbP(x_i \ | \ x_j, x_1, \ldots, x_k) = \mbbP(x_i \ | \ x_1, \ldots, x_k)
\]
whenever both sides are defined.\footnote{Or equivalently, we could require that
$\mbbP(x_i, x_j \ | \ x_1, \ldots, x_k) =  \mbbP(x_i \ | \ x_1, \ldots, x_k)  \mbbP(x_j \ | \ x_1, \ldots, x_k)$
whenever all factors are defined.}

\item[(b)] The joint probability distribution on $X_1, \ldots, X_n$ 
is determined by the conditional probabilities associated with the vertices of $\mfG$.
More precisely: The event $(x_1, \ldots, x_n)$ can be computed recursively by repeatedly using the
following identities which hold for any choice of distinct $i_1, \ldots, i_k \in \{1, \ldots, n\}$:
\begin{align*}
&\mbbP(x_{i_1}, \ldots, x_{i_k}) = 
\mbbP(x_{i_k} \ | \ x_{i_1}, \ldots, x_{i_{k-1}}) \cdot \mbbP(x_{i_1}, \ldots, x_{i_{k-1}}), \\
&\mbbP(x_{i_k} \ | \ x_{i_1}, \ldots, x_{i_{k-1}}) = \mbbP(x_{i_k} \ | \ x_{j_1}, \ldots, x_{j_m})\\
&\text{ if } j_1, \ldots, j_m \in \{i_1, \ldots, i_{k-1}\} \text{ and } par(X_{i_k}) = \{X_{j_1}, \ldots, X_{j_m}\}.
\end{align*}
\end{itemize}
\end{enumerate}

\noindent
If $\mfG$ is a Bayesian network as defined above, then it follows 
(from e.g.  \cite[Definition~1.2.1 and Theorems 1.2.6, 1.2.7]{Pea})
that
\begin{itemize}
\item[(i)] For every $X_j \in V$, $X_j$ and the set of all predecessors of $X_j$ are conditionally independent over
$par(X_j)$.
\item[(ii)] For every $X_j \in V$, $X_j$ and the set of all nondescendants of $X_j$
(except $X_j$ itself) are conditionally independent over $par(X_j)$.
\end{itemize}
Moreover: if condition~(i) or condition~(ii) holds, then $\{X_1, \ldots, X_n\}$ can be ordered so that conditions~(a)
and~(b) above hold without changing the arrows of the DAG.

Graphical models are used in machine learning, data mining and artificial intelligence in (probability based) learning and inference making.
To illustrate this by a very simple example, suppose that we have a finite set $A$ of some kind of objects
and properties $P, Q$ and $R$ which objects in $A$ may, or may not, have.
We can view $A$ as a ``training set''. The training set can be formalized as a $\sigma$-structure with domain $A$
where $\sigma = \{P, Q, R\}$ and $P, Q$ and $R$ are also viewed as unary relation symbols.
Let $\mu$ be a probability distribution on $A$ and let binary random variables $X, Y, Z : A \to \{0, 1\}$ be defined by
$X(a) = 1$ if $a$ has the property $P$ and $X(a) = 0$ otherwise (for every $a \in A$); 
$Y(a) = 1$ if $a$ has the property $Q$ and $Y(a) = 0$ otherwise; and analogously for $Z$ and $R$.
Suppose that, after some ``learning'', we have found a Bayesian network $\mfG$ for $(A, \mu)$ and $X, Y, Z$ such that its DAG is as illustrated
\begin{figure}[h!]
\begin{center}
\includegraphics[scale=1]{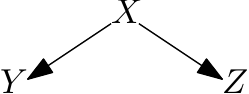}
\end{center}
\end{figure}
and the (conditional) probabilities $\mu(X = 1)$, $\mu(Y = 1 \ | \ X = 1)$, $\mu(Y = 1 \ | \ X = 0)$,
$\mu(Z = 1 \ | \ X = 1)$ and $\mu(Z = 1 \ | \ X = 0)$ are specified.
(In real applications, it is unlikely that a relatively simple probabilistic model,which is desirable for computational efficiency, 
fits the training data completely
and usually this is not even the goal because one wants to avoid so-called ``overfitting''; so one can view the 
Bayesian network as a reasonable approximation of the training data.)
An application of the Bayesian network $\mfG$ is to make predictions about probabilities on some other finite domain $B$.
Let us now make the following assumptions, partly based on $\mfG$ but where the independency assumptions between different objects are imposed.
Every $b \in B$ has the probability $\mu(X = 1)$ of having property $P$, independently of what the case is for other $b' \in B$.
For every $b \in B$, if $b$ has the property $P$ then the probability that $b$ also has the property $Q$ ($R$) is $\mu(Y = 1 \ | \ X = 1)$
($\mu(Z = 1 \ | \ X = 1)$),
independently of what the case is for other elements in $B$,
and if $b$ does not have the property $P$ then the probability that $b$ has $Q$ ($R$)  is $\mu(Y = 1 \ | \ X = 0)$ 
($\mu(Z = 1 \ | \ X = 0)$), independently of what the case
is for other elements.
Based on this we can define a probability distribution 
(as in Definition~\ref{the probability distribution induced by the BN}) 
on the set $\mbW_B$ of all $\sigma$-structures with domain $B$, where
each member of $\mbW_B$ represents a ``possible scenario'' or ``possible world''.
For every formula $\varphi(x_1, \ldots, x_k)$ of conditional probability logic and any choice
of $b_1, \ldots, b_k \in B$ we can now ask what the probability is that $\varphi(x_1, \ldots, x_k)$ is satisfied by $b_1, \ldots, b_k$.

When using a Bayesian network $\mfG$ for prediction as in the example we have ``lifted'' it from its original context (the set $A$) and used it on 
a new domain of objects. 
Also when moving from the fixed domain $A$ to an arbitrary domain $B$ we have, in a sense,
``lifted'' our reasoning from propositional logic to first-order logic, or some extension of it.
Perhaps this is the reason why the term ``lifted graphical model'' is used by
some authors when a graphical model is used to describe or predict (conditional) probabilities of events on an
arbitrary or unknown domain; see~\cite{KMG} for a survey of lifted graphical models.
In the subfield of machine learning, data mining and artificial intelligence called {\em statistical relational learning}
(or sometimes {\em probabilistic logic learning}) the ``lifted'' perspective is central as one here considers
general domains of objects and properties and relations that may, or may not, hold for, or between, the objects.
(See for example \cite{DKNP, GT}.)
There is no consensus regarding what, exactly, a lifted Bayesian network (let alone lifted graphical model) is or how it determines
a probability distribution on a set of ``possible worlds''.
Different approaches have been considered.
A key question is how the probability that a random variable takes a particular value is influenced by
its parents in the DAG of the Bayesian network.
The above example uses the most simple form of {\em aggregation/combination rules}.
Another approach is to use {\em aggregation/combination functions}.
Some explanation of these notions are found in e.g. \cite[p. 31, 54]{DKNP}, \cite[p. 18]{KMG}, \cite{Jae98a}.
I have not seen any definition of the notions of aggregation rule and aggregation function, but aggregation rules tend to be mainly
linguistic descriptions (think of formal logic) of how the value of a random variable depends on the values of other random variables in the network,
while aggregation functions are specifications of the dependence mainly in terms of functions.
From a practical point of view it probably makes sense to have the freedom to adapt one's lifted graphical model to 
the application at hand, so uniformity may not be a primary concern for practicians.
But to prove mathematical theorems about lifted graphical models, and the probability distributions that they induce, 
we need (of course) to 
make precise what we mean, which is done in Section~\ref{conditional probability logic}.

In this article we use aggregation rules expressed by formulas of conditional probability logic (CPL).
The idea is that for any relation symbol $R$, of arity $k$ say, 
there are an integer $\nu_R$, numbers $\alpha_{R, i} \in [0, 1]$, and 
CPL-formulas $\chi_{R, i}(x_1, \ldots, x_k)$ for $i = 1, \ldots, \nu_R$ 
such that if $\chi_{R, i}(x_1, \ldots, x_k)$ holds, then the probability that $R(x_1, \ldots, x_k)$ holds is $\alpha_{R, i}$.
This formalism is strong enough to express, for example, aggregation rules of the following kind for arbitrary $m$, any
CPL-formula $\psi(x_1, \ldots, x_k)$ and any $\alpha_i \in [0, 1]$, $i = 0, \ldots, m$:
For all $i = 0, \ldots, m$, if the proportion of $k$-tuples that satisfy $\psi(x_1, \ldots, x_k)$ is in the interval
$[i/m, (i+1)/m]$, then the probability that $R(x_1, \ldots, x_k)$ holds is $\alpha_i$.

Once we have made precise (as in Definition~\ref{definition of BN}) what we mean by a lifted Bayesian network $\mfG$
for a finite relational signature $\sigma$ (i.e. a finite set of relation symbols, possibly of different arities)
and also made precise 
(as in Definition~\ref{the probability distribution induced by the BN}) 
how $\mfG$ determines a probability distribution $\mbbP_D$ on the set of all 
$\sigma$-structures with domain $D$ (for some finite set $D$), then we can ask questions like this:
Given a CPL-formula, $\varphi(x_1, \ldots, x_k)$ and $d_1, \ldots, d_k \in D$ what is the probability that 
$\varphi(x_1, \ldots, x_k)$ is satisfied by the sequence $d_1, \ldots, d_k$? Or more formally, what is
$\mbbP_D\big(\{\mcD \in \mbW_D : \mcD \models \varphi(d_1, \ldots, d_k)\}\big)$?
It is computationally very expensive to answer the question by analyzing all members of $\mbW_D$, since, in general, the
cardinality of $\mbW_D$ is in the order of $2^{|D|^r}$ where $r$ is the maximal arity of relation symbols in $\sigma$
and $|D|$ is the cardinality of $D$.
However, our first theorem 
(Theorem~\ref{main result on quantifier elimination})
says that if $\varphi$ is ``noncritical'' in the sense that its conditional probability quantifiers (if any) avoid ``talking about'' 
certain finitely many critical numbers, then there is a {\em quantifier-free} formula $\varphi^*(x_1, \ldots, x_k)$
such that, with probability approaching 1 as $|D| \to \infty$, $\varphi$ and $\varphi^*$ are equivalent.
If we are given such $\varphi^*$ then we can easily compute the probability
$\alpha^* = \mbbP_D\big(\{\mcD \in \mbW_D : \mcD \models \varphi^*(d_1, \ldots, d_k)\}\big)$
by using only the lifted Bayesian network $\mfG$, so in particular this computation is independent of the cardinality of $D$.
Moreover, $\alpha^*$ only depends on the quantifier-free formula $\varphi^*$ and not on the choice of elements $d_1, \ldots, d_k$.
We also get that, as $|D| \to \infty$, 
$\mbbP_D\big(\{\mcD \in \mbW_D : \mcD \models \varphi(d_1, \ldots, d_k)\}\big) \to \alpha^*$.

But of course, given a noncritical $\varphi$, we have to first find a quantifier-free $\varphi^*$ which is 
``almost surely'' equivalent to $\varphi$.
The proof of 
Theorem~\ref{main result on quantifier elimination}
produces an algorithm for doing this.
At one step in the algoritm one may need to transform a quantifier-free formula into an equivalent disjunctive normal form and this
computational task is, in general, NP-hard.
But if one assumes that all quantifier-free subformulas of $\varphi$ are disjunctive normal forms,
then the algorithm that produces $\varphi^*$ works in quadratic time in the length of $\varphi$ if we 
assume that an arithmetic operation, a comparison of two numbers and a comparison of two literals is completed in one time step
(more details in Remark~\ref{remark about computational complexity}). 

The proof of Theorem~\ref{main result on quantifier elimination}
gives some by-products such as a ``logical limit/convergence law'' 
(Theorem~\ref{corollary to the main result on quantifier elimination})
and a result 
(Theorem~\ref{asymptotically equivalent bayesian network})
saying that for every lifted Bayesian network as in 
Definition~\ref{definition of BN}
there is an ``almost surely equivalent'' lifted Bayesian network in which all aggregation formulas
(as in Definition~\ref{definition of BN})
are quantifier-free.
The original zero-one law for first-order logic, independently of Glebskii et. al. \cite{Gleb} and Fagin \cite{Fag},
becomes a special case of Theorem~\ref{corollary to the main result on quantifier elimination}
when we restrict attention to first-order sentences and 
the DAG of the lifted Bayesian network has no edges and all the probabilities associated to the vertices are $1/2$.

A couple of earlier results exist which are similar to the results of this article.
Jaeger \cite{Jae98a} has considered another sort of lifted Bayesian network which he calls {\em relational Bayesian network}.
Instead of using using aggregation/combination rules (as we do in this article) relational Bayesian networks use 
aggregation/combination functions.
Theorem~3.9 in \cite{Jae98a} is as analogoue of Theorem~\ref{corollary to the main result on quantifier elimination} below
for {\em first-order} formulas in the setting of relational Bayesian networks which use only ``exponentially convergent'' combination functions.
Theorem~4.7 in \cite{Jae98b} has a similar flavour as Theorem~\ref{asymptotically equivalent bayesian network} below,
but \cite{Jae98b} considers ``admissible'' relational Bayesian networks and a probability measure defined by such on the set of structures with
a common infinite countable domain.

The results of this article are mainly motivated by concepts and methods 
in machine learning, data mining and artificial intelligence, but if the results are seen from the perspective of
finite model theory and random discrete structures, then they join a long tradition of results concerning logical limit laws
and almost sure elimination of quantifiers.
For a very small and eclectic selection of work in this field, 
ranging from the first to some of the last, see for example
\cite{Fag, Gleb, Hill, Kai, KPR, Lyn, MT, SS, Spen}.

The organization of this article is as follows.
Section~\ref{preliminaries} introduces the basic conventions used in this article as well as some basic definitions.
Section~\ref{conditional probability logic} defines the main notions of the article and states the main results.
Section~\ref{proofs} gives the proofs of these results.
The last section is a brief discussion about further research in the topics of formal logic, probabilistic graphical models, almost sure elimination 
of quantifiers and convergence laws.

\section{Preliminaries}\label{preliminaries}

\noindent
Basic knowledge of first-order logic and first-order structures is expected and there are many sources in which the 
reader can find this background, for example~\cite{Lib}.
In this section we clarify and define some basic notation and terminology concerning logic and graph theory.
Formulas of a formal logic will usually be denoted by $\varphi$, $\psi$, $\theta$ or $\chi$, possibly with sub- or superscripts.
Logical variables will be denoted $x, y, z, u, v, w$ possibly with sub- or superscripts.
Finite sequences/tuples of variables are similarly denoted $\bar{x}, \bar{y}, \bar{z}$, etc.
If a formula is denoted by $\varphi(\bar{x})$ then it is, as usual, assumed that all free variables of $\varphi$ occur in the sequence $\bar{x}$
(but we do not insist that every variable in $\bar{x}$ occurs in the formula denoted by $\varphi(\bar{x})$);
moreover in this context we will assume that all variables in $\bar{x}$ are different although this is occasionally restated.
In general, finite sequences/tuples of elements are denoted by $\bar{a}, \bar{b}, \bar{c}$, etc.
For a sequence $\bar{a}$, $\rng(\bar{a})$ denotes the set of elements occuring in $\bar{a}$.
For a sequence $\bar{a}$, $|\bar{a}|$ denotes its length.
For a set $A$, $|A|$ denotes its cardinality.
In particular, if $\varphi$ is a formula of some formal logic (so $\varphi$ is a sequence of symbols), then $|\varphi|$ denotes its length.
Sometimes we abuse notation by writing `$\bar{a} \in A$' when we actually mean that $\rng(\bar{a}) \subseteq A$.

By a {\em signature} (or {\em vocabulary}) we mean a set of relation symbols, function symbols and constant symbols.
A signature $\sigma$ is called {\em finite relational} if it is finite as a set and all symbols in it are relation symbols.
We use the terminology `{\em $\sigma$-structure}', or just {\em structure} if we omit mentioning the signature, in the sense of first-order logic.
Structures in this sense will be denoted by calligraphic letters $\mcA$, $\mcB$, $\mcC$, etc.
The domain (or universe) of a structure $\mcA$ will often be denoted by the corresponding non-calligraphic letter $A$.
A structure is called finite if its domain is finite.
If $\sigma' \subset \sigma$ are signatures and $\mcA$ is $\sigma$-structure, then $\mcA \uhrc \sigma'$ denotes 
the {\em reduct} of $\mcA$ to the signature $\sigma'$.
We let $[n]$ denote the set $\{1, \ldots, n\}$.
We use the terminology {\em atomic ($\sigma$-)formula} in the sense of first-order logic with equality, so in particular,
the expression `$x = y$' is an atomic $\sigma$-formula for every signature $\sigma$, including the empty signature $\sigma = \es$.
It will also be convenient to have a special symbol $\top$ which is viewed as an atomic $\sigma$-formula for every signature $\sigma$;
the formula $\top$ is interpreted as being true in every structure.

\begin{defin}\label{definition of atomic types}{\rm
Let $\sigma$ be a finite relational signature and $\bar{x}$ a sequence of different variables.
\begin{itemize}
\item[(i)] If $\varphi(\bar{x})$ is an atomic $\sigma$-formula, then $\varphi(\bar{x})$ and $\neg\varphi(\bar{x})$ are called {\em $\sigma$-literals}.

\item[(ii)] A consistent set of $\sigma$-literals is called an {\em atomic $\sigma$-type}.
When denoting an atomic $\sigma$-type by $p(\bar{x})$ it is assumed (as for formulas) that if a variable occurs in a formula in $p(\bar{x})$,
then it belongs to the sequence $\bar{x}$.

\item[(iii)] If $p(\bar{x})$ is an atomic $\sigma$-type, then the {\em identity fragment of $p(\bar{x})$} is the set
of formulas of the form $x_i = x_j$ or $x_i \neq x_j$ that belong to $p(\bar{x})$.

\item[(iv)] If $p(\bar{x})$ denotes an atomic $\sigma$-type and for every $\sigma$-literal $\varphi(\bar{x})$,
either $\varphi(\bar{x}) \in p(\bar{x})$ or $\neg\varphi(\bar{x}) \in p(\bar{x})$, then $p(\bar{x})$ is called a 
{\em complete atomic $\sigma$-type (with respect to $\sigma$)}.
An atomic $\sigma$-type which is not complete is sometimes called {\em partial}.

\item[(v)] Let $p(\bar{x}, \bar{y})$ be an atomic $\sigma$-type.
The {\em $\bar{y}$-dimension} of $p(\bar{x}, \bar{y})$, denoted $\dim_{\bar{y}}(p(\bar{x}, \bar{y}))$, is the maximal $d \in \mbbN$
such that there are a $\sigma$-structure $\mcA$ and $\bar{a}, \bar{b} \in A$ such that $\mcA \models p(\bar{a}, \bar{b})$
and $\big| \rng(\bar{b}) \setminus \rng(\bar{a})\big| \geq d$.

\item[(vi)] Let $\sigma' \subseteq \sigma$ and let $p$ be an atomic $\sigma$-type.
Then $p \uhrc \sigma' = \{\varphi \in p : \varphi \text{ is a $\sigma'$-formula}\}$ and 
$p \uhrc \bar{x} = \{\varphi \in p : \text{ all free variables of $\varphi$ occur in $\bar{x}$}\}$.
\end{itemize}
}\end{defin}

\begin{rem}{\rm
Note that if $p(\bar{x})$ is complete atomic $\sigma$-type where $\bar{x} = (x_1, \ldots, x_m)$, then this implies that for all $1 \leq i, j \leq m$, 
either $x_i = x_j$ or $x_i \neq x_j$ belongs to $p(\bar{x})$.
(Also observe that if $p(\bar{x}, \bar{y})$ is a complete atomic $\sigma$-type and $\dim_{\bar{y}}(p(\bar{x}, \bar{y})) = d$,
then for every $\sigma$-structure $\mcA$ and for all $\bar{a}, \bar{b}$ such that $\mcA \models p(\bar{a}, \bar{b})$, 
we have $\big|\rng(\bar{b} \setminus \rng(\bar{a})\big| = d$.
}\end{rem}

\begin{notation}\label{notation about realizations of atomic types}{\rm
Let $\sigma$ be a signature, $\bar{x}$ a sequence of different variables,
$\mcA$ a $\sigma$-structure with domain $A$ and $\bar{a} \in A^{|\bar{x}|}$.
\begin{itemize}
\item[(i)] If $p(\bar{x})$ is an atomic $\sigma$-type,
then the notation `$\mcA \models p(\bar{a})$' means that $\mcA \models \varphi(\bar{a})$ for every formula $\varphi(\bar{x}) \in p(\bar{x})$,
or in other words that $\bar{a}$ satisfies every formula in $p(\bar{x})$ with respect to the structure $\mcA$, or 
(to use model theoretic language) that $\bar{a}$ {\em realizes} $p(\bar{x})$ with respect to the structure $\mcA$.

\item[(ii)] If $\bar{y}$ is a sequence of different variables 
(such that no variable occurs in both $\bar{x}$ and $\bar{y}$) and $q(\bar{x}, \bar{y})$   is an atomic $\sigma$-type, 
then $q(\bar{a}, \mcA) = \{\bar{b} \in A^{|\bar{y}|} : \mcA \models q(\bar{a}, \bar{b})\}$.
\end{itemize}
}\end{notation}

\noindent
By a {\em directed graph} we mean a pair $(V, E)$ where $V$ is a (vertex) set and $E \subseteq V \times V$.
A {\em directed acyclic graph}, abbreviated {\em DAG}, 
is a directed graph $(V, E)$ such that $(v, v) \notin E$ for all $v \in V$ and such that 
there do not exist distinct $v_0, \ldots, v_k \in V$ for any $k \geq 1$
such that $(v_i, v_{i+1}) \in E$ for all $i = 0, \ldots, k-1$ and $(a_k, a_0) \in E$.
A {\em directed path} in a directed graph $(V, E)$ is a sequence of distinct vertices $v_0, \ldots, v_k \in V$ such that
$(v_i, v_{i+1})$ for all $i = 0, \ldots, k-1$; the {\em length} of this path is the number of edges in it, in other words, the length is $k$.

\begin{defin}\label{definitions about a DAG} {\bf (About directed acyclic graphs)} {\rm
Suppose that $\mfG$ is a DAG with nonempty and finite vertex set $V$.
Let $a \in V$. 
\begin{itemize}
\item[(i)] A vertex $b \in V$ is a {\em parent} of $a$ if $(b, a)$ is a directed edge of $\mfG$.
We let $par(a)$ denote the set of parents of $a$.

\item[(ii)] We define the {\em maximal path rank of $a$}, or just {\em mp-rank of $a$}, denoted $\mr{mpr}(a)$, to be
the length of the longest directed path having $a$ as its first vertex
(i.e. the length of the longest path $a_0, a_1, \ldots, a_k$ where $a = a_0$ and 
$(a_i, a_{i+1})$ is a directed edge for each $i = 0, \ldots, k-1$).

\item[(iii)] The {\em maximal path rank of $\mfG$}, or just {\em mp-rank of $\mfG$}, denoted $\mr{mpr}(\mfG)$
is defined as $\mr{mpr}(\mfG) = \max\{\mr{mpr}(a) : a \in V\}$.
\end{itemize}
}\end{defin}

\noindent
Observe that if $\mfG$ is a DAG with vertex set $V$ and $\mr{mpr}(\mfG) = r$ and $\mfG'$ is the induced subgraph of $\mfG$
with vertex set $V' = \{a \in V : \mr{mpr}(a) < r\}$, then, for every $a \in V'$, the mp-rank of $a$ is the same no matter if we compute
it with respect to $\mfG'$ or with respect to $\mfG$; it follows that
$\mr{mpr}(\mfG') = r-1$.

We call a random variable {\em binary} if it can only take the value $0$ or $1$.
The following is a direct consequence of \cite[Corollary~A.1.14]{AlonSpencer} which in turn follows from the
Chernoff bound \cite{Chernoff}:

\begin{lem}\label{independent bernoulli trials}
Let $Z$ be the sum om $n$ independent binary random variables, each one with probability $p$ of having the value 1.
For every $\varepsilon > 0$ there is $c_\varepsilon > 0$, depending only on $\varepsilon$, such that the probability that
$|Z - pn| > \varepsilon p n$ is less than $2 e^{-c_\varepsilon p n}$.
\end{lem}

\section{Conditional probability logic and lifted Bayesian networks}\label{conditional probability logic}

\noindent
In this section we define the main concepts of this article and state the main results.

\begin{defin}\label{definition of CPL} {\bf (Conditional probability logic)} {\rm
Suppose that $\sigma$ is a signature.
Then the set of {\em conditional probability formulas over $\sigma$}, denoted $CPL(\sigma)$, is defined inductively as follows:
\begin{enumerate}
\item Every atomic $\sigma$-formula belongs to $CPL(\sigma)$ (where `atomic' has the same meaning as in first-order logic with equality).

\item If $\varphi, \psi \in CPL(\sigma)$ then $(\neg\varphi), (\varphi\wedge\psi), (\varphi\vee\psi), (\varphi\rightarrow\psi), 
(\varphi\leftrightarrow\psi), (\exists x \varphi) \in CPL(\sigma)$ where $x$ is a variable.
(As usual, in practice we do not necessarily write out all parentheses.)
We consider $\forall x \varphi$ to be an abbreviation of $\neg \exists x \neg \varphi$.

\item If $r \geq 0$ is a real number, $\varphi, \psi, \theta, \tau \in CPL(\sigma)$ and $\bar{y}$ is a sequence of distinct variables, then
\begin{align*}
&\Big( r + \| \varphi \ |  \ \psi \|_{\bar{y}} \ \geq \ 
\| \theta \ |  \ \tau \|_{\bar{y}} \Big) \in CPL(\sigma)
\ \ \text{ and} \\
&\Big( \| \varphi \ |  \ \psi \|_{\bar{y}} \ \geq \ 
\| \theta \ |  \ \tau \|_{\bar{y}} + r \Big) \in CPL(\sigma).
\end{align*}
In both these new formulas all variables of $\varphi, \psi, \theta$ and $\tau$ that appear in the sequence $\bar{y}$ become {\em bound}.
So this construction can be seen as a sort of quantification, which may become more clear by the provided semantics below.
\end{enumerate}
}\end{defin}

\noindent
A formula $\varphi \in CPL(\sigma)$ is called {\em quantifier-free} if it contains no quantifier, that is, if it is constructed from atomic formulas by
using only connectives $\neg, \wedge, \vee, \rightarrow, \leftrightarrow$.

\begin{defin}\label{semantics of CPL} {\bf (Semantics)} {\rm
\begin{enumerate}
\item The interpretations of $\neg, \wedge, \vee, \rightarrow, \leftrightarrow$ and $\exists$ are as in first-order logic.

\item Suppose that $\mcA$ is a {\em finite} $\sigma$-structure and let $\varphi(\bar{x}, \bar{y}), \psi(\bar{x}, \bar{y}), 
\theta(\bar{x}, \bar{y}), \tau(\bar{x}, \bar{y}) \in CPL(\sigma)$.
Let $\bar{a} \in A^{|\bar{x}|}$.
\begin{enumerate}
\item We define $\varphi(\bar{a}, \mcA) = \big\{\bar{b} \in A^{|\bar{y}|} : \mcA \models \varphi(\bar{a}, \bar{b}) \big\}$.

\item The expression 
\[
\mcA \  \models \ 
\Big( r + \| \varphi(\bar{a}, \bar{y}) \ | \ \psi(\bar{a}, \bar{y}) \|_{\bar{y}} \ \geq \ 
\| \theta(\bar{a}, \bar{y}) \ | \ \tau(\bar{a}, \bar{y}) \|_{\bar{y}} \Big)
\]
means that $\psi(\bar{a}, \mcA) \neq \es$, $\tau(\bar{a}, \mcA) \neq \es$ and
\[
r + \frac{\big| \varphi(\bar{a}, \mcA) \cap \psi(\bar{a}, \mcA) \big|}{\big| \psi(\bar{a}, \mcA) \big|} \ \geq \ 
\frac{\big| \theta(\bar{a}, \mcA) \cap \tau(\bar{a}, \mcA) \big|}{\big| \tau(\bar{a}, \mcA) \big|}
\]
and in this case we say that 
$\Big( r + \| \varphi(\bar{a}, \bar{y}) \ | \ \psi(\bar{a}, \bar{y}) \|_{\bar{y}} \ \geq \ 
\| \theta(\bar{a}, \bar{y}) \ | \ \tau(\bar{a}, \bar{y}) \|_{\bar{y}} \Big)$
is true (or holds) in $\mcA$.
If $\psi(\bar{a}, \mcA) = \es$ or $\tau(\bar{a}, \mcA) = \es$ or
\[
r + \frac{\big| \varphi(\bar{a}, \mcA) \cap \psi(\bar{a}, \mcA) \big|}{\big| \psi(\bar{a}, \mcA) \big|} \ < \ 
\frac{\big| \theta(\bar{a}, \mcA) \cap \tau(\bar{a}, \mcA) \big|}{\big| \tau(\bar{a}, \mcA) \big|}
\]
then we write
\[
\mcA \  \not\models \ 
\Big( r + \| \varphi(\bar{a}, \bar{y}) \ | \ \psi(\bar{a}, \bar{y}) \|_{\bar{y}} \ \geq \ 
 \| \theta(\bar{a}, \bar{y}) \ | \ \tau(\bar{a}, \bar{y}) \|_{\bar{y}} \Big)
\]
and say that 
$\Big( r + \| \varphi(\bar{a}, \bar{y}) \ | \ \psi(\bar{a}, \bar{y}) \|_{\bar{y}} \ \geq \ 
\| \theta(\bar{a}, \bar{y}) \ | \ \tau(\bar{a}, \bar{y}) \|_{\bar{y}} \Big)$
is false in $\mcA$.

\item The meaning of 
\[
\mcA \  \models \ 
\Big( \| \varphi(\bar{a}, \bar{y}) \ | \ \psi(\bar{a}, \bar{y}) \|_{\bar{y}} \ \geq \ 
\| \theta(\bar{a}, \bar{y}) \ | \ \tau(\bar{a}, \bar{y}) \|_{\bar{y}} + r \Big)
\]
is defined similarly.
\end{enumerate}
\end{enumerate}
}\end{defin}

\begin{rem}\label{warning about semantics} {\bf (A warning)} {\rm
Observe that with the given semantics,
\[
\mcA \  \not\models \ 
\Big( r + \| \varphi(\bar{a}, \bar{y}) \ | \ \psi(\bar{a}, \bar{y}) \|_{\bar{y}} \ \geq \ 
 \| \theta(\bar{a}, \bar{y}) \ | \ \tau(\bar{a}, \bar{y}) \|_{\bar{y}} \Big)
\]
does {\em not} necessarily imply
\[
\mcA \  \models \ 
\Big( r + \| \varphi(\bar{a}, \bar{y}) \ | \ \psi(\bar{a}, \bar{y}) \|_{\bar{y}} \ \leq \ 
 \| \theta(\bar{a}, \bar{y}) \ | \ \tau(\bar{a}, \bar{y}) \|_{\bar{y}} \Big)
\]
because the first formula may fail to be true for $\bar{a}$ because $\psi(\bar{a}, \mcA) = \es$ or $\tau(\bar{a}, \mcA) = \es$
in which case the corresponding fraction is undefined and then also the other formula is false for $\bar{a}$.
}\end{rem}

\begin{rem}\label{remark on abbreviations} {\bf (Expressing conditional probabilities, or just probabilities)} {\rm
Let $\bar{x} = (x_1, \ldots, x_k)$ and $\bar{y} = (y_1, \ldots, y_l)$.
If $\tau(\bar{x}, \bar{y})$ denotes the formula $y_1 = y_1$ and $\theta(\bar{x}, \bar{y})$ denotes the formula $y_1 \neq y_1$,
then 
\begin{equation}\label{formula expressing conditional probability}
\Big( \| \varphi(\bar{x}, \bar{y}) \ | \ \psi(\bar{x}, \bar{y}) \|_{\bar{y}} \ \geq \ 
\| \theta(\bar{x}, \bar{y}) \ | \ \tau(\bar{x}, \bar{y}) \|_{\bar{y}} + r \Big)
\end{equation}
expresses that the proportion of tuples $\bar{y}$ that satisfy $\varphi(\bar{x}, \bar{y})$ among those $\bar{y}$ that
satisfy $\psi(\bar{x}, \bar{y})$ is at least $r$.
Thus the formula expresses a conditional probability if we assume that all $l$-tuples have the same probability.
Under the stated assumptions, let us abbreviate~(\ref{formula expressing conditional probability}) by
\begin{equation}\label{abbreviation of formula expressing conditional probability}
\Big(  \| \varphi(\bar{x}, \bar{y}) \ | \ \psi(\bar{x}, \bar{y}) \|_{\bar{y}} \ \geq \ r \Big).
\end{equation}
If we assume, in addition, that $\psi(\bar{x}, \bar{y})$ is the formula $y_1 = y_1$, then each 
of~(\ref{formula expressing conditional probability}) and~(\ref{abbreviation of formula expressing conditional probability}) 
expresses
that the proportion of $l$-tuples $\bar{y}$ that satisfy $\varphi(\bar{x}, \bar{y})$ is at least $r$.
}\end{rem}

\begin{exam}\label{example 1}{\rm
Suppose that $M$ is a unary relation symbol and $F$ a binary relation symbol. 
Consider the statement ``For at least half of all persons $x$, if at least one third of the friends of $x$ are mathematicians,
then $x$ is a mathematician''. 
If $M(x)$ expresses that ``$x$ is a mathematician'' and $F(x, y)$ expresses that ``$x$ and $y$ are friends'',
then this statement can be formulated in CPL, using the abbreviation~(\ref{abbreviation of formula expressing conditional probability}), as
\[
\Big( \big\| \big( \| M(y) \ | \ F(x, y) \|_y \geq 1/3 \big) \ \rightarrow \ M(x) \ \big| \ x = x \big\|_x \geq 1/2 \Big).
\]
}\end{exam}

\begin{rem}\label{remark on expressing independence} {\bf (Expressing independence)} {\rm
Suppose that $\mcA$ is a finite $\sigma$-structure, $\theta(\bar{x}, \bar{y})$ is the formula $y_1 = y_1$ and $\bar{a} \in A^{|\bar{x}|}$.
If $r = 0$ and
\begin{align*}
\mcA \models &\Big( r + \| \varphi(\bar{a}, \bar{y}) \ | \ \psi(\bar{a}, \bar{y}) \|_{\bar{y}} \ \geq \ 
\| \varphi(\bar{a}, \bar{y}) \ | \ \theta(\bar{a}, \bar{y}) \|_{\bar{y}} \Big)
\ \wedge \\ 
&\Big( \| \varphi(\bar{a}, \bar{y}) \ | \ \theta(\bar{a}, \bar{y}) \|_{\bar{y}} \ \geq \ 
\| \varphi(\bar{a}, \bar{y}) \ | \ \psi(\bar{a}, \bar{y}) \|_{\bar{y}} + r\Big),
\end{align*}
then the event $\mbX = \{\bar{b} \in A^{|\bar{y}|} : \mcA \models \varphi(\bar{a}, \bar{b})\}$ is independent from 
the event $\mbY = \{\bar{b} \in A^{|\bar{y}|} : \mcA \models \psi(\bar{a}, \bar{b})\}$ if all $|\bar{y}|$-tuples have the same probability.

If $\mcA$ represents a database from the real world, then it is unlikely that events of interest are (conditionally) independent according
the precise mathematical definition. Instead one may look for ``approximate (conditional) independencies''.
If $r$ is changed to be a small positive number and if
\begin{align*}
\mcA \models &\Big( r + \| \varphi(\bar{a}, \bar{y}) \ | \ \psi(\bar{a}, \bar{y}) \|_{\bar{y}} \ \geq \ 
\| \varphi(\bar{a}, \bar{y}) \ | \ \theta(\bar{a}, \bar{y}) \|_{\bar{y}} \Big)
\ \wedge \\ 
&\Big( r + \| \varphi(\bar{a}, \bar{y}) \ | \ \theta(\bar{a}, \bar{y}) \|_{\bar{y}} \ \geq \ 
\| \varphi(\bar{a}, \bar{y}) \ | \ \psi(\bar{a}, \bar{y}) \|_{\bar{y}} \Big) \ \wedge \\
&\Big( r + \| \psi(\bar{a}, \bar{y}) \ | \ \varphi(\bar{a}, \bar{y}) \|_{\bar{y}} \ \geq \ 
\| \psi(\bar{a}, \bar{y}) \ | \ \theta(\bar{a}, \bar{y}) \|_{\bar{y}} \Big)
\ \wedge \\ 
&\Big( r + \| \psi(\bar{a}, \bar{y}) \ | \ \theta(\bar{a}, \bar{y}) \|_{\bar{y}} \ \geq \ 
\| \psi(\bar{a}, \bar{y}) \ | \ \varphi(\bar{a}, \bar{y}) \|_{\bar{y}} \Big)
\end{align*}
then the dependency between $\mbX$ and $\mbY$ is weak,
or one could say that they are ``approximately independent up to an error of $r$''.
The reason for the more complicated formula is to make ``$r$-approximate independence'' symmetric.
}\end{rem}

\begin{defin}\label{definition of quantifier rank}{\rm
The {\em quantifier rank}, $\mr{qr}(\varphi)$, of formulas $\varphi \in CPL(\sigma)$ is defined inductively as follows:
\begin{enumerate}
\item For atomic $\varphi$, $\mr{qr}(\varphi) = 0$.

\item $\mr{qr}(\neg\varphi) = \mr{qr}(\varphi)$, $\mr{qr}(\varphi \star \psi) = \max\{\mr{qr}(\varphi), \mr{qr}(\psi)\}$
if $\star$ is one of $\wedge, \vee, \rightarrow$ or $\leftrightarrow$.

\item $\mr{qr}(\exists x \varphi) = \mr{qr}(\varphi) + 1$

\item
$\mr{qr}\Big( \big( r + \| \varphi \ |  \ \psi \|_{\bar{y}} \ \geq \ \| \theta \ |  \ \tau \|_{\bar{y}} \big) \Big) \ = \ 
\mr{qr}\Big( \big( \| \varphi \ |  \ \psi \|_{\bar{y}} \ \geq \ \| \theta \ |  \ \tau \|_{\bar{y}} + r \big) \Big) \ = $ \\
$\max\{\mr{qr}(\varphi), \mr{qr}(\psi), \mr{qr}(\theta), \mr{qr}(\tau)\} + |\bar{y}|.$
\end{enumerate}
}\end{defin}

\begin{defin}\label{definition of BN} {\bf (Lifted Bayesian network)} {\rm
Let $\sigma$ be a finite relational signature.
In this article we define a {\em lifted Bayesian network for $\sigma$} to consist of the following components:
\begin{itemize}
\item[(a)] An acyclic directed graph (DAG) $\mfG$ with vertex set $\sigma$.

\item[(b)] For each $R \in \sigma$, a number $\nu_R \in \mbbN^+$, formulas $\chi_{R, i}(\bar{x}) \in CPL(\mr{par}(R))$,
for $i = 1, \ldots, \nu_R$, where $|\bar{x}|$ equals the arity of $R$, such that
$\forall \bar{x} \big( \bigvee_{i = 1}^{\nu_R} \chi_{R, i}(\bar{x})\big)$ is valid (i.e. true in all $\mr{par}(R)$-structures) and if
$i \neq j$ then $\exists \bar{x} \big(\chi_{R, i}(\bar{x}) \wedge \chi_{R, j}(\bar{x})\big)$
is unsatisfiable.
Each $\chi_{R, i}$ will be called an {\em aggregation formula (of $\mfG$)}.

\item[(c)] For each $R \in \sigma$ and each $1 \leq i \leq \nu_R$, 
a number denoted $\mu(R \ | \ \chi_{R, i})$ (or $\mu(R(\bar{x}) \ | \ \chi_{R, i}(\bar{x}))$)
in the interval $[0, 1]$.
\end{itemize}
}\end{defin}

\noindent
We will use the same symbol (for example $\mfG$) to denote a lifted Bayesian network and its underlying DAG.
The intuitive meaning of $\mu(R \ | \ \chi_{R, i})$ in part~(c) is that if $\bar{a}$ is a sequence of elements from a structure and
$\bar{a}$ satisfies $\chi_{R, i}(\bar{x})$, then the probability that $\bar{a}$ satisfies $R(\bar{x})$ is $\mu(R \ | \ \chi_{R, i})$.

\begin{rem}\label{remark on subnetworks}{\bf (Subnetworks)} {\rm
Let $\mfG$ denote a lifted Bayesian network for $\sigma$.
Suppose that $\sigma' \subset \sigma$ is such that if $R \in \sigma'$ then $\mr{par}(R) \subseteq \sigma'$.
Then it is easy to see that $\sigma'$ determines a lifted Bayesian network $\mfG'$ for $\sigma'$ such that 
\begin{itemize}
\item the vertex set of the underlying DAG of $\mfG'$ is $\sigma'$,

\item for every $R \in \sigma'$, the number $\nu_R$ and the formulas $\chi_{R, i}$, $i = 1, \ldots, \nu_R$, are the same as those for $\mfG$,

\item for every $R \in \sigma'$ and every $1 \leq i \leq \nu_R$, the numbers $\mu(R \ | \ \chi_{R, i})$ are the same as those for $\mfG$.
\end{itemize}
We call the so defined lifted Bayesian network $\mfG'$ for $\sigma'$ the {\em subnetwork (of $\mfG$) induced by $\sigma'$}.
}\end{rem}

\begin{defin}\label{empty BN}{\bf (The case of an empty signature)} {\rm
(i) As a technical convenience we will also consider a lifted Bayesian network, denoted $\mfG^\es$, for the empty signature $\es$.
According to Definition~\ref{definition of BN} the vertex set of the underlying DAG of $\mfG^\es$ is $\es$, the empty set.
It follows that no formulas or numbers as in parts~(b) and~(c) of Definition~\ref{definition of BN} need to be specified for $\mfG^\es$.\\
(ii) For every $n \in \mbbN^+$, let $\mbW^\es_n$ denote the set of all $\es$-structures with domain $[n]$ and note that every 
$\mbW^\es_n$ has only one member which is just the set $[n]$.\\
(iii) For every $n \in \mbbN^+$, let $\mbbP^\es_n$ be the unique probability distribution on $\mbW^\es_n$.
}\end{defin}

\begin{defin}\label{the probability distribution induced by the BN}{\bf (The probability distribution in the general case)} {\rm
Let $\sigma$ be a finite nonempty relational signature and let $\mfG$ denote a lifted Bayesian network for $\sigma$.
Suppose that the underlying DAG of $\mfG$ has mp-rank $\rho$.
For each $0 \leq r \leq \rho$ let $\mfG_r$ be the subnetwork (in the sense of Remark~\ref{remark on subnetworks})
induced by $\sigma_r = \{R \in \sigma : \mr{mpr}(R) \leq r\}$ and note that $\mfG_\rho = \mfG$.
Also let $\sigma_{-1} = \es$, $\mfG_{-1} = \mfG^\es$ and let $\mbbP^{-1}_n$ be the unique probability distribution on $\mbW^{-1}_n = \mbW^\es_n$.
By induction on $r$ we define, for every $r = 0, 1, \ldots, \rho$, a probability distribution $\mbbP^r_n$ on the set $\mbW^r_n$ of all
$\sigma_r$-structures with domain $[n]$ as follows:
For every $\mcA \in \mbW^r_n$,
\[
\mbbP^r_n(\mcA) \ = \ \mbbP^{r-1}_n(\mcA \uhrc \sigma_{r-1}) 
\prod_{R \in \sigma_r \setminus \sigma_{r-1}} \ \prod_{i=1}^{\nu_R} \ \prod_{\bar{a} \in \chi_{R, i}(\mcA \uhr \sigma_{r-1})} 
\lambda(\mcA, R, i, \bar{a})
\]
where
\[
\lambda(\mcA, R, i, \bar{a}) = 
\begin{cases}
\mu(R \ | \ \chi_{R, i}) \ \ \ \ \ \ \ \text{ if } \mcA \models \chi_{R, i}(\bar{a}) \wedge R(\bar{a}),\\
1 - \mu(R \ | \ \chi_{R, i}) \ \ \text{ if } \mcA \models \chi_{R, i}(\bar{a}) \wedge \neg R(\bar{a}),\\
\text{0 \ \ \ \ \ \ \ \ \ \ \ \ \ \ \ \ \ \ \ \ otherwise.}
\end{cases}
\]
Finally we let $\mbW_n = \mbW^\rho_n$ and $\mbbP_n = \mbbP^\rho_n$, so $\mbbP_n$ is a probability distribution on the set
of all $\sigma$-structures with domain $[n]$.
}\end{defin}

\begin{rem}\label{remark about reflexive and/or symmetric relations} {\bf ((Ir)reflexive and/or symmetric relations)} {\rm 
Let $A$ be a set and let $R \subseteq A^k$ be a $k$-ary relation on $A$.
We call $R$ {\em reflexive} if for all $a \in A$ the $k$-tuple containing $a$ in each coordinate belongs to $R$.
We call $R$ {\em irreflexive} if for every $(a_1, \ldots, a_k) \in R$ we have $a_i \neq a_j$ if $i \neq j$.
We call $R$ {\em symmetric} if for every $(a_1, \ldots, a_k) \in R$, every permutation of $(a_1, \ldots, a_k)$ also belongs to $R$.
Consider Definition~\ref{the probability distribution induced by the BN}
and let $R \in \sigma$.
We can make sure that $\mbbP_n(\mcA) > 0$ only if the interpretation of $R$ in $\mcA$ is
reflexive (respectively irreflexive) by choosing the formulas $\chi_{R, i}$ and associated (conditional) probabilities in an appropriate way.
To achieve that $\mbbP_n(\mcA) > 0$ only if the interpretation of $R$ in $\mcA$ is symmetric we can do like this:
In the definition of $\lambda(\mcA, R, i, \bar{a})$ (in Definition~\ref{the probability distribution induced by the BN})
we interpret $R(\bar{a})$ as meaning that $R$ is satisfied by every permutation of $\bar{a}$ and we interpret $\neg R(\bar{a})$ as
meaning that $R$ is not satisfied by any permutation of $\bar{a}$.
We also need to assume that for every $k$-tuple $\bar{a}$, either every permutation of $\bar{a}$ satisfies $\chi_{R, i}(\bar{x})$
or no permutation of $\bar{a}$ satisfies $\chi_{R, i}(\bar{x})$.
Then the proof of 
Theorems~\ref{main result on quantifier elimination} --~\ref{asymptotically equivalent bayesian network}
still works out with very small modifications.
}\end{rem}

\begin{defin}\label{definition of probability of a sentence}{\rm
Let $\sigma$, $\mbW_n$ and $\mbbP_n$ be as in Definition~\ref{the probability distribution induced by the BN}.\\
(i) If $\varphi(\bar{x}) \in CPL(\sigma)$ and $\bar{a} \in [n]^{|\bar{x}|}$, then we define
$\mbbP_n(\varphi(\bar{a})) = \mbbP_n\big(\{\mcA \in \mbW_n : \mcA \models \varphi(\bar{a})\}\big)$.\\
(ii) If $\varphi \in CPL(\sigma)$ has no free variables (i.e. is a sentence),
then we define $\mbbP_n(\varphi) = \mbbP_n\big(\{\mcA \in \mbW_n : \mcA \models \varphi\}\big)$.
}\end{defin}

\noindent
Now we can state the main results.
They use the notion of {\em noncritical formula} which depends on the lifted Bayesian network under consideration.
Since this notion is quite technical and relies on some technical results (concerning the convergence of the probability that
an atomic type is realized) which will be proved later, we give the precise definition later in 
Definition~\ref{definition of noncritical formula}; in that context it will be more evident why the definition of noncritical formula 
looks as it looks.
For now I only say this: For every $m \in \mbbN^+$ there are finitely many numbers (depending only on $\mfG$) which are called
{\em $m$-critical} (according to Definition~\ref{definition of critical number}).
Roughly speaking, a formula $\varphi(\bar{x}) \in CPL(\sigma)$ is {\em noncritical} 
(details in Definition~\ref{definition of noncritical formula}) 
if for every subformula (of $\varphi(\bar{x})$)
of the form 
\begin{align*}
\Big( r + \| \chi \ |  \ \psi \|_{\bar{y}} \ \geq \ 
\| \theta \ |  \ \tau \|_{\bar{y}} \Big) \
\ \ \text{or} \ \ 
\Big( \| \chi \ |  \ \psi \|_{\bar{y}} \ \geq \ 
\| \theta \ |  \ \tau \|_{\bar{y}} + r \Big) 
\end{align*}
the number $r$ is {\em not} the difference of two $m$-critical numbers where $m = |\bar{x}| + \mr{qr}(\varphi)$.
It follows that {\em every first-order formula is noncritical}.
For a longer discussion on the topics of critical formulas and (non)convergence see 
Remark~\ref{remark about the necessity of noncriticality}.

\begin{theor}\label{main result on quantifier elimination}{\bf (Almost sure elimination of quantifiers for noncritical formulas)}
Let $\sigma$ be a finite relational signature, let $\mfG$ be a lifted Bayesian network and, for each $n \in \mbbN^+$,
let $\mbbP_n$ be the probability distribution induced by $\mfG$
(according to Definition~\ref{the probability distribution induced by the BN}) 
on the set $\mbW_n$ of all $\sigma$-structures with domain $[n]$.
Suppose that every aggregation formula $\chi_{R, i}$ of $\mfG$ is noncritical.
If $\varphi(\bar{x}) \in CPL(\sigma)$ is noncritical, 
then there are a quantifier free formula $\varphi^*(\bar{x}) \in CPL(\sigma)$ and $c > 0$, which depend only on
$\varphi(\bar{x})$ and $\mfG$, such that for all sufficiently large~$n$
\[
\mbbP_n\big(\forall \bar{x} (\varphi(\bar{x}) \leftrightarrow \varphi^*(\bar{x}))\big) \ \geq \ 1 - e^{-cn}.
\]
\end{theor}

\begin{theor}\label{corollary to the main result on quantifier elimination} {\bf (Convergence for noncritical formulas)}
Let $\sigma$, $\mfG$, $\mbW_n$ and $\mbbP_n$ be as in Theorem~\ref{main result on quantifier elimination}.
For every noncritical $\varphi(\bar{x}) \in CPL(\sigma)$ 
there are $c > 0$ and $0 \leq d \leq 1$, depending only on $\varphi(\bar{x})$ and $\mfG$,
such that for every $m \in \mbbN^+$ and every $\bar{a} \in [m]^{|\bar{x}|}$,
\begin{align*}
\big| \mbbP_n(\varphi(\bar{a})) - d \big| \ \leq \ e^{-cn} \quad \text{for all sufficiently large $n \geq m$.}
\end{align*}
The number $d$ is always critical (i.e. $l$-critical for some $l$).
Moreover, if $\varphi$ has no free variable (i.e. is a sentence), then $\mbbP_n(\varphi)$ converges to either~0 or~1.
\end{theor}

\begin{theor}\label{asymptotically equivalent bayesian network} {\bf (An asymptotically equivalent ``quantifier-free'' network)}
Let $\sigma$, $\mfG$, $\mbW_n$ and $\mbbP_n$ be as in Theorem~\ref{main result on quantifier elimination}.
Then for every aggregation formula $\chi_{R, i}(\bar{x})$ of $\mfG$ there is a quantifier-free formula 
$\chi_{R, i}^*(\bar{x})$ containing only relation symbols that occur in $\chi_{R, i}$ such that if $\mfG^*$ is the
lifted Bayesian network 
\begin{itemize}
\item with the same underlying DAG as $\mfG$,

\item where, for every $R \in \sigma$ and every $1 \leq i \leq \nu_R$, the aggregation formula $\chi_{R, i}$ is replaced by $\chi_{R, i}^*$, and

\item where $\mu^*(R \ | \ \chi_{R, i}^*) = \mu(R \ | \ \chi_{R, i})$ for every $R \in \sigma$ and every $1 \leq i \leq \nu_R$,
\end{itemize}
then for every noncritical $\varphi(\bar{y}) \in CPL(\sigma)$
there is $d > 0$, depending only on $\varphi(\bar{y})$ and $\mfG$, 
such that for every $m \in \mbbN^+$ and every $\bar{a} \in [m]^{|\bar{y}|}$, 
\[
\big| \mbbP_n(\varphi(\bar{a})) - \mbbP_n^*(\varphi(\bar{a})) \big| \ \leq \ e^{-dn}, 
\ \text{ for all sufficiently large } n \geq m,
\]
where $\mbbP_n^*$ is the the probability distribution on $\mbW_n$ according to 
Definition~\ref{the probability distribution induced by the BN}
if $\mfG$ is replaced by $\mfG^*$ and $\mbbP_n$ is replaced by $\mbbP_n^*$.
\end{theor}

\begin{rem}\label{remark about computational complexity} {\bf (Computational complexity)} {\rm 
The proof of Theorem~\ref{main result on quantifier elimination} indicates an algorithm for 
finding the quantifier-free $\varphi^*$ from $\varphi$.
Suppose that we fix the lifted Bayesian network (so $\sigma$ is also fixed)
and try to understand how efficient the algorithm is with respect to the length of $\varphi$.
The crucial step is Definition~\ref{definition of I}
and Lemma~\ref{the case when I is nonempty}
which together show how to eliminate a quantifier of the form
constructed in part~(3) of Definition~\ref{definition of CPL} in a satisfiable formula.
However, at this step in the proof we assume that the formulas inside the latest quantification are
written as disjunctions of complete atomic types.
The problem of transforming an arbitrary quantifier-free formula into an equivalent disjunctive normal form is NP-hard
so the algorithm is not necessarily efficient in general (given the current state of affairs in computational complexity theory).
But if we assume that every quantifier-free subformula of $\varphi$ is a disjunctive normal form,
then the number ``steps'' that the indicated algorithm needs to find $\varphi^*$ is $O(|\varphi|^2)$ if $|\varphi|$ 
denotes the length of $\varphi$ and  ``step'' means an arithmetic operation\footnote{
More precisely, adding, multiplying or dividing two numbers.}, a comparison of two numbers or a
comparison of two literals.
This essentially follows from Remark~\ref{remark on computing I}
because the number of times that a quantifier needs to be eliminated is bounded by $|\varphi|$.
}\end{rem}

\begin{rem}\label{remark about the necessity of noncriticality} {\bf (Necessity of noncriticality)} {\rm
It follows from Remark~\ref{remark on abbreviations} that for every sentence $\psi$ of the language $\mcL_{\omega P}$
considered in~\cite{KL} there is a sentence of CPL which has exactly the same finite models as $\psi$.
Therefore it follows from \cite[Proposition~3.1]{KL} that the assumption that $\varphi$ is noncritical in
Theorems~\ref{main result on quantifier elimination}
and~\ref{corollary to the main result on quantifier elimination}
is necessary. More precisely, let $\sigma$ contain one binary relation symbol and no other symbols and let
$\mfG$ be a lifted Bayesian network for $\sigma$ where $\mu(R(x, y) \ | \ x = x) = 1/2$ and `$x = x$' is the
only aggregation formula associated to $R$. Then, according to \cite[Proposition~3.1]{KL} interpreted in 
the present context, there is a
(``critical'') sentence $\psi \in CPL(\sigma)$ such that $\mbbP_n(\psi)$ does not converge.

We now generalize the idea of \cite[Proposition~3.1]{KL} to show that nonconvergence for
at least some ``critical'' formulas is the case for many (if not all) lifted Bayesian networks. 
Let $\sigma$ be a finite relational signature, let $\mfG$ be a lifted Bayesian network for $\sigma$ such that
every aggregation formula of $\mfG$ is noncritical. 
Let $\varphi(\bar{x}) \in CPL(\sigma)$ be a noncritical 
formula where $\bar{x} = (x_1, \ldots, x_k)$ and $k \geq 2$.
By Theorem~\ref{corollary to the main result on quantifier elimination}
there is $0 \leq d \leq 1$ such that for every $n_0 \in \mbbN^+$ and every $\bar{a} \in [n_0]^k$, 
$\lim_{n\to\infty} \mbbP_n(\varphi(\bar{a})) = d$.
By the same theorem, $d$ is a ``critical'' number, that is,  $l$-critical for some $l$ according to
Definition~\ref{definition of critical number}.
Suppose that $0 < d < 1$ (which would typically be the case if $\varphi(\bar{x})$ is atomic).
Furthermore, suppose that all numbers of the form $\mu(R \ | \ \chi_{R, i})$ associated with $\mfG$
(as in Definition~\ref{definition of BN}) are rational.
It then follows (from Definition~\ref{definition of critical number} and results and definitions preceeding it) that $d$ is rational.
Now suppose that for all $n$ and any choice of distinct $\bar{a}_1, \ldots, \bar{a}_m \in [n]$,
the binary random variables $X_1, \ldots, X_m$ with domain $\mbW_n$ are independent, 
where $X_i(\mcA) = 1$ if $\mcA \models \varphi(\bar{a}_i)$ and $X_i(\mcA) = 0$ otherwise.
From item~(b) of Remark~\ref{remark on P restricted to A'} one can derive that if $\varphi(\bar{x})$ is atomic 
then this independence assumption holds.

Now we consider the following formula, denoted $\psi(x_1, \ldots, x_{k-1})$:
\begin{align*}
\big(\| \varphi(x_1, \ldots, x_{k-1}, y) \ | \ y = y \|_y \ \geq \ d \big) \ \wedge \
\big(\| \neg\varphi(x_1, \ldots, x_{k-1}, y) \ | \ y = y \|_y \ \geq \ 1 - d \big).
\end{align*}
In structures in $\mbW_n$ it expresses that ``there are exactly $dn$ elements $y$ such that 
$\varphi(\bar{x}', y)$ is satisfied'' where  $\bar{x}' = (x_1, \ldots, x_{k-1})$.
We will show that $\mbbP_n(\exists \bar{x}' \psi(\bar{x}'))$ does not converge as $n \to \infty$.

If $dn$ is not an integer and $\mcA \in \mbW_n$, then $\mcA \not\models \exists\bar{x}' \varphi(\bar{x}')$,
so $\mbbP_n(\exists \bar{x}' \psi(\bar{x}')) = 0$.
Note that as $d$ is rational there are infinitely many $n$ such that $dn$ is an integer and infinitely many
$n$ such that $dn$ is not an integer.
Hence it suffices to show that $\mbbP_n(\exists \bar{x}' \psi(\bar{x}'))$ gets arbitrarily close to 1
for sufficiently large $n$ such that $dn$ is an integer.
Fix a large $n$ such that $dn$ is an integer and $\bar{a} \in [n]^{k-1}$.
Then, for every $b \in [n]$, $\mbbP_n(\varphi(\bar{a}, b))$ is very close to $d$.
By the assumption about independence above, it follows that the probability that ``there are exactly
$dn$ elements $b$ such that $\varphi(\bar{a}, b)$ holds'' is close to
$\binom{n}{dn} dn (1-d)^{(1-d)n}$.
Thus the probability of the negation of this statement is close to
$1 - \binom{n}{dn} dn (1-d)^{(1-d)n}$ and, using the assumption about independence again,
it follows that $\mbbP_n(\neg\exists\bar{x}' \varphi(\bar{x}'))$ is close to
\begin{align*}
&\bigg( 1 - \binom{n}{dn} dn (1-d)^{(1-d)n}\bigg)^{n^{k-1}} \ \leq \ 
\bigg( 1 - \binom{n}{dn} dn (1-d)^{(1-d)n}\bigg)^n \\
\sim \ &\bigg( 1 - \frac{1}{\sqrt{2\pi d(1-d)n}} \bigg)^n \ = \ \bigg( 1 - \frac{c}{\sqrt{n}}\bigg)^n
\ (\text{by Stirlings approximation}) \\
\sim \ &e^{-c\sqrt{n}} \to 0 \  \ \text{ as } n \to \infty \ (\text{where } c = 1/\sqrt{2\pi d(1-d)}).
\end{align*}

Finally let us take a broader look at a formula of the form
\begin{equation}\label{a formula exemplifying necessity of noncriticality}
r \ + \ \| \varphi_1(\bar{y}) \ | \ \varphi_2(\bar{y}) \|_{\bar{y}} \ \geq \ 
\| \varphi_3(\bar{y}) \ | \ \varphi_4(\bar{y}) \|_{\bar{y}}.
\end{equation}
Suppose that $\varphi_i$ is noncritical for every $i$.
By Theorem~\ref{corollary to the main result on quantifier elimination},
there are numbers $d_1, d_2, d_3, d_4$ such that for any $n_0$, every $\bar{a} \in [n_0]^{|\bar{y}|}$
and every $i$,
$\lim_{n\to\infty}\mbbP_n(\varphi_i(\bar{a})) = d_i$.
Suppose that $r$ is chosen so that $r + d_1/d_2 = d_3/d_4$.
Then the formula~(\ref{a formula exemplifying necessity of noncriticality}) is critical
(see Definition~\ref{definition of noncritical formula}).
The intuition is now that for all large enough $n$ and almost all $\mcA \in \mbW_n$,
the numbers $r + |\varphi_1(\mcA) \cap \varphi_2(\mcA)|/|\varphi_2(\mcA)|$  and
$|\varphi_3(\mcA) \cap \varphi_4(\mcA)|/|\varphi_4(\mcA)|$ are very close to each other,
but this does not exclude the possibility that for infinitely many $n$ the first number is at least as large as the second 
and for infinitely many $n$
the second number is larger. In this case the truth value of the 
formula~(\ref{a formula exemplifying necessity of noncriticality})
will alternate between true and false infinitely many times as $n$ tends to infinity
in typical (or ``almost all'') members of $\mbW_n$.

One may also ask if it is necessary in the above theorems that all aggregation formulas $\chi_{R, i}$ are noncritical.
I do not currently know but I assume that the answer is yes.
}\end{rem}

\section{Proof of 
Theorems~\ref{main result on quantifier elimination},
\ref{corollary to the main result on quantifier elimination}
and~\ref{asymptotically equivalent bayesian network}}\label{proofs}

\noindent
Let $\sigma$ be a finite relational signature and $\mfG$ a lifted Bayesian network for $\sigma$.
The proof proceeds by induction on the mp-rank of the underlying DAG of $\mfG$.
The base case will {\em not} be when the mp-rank of $\mfG$ is 0.
Instead the base case will be the ``empty'' lifted Bayesian network for the empty signature $\es$, as described in
Definition~\ref{empty BN}.
In the case of an empty signature (and consequently empty lifted Bayesian network) 
Theorems~\ref{main result on quantifier elimination} --~\ref{asymptotically equivalent bayesian network} 
are a direct consequence of
Lemma~\ref{the case of empty sigma'} below.

The rest of the proof concerns the induction step.
The induction step is proved by Proposition~\ref{completing the inductive step with regard to the BN}
and Corollary~\ref{corollary to completion of the inductive step}
which rely (only) on Assumption~\ref{basic assumptions in the proofs} below 
which states the general assumptions related to the lifted Bayesian network
and Assumption~\ref{inductive assumptions} below which states the induction hypothesis.
Theorems~\ref{main result on quantifier elimination} --~\ref{asymptotically equivalent bayesian network}
 follow from the arguments in this section,
in particular Proposition~\ref{completing the inductive step with regard to the BN}
and Corollary~\ref{corollary to completion of the inductive step},
because 
\begin{itemize}
\item $k \in \mbbN^+$ can be chosen arbitrarily large in 
Lemma~\ref{the case of empty sigma'} and in Assumption~\ref{inductive assumptions},
\item $\varepsilon' > 0$ can be chosen arbitrarily small in 
Lemma~\ref{the case of empty sigma'} and in Assumption~\ref{inductive assumptions}, and
\item because we can choose $\delta'(n) = e^{-dn}$ for any $d > 0$ in Lemma~\ref{the case of empty sigma'} and 
because of the lower bound in 
Lemma~\ref{Y-n has large probability}.
\end{itemize}

\noindent
For the rest of this section we assume the following:

\begin{assump}\label{basic assumptions in the proofs} {\bf (Relationship to a lifted Bayesian network)} {\rm 
\begin{itemize}
\item $\sigma$ is a finite relational signature and $\sigma'$ is a proper subset of $\sigma$.

\item For each $R \in \sigma \setminus \sigma'$, of arity $m$ say, there are a number $\nu_R \in \mbbN$, a sequence of variables
$\bar{x} = (x_1, \ldots, x_m)$ and formulas 
$\chi_{R, i}(\bar{x}) \in CPL(\sigma')$, for $i = 1, \ldots, \nu_R$, such that 
$\forall \bar{x} \big( \bigvee_{i = 1}^{\nu_R} \chi_{R, i}(\bar{x})\big)$ is valid (i.e. true in all $\sigma'$-structures) and if
$i \neq j$ then $\exists \bar{x} \big(\chi_{R, i}(\bar{x}) \wedge \chi_{R, j}(\bar{x})\big)$
is unsatisfiable.

\item For every $R \in \sigma \setminus \sigma'$ and every $1 \leq i \leq \nu_R$, $\mu(R \ | \ \chi_{R, i})$ denotes a real number
in the interval $[0, 1]$. (Sometimes we write $\mu(R(\bar{x}) \ | \ \chi_{R, i}(\bar{x}))$ where $\bar{x}$ is a sequence of variables
the length of which equals the arity of $R$.)

\item For every $\sigma$-structure $\mcA$, every $R \in \sigma \setminus \sigma'$, every $1 \leq i \leq \nu_R$ and every
$\bar{a} \in A^r$ where $r$ is the arity of $R$, let
\[
\lambda(\mcA, R, i, \bar{a}) = 
\begin{cases}
\mu(R \ | \ \chi_{R, i}) \ \ \ \ \ \ \ \text{ if } \mcA \models \chi_{R, i}(\bar{a}) \wedge R(\bar{a}),\\
1 - \mu(R \ | \ \chi_{R, i}) \ \ \text{ if } \mcA \models \chi_{R, i}(\bar{a}) \wedge \neg R(\bar{a}),\\
\text{0 \ \ \ \ \ \ \ \ \ \ \ \ \ \ \ \ \ \ \ \ otherwise.}
\end{cases}
\]

\item For every $n \in \mbbN^+$, $\mbW'_n$ is the set of all $\sigma'$-structures with domain $[n] = \{1, \ldots, n\}$
and $\mbbP'_n$ is a probability distribution on $\mbW'_n$.

\item For every $n \in \mbbN^+$, $\mbW_n$ is the set of all $\sigma$-structures with domain $[n]$.
\end{itemize}
}\end{assump}

\noindent
Recall that, according to 
Definition~\ref{semantics of CPL},
if $\psi(\bar{x}) \in CPL(\sigma')$ and $\mcA \in\mbW_n$ then 
$\psi(\mcA \uhrc \sigma') = \{\bar{b} : \mcA \uhrc \sigma' \models \psi(\bar{b})\}$.

\begin{defin}\label{definition of conditional probability distribution}{\rm 
For every $n \in \mbbN$ and every $\mcA \in \mbW_n$ we define
\[
\mbbP_n(\mcA) \ = \ \mbbP'_n(\mcA \uhrc \sigma') 
\prod_{R \in \sigma \setminus \sigma'} \ \prod_{i=1}^{\nu_R} \ \prod_{\bar{a} \in \chi_{R, i}(\mcA \uhr \sigma')} 
\lambda(\mcA, R, i, \bar{a}).
\]
Then $\mbbP_n$ is a probability distribution on $\mbW_n$ which we may call the 
{\em $\mbbP'_n$-conditional probability distribution on $\mbW_n$}.
}\end{defin}

\begin{notation}\label{remark on notation concerning primes}{\rm
The notation in this section will follow the following pattern:
$\sigma'$-structures, in particular members of $\mbW'_n$, will be denoted $\mcA', \mcB'$, etcetera;
subsets of $\mbW'_n$ will be denoted $\mbX'$ (or $\mbX'_n$), $\mbY'$ (or $\mbY'_n$), etcetera;
$\sigma$-structures and subsets of $\mbW_n$ will be denoted similarly but without the (symbol for) ``prime''.
}\end{notation}

\noindent
In the proofs that follow we will consider ``restrictions'' of $\mbbP_n$ to some subsets of $\mbW_n$ according to the next definition.

\begin{defin}{\rm
(i) If $\mbY' \subseteq \mbW'_n$ then we define
\begin{align*}
&\mbW^{\mbY'} = \{\mcA \in \mbW_n : \mcA \uhrc \sigma' \in \mbY'\} \quad \text{ and} \\
&\mbbP^{\mbY'}(\mcA) \ = \
\frac{\mbbP'_n(\mcA \uhrc \sigma')}{\mbbP'_n(\mbY')} 
\prod_{R \in \sigma \setminus \sigma'} \ \prod_{i=1}^{\nu_R} \ \prod_{\bar{a} \in \chi_{R, i}(\mcA)} 
\lambda(\mcA, R, i, \bar{a}).
\end{align*}
(ii) If $\mcA' \in \mbW'_n$, then we let
\begin{align*}
&\mbW^{\mcA'} = \mbW^{\{\mcA'\}} \quad \text{ and, for every $\mcA \in \mbW^{\mcA'}$,} \\
&\mbbP^{\mcA'}(\mcA) = \mbbP^{\{\mcA'\}}(\mcA) = 
\prod_{R \in \sigma \setminus \sigma'} \ \prod_{i=1}^{\nu_R} \ \prod_{\bar{a} \in \chi_{R, i}(\mcA)} 
\lambda(\mcA, R, i, \bar{a}).
\end{align*}
}\end{defin}

\noindent
Then $\mbbP^{\mbY'}$ and $\mbbP^{\mcA'}$ are probability distributions on $\mbW^{\mbY'}$ and $\mbW^{\mcA'}$, respectively;
if this is not clear see Remark~\ref{remark on P restricted to A'} below.
Note also that if $\mbY' \subseteq \mbW'_n$, $\mcA' \in \mbY'$ and $\mcA \in \mbW^{\mcA'}$, then
\begin{equation}\label{first basic equality about restricted distributions}
\mbbP^{\mbY'}(\mcA) = \frac{\mbbP'_n(\mcA')}{\mbbP'_n(\mbY')}\mbbP^{\mcA'}(\mcA),
\end{equation}
and in particular, taking $\mbY' = \mbW'_n$, we have, for every $\mcA \in \mbW_n$,
\begin{equation}\label{second basic equality about restricted distributions}
\mbbP_n(\mcA) = \mbbP'_n(\mcA \uhrc \sigma')\mbbP^{\mcA \uhr \sigma'}(\mcA).
\end{equation}

\noindent
We now state a few basic lemmas which will be useful.

\begin{lem}\label{P and P' agree when constraints speak only about W'}
For every $n$, if $\mbY' \subseteq \mbW'_n$ then $\mbbP_n(\mbW^{\mbY'}) = \mbbP'_n(\mbY')$.
\end{lem}

\noindent
{\bf Proof.} By using~(\ref{second basic equality about restricted distributions}) in the first line below we get
\begin{align*}
&\mbbP_n(\mbW^{\mbY'}) \ = \ \sum_{\mcA' \in \mbY'} \sum_{\mcA \in \mbW^{\mcA'}} \mbbP_n(\mcA) \ = \ 
\sum_{\mcA' \in \mbY'} \sum_{\mcA \in \mbW^{\mcA'}} \mbbP'_n(\mcA') \mbbP^{\mcA'}(\mcA) \ = \\ 
&\sum_{\mcA' \in \mbY'}  \mbbP'_n(\mcA')  \sum_{\mcA \in \mbW^{\mcA'}} \mbbP^{\mcA'}(\mcA) \ = \ 
\sum_{\mcA' \in \mbY'}  \mbbP'_n(\mcA') \ = \ \mbbP'_n(\mbY').
\end{align*}
\hfill $\square$

\begin{lem}\label{P(X given Y) equals P_Y(X cut Y)}
For every $n$, \\
(i) if $\mbX \subseteq \mbW_n$ and $\mcA' \in \mbW'_n$, 
then $\mbbP_n(\mbX \ | \ \mbW^{\mcA'}) = \mbbP^{\mcA'}(\mbX \cap \mbW^{\mcA'})$, and\\
(ii) if $\mbX \subseteq \mbW_n$ and $\mbY' \subseteq \mbW'_n$, 
then $\mbbP_n(\mbX \ | \ \mbW^{\mbY'}) = \mbbP^{\mbY'}(\mbX \cap \mbW^{\mbY'})$.
\end{lem}

\noindent
{\bf Proof.} Let $\mbX \subseteq \mbW_n$.

(i) Let $\mcA' \in \mbW'_n$. Using Lemma~\ref{P and P' agree when constraints speak only about W'} in the first line below
and~(\ref{second basic equality about restricted distributions})) in the second line below, we get
\begin{align*}
&\mbbP_n(\mbX \ | \ \mbW^{\mcA'}) \ = \
\frac{\mbbP_n(\mbX \cap \mbW^{\mcA'})}{\mbbP_n(\mbW^{\mcA'})} \ = \ 
\frac{\mbbP_n(\mbX \cap \mbW^{\mcA'})}{\mbbP'_n(\mcA')} \ = \\ 
&\frac{\mbbP'_n(\mcA') \sum_{\mcA \in \mbX \cap \mbW^{\mcA'}} \mbbP^{\mcA'}(\mcA)}{\mbbP'_n(\mcA')} \ = \ 
\mbbP^{\mcA'}(\mbX \cap \mbW^{\mcA'}).
\end{align*}

(ii) Let $\mbY' \subseteq \mbW'_n$. Using that $\mbX \cap \mbW^{\mbY'}$ is the disjoint union of all
$\mbX \cap \mbW^{\mcA'}$ such that $\mcA' \in \mbY'$, 
Lemma~\ref{P and P' agree when constraints speak only about W'},
part~(i) of this lemma 
and~(\ref{first basic equality about restricted distributions}), we get
\begin{align*}
&\mbbP_n(\mbX \ | \ \mbW^{\mbY'}) \ = \ 
\frac{\mbbP_n(\mbX \cap \mbW^{\mbY'})}{\mbbP_n(\mbW^{\mbY'})} \ = \ 
  \sum_{\mcA' \in \mbY'} \frac{\mbbP_n(\mbX \cap \mbW^{\mcA'})}{\mbbP_n(\mbW^{\mbY'})} \ = \\ 
&\sum_{\mcA' \in \mbY'} \frac{\mbbP_n(\mbW^{\mcA'})}{\mbbP_n(\mbW^{\mbY'})} \ \mbbP_n(\mbX \ | \ \mbW^{\mcA'}) \ = \ 
\sum_{\mcA' \in \mbY'} \frac{\mbbP'_n(\mcA')}{\mbbP'_n(\mbY')} \mbbP^{\mcA'}(\mbX \cap \mbW^{\mcA'}) \ = \\ 
&\sum_{\mcA' \in \mbY'} \mbbP^{\mbY'}(\mbX \cap \mbW^{\mcA'}) \  = \ 
\mbbP^{\mbY'}(\mbX \cap \mbW^{\mbY'}).
\end{align*}
\hfill $\square$

\begin{rem}\label{remark on P restricted to A'}{\rm ({\bf About $\mbbP^{\mcA'}$})
Fix any $n$ and any $\mcA' \in \mbW'_n$.
For every  $R \in \sigma \setminus \sigma'$, every $1 \leq i \leq \nu_R$ and every $\bar{a} \in \chi_{R, i}(\mcA')$,
let $\Omega(R, i, \bar{a}) = \{0, 1\}$ and let $\mbbP_{R, i, \bar{a}}$ be the probability distribution on $\Omega(R, i, \bar{a})$ with
$\mbbP_{R, i, \bar{a}}(1) = \mu(R \ | \ \chi_{R, i})$.
Then let $\mbbP_\Omega$ be the product measure on 
\[
\Omega \ = \ \prod_{\substack{R \in \sigma \setminus \sigma' \\ 1 \leq i \leq \nu_R \\ \bar{a} \in \chi_{R, i}(\mcA')}} \Omega(R, i, \bar{a}).
\]
Consider the map which sends $\mcA \in \mbW^{\mcA'}$ to the finite sequence
\[
\bar{\kappa}_\mcA \ = \ 
\big(\kappa(R, i, \bar{a}) : R \in \sigma \setminus \sigma', 1 \leq i \leq \nu_R, \bar{a} \in \chi_{R, i}(\mcA') \big)
\]
where $\kappa(R, i, \bar{a}) = 1$ if $\mcA \models R(\bar{a})$ and $\kappa(R, i, \bar{a}) = 0$ otherwise.
This map is clearly a bijection from $\mbW^{\mcA'}$ to $\Omega$ and, for every $\mcA \in \mbW^{\mcA'}$,
$\mbbP^{\mcA'}(\mcA) = \mbbP_\Omega(\bar{\kappa}_\mcA)$.

For every $\alpha \in \{0, 1\}$, every $R \in \sigma \setminus \sigma'$ and every $\bar{a} \in [n]$ (having the same length as the arity of $R$),
let
$\mbE^\alpha_{R, \bar{a}} = \{\mcA \in \mbW^{\mcA'} : \mcA \models R^\alpha(\bar{a})\}$
where $R^0$ and $R^1$ denote $\neg R$ and $R$, respectively.
From the connection to the product measure it follows that 
\begin{itemize}
\item[(a)] for every $R \in \sigma \setminus \sigma'$, every $1 \leq i \leq \nu_R$ and every $\bar{a} \in \chi_{R, i}(\mcA')$,
$\mbbP^{\mcA'}(\mbE^1_{R, \bar{a}}) = \mu\big(R \ | \ \chi_{R, i}\big)$, and
\item[(b)] if $\alpha_1, \ldots, \alpha_m \in \{0, 1\}$, $R_1, \ldots, R_m \in \sigma \setminus \sigma'$
and $\bar{a}_1, \ldots, \bar{a}_m$ are tuples where $|\bar{a}_i|$ is the arity of $R_i$ for each $i$,
and for all  $1 \leq i < j \leq m$,  $R_i \neq R_j$ or $\bar{a}_i \neq \bar{a}_j$, then
the events $\mbE^{\alpha_1}_{R_1, \bar{a}_1}, \ldots, \mbE^{\alpha_m}_{R_m, \bar{a}_m}$
are independent.
\end{itemize}
}\end{rem}

\noindent
The next lemma is a direct consequence of~(b) of Remark~\ref{remark on P restricted to A'}.

\begin{lem}\label{independence in P restricted to expansions of A', generalized}
Suppose that $p(x_1, \ldots, x_m)$ and $q(x_1, \ldots, x_m)$ are (possibly partial) atomic $(\sigma \setminus \sigma')$-types.
Also assume that if 
$\varphi$ is an atomic $\sigma$-formula which does not have the form $x = x$ or the form $\top$ and $\varphi \in p$ or $\neg\varphi \in p$,
then neither $\varphi$ nor $\neg\varphi$ belongs to $q$.
Then, for every $n$, every $\mcA' \in \mbW'_n$ and all distinct $a_1, \ldots, a_m \in [n]$, 
the event $\{\mcA \in \mbW^{\mcA'} : \mcA \models p(a_1, \ldots, a_m)\}$ is independent from the event
$\{\mcA \in \mbW^{\mcA'} : \mcA \models q(a_1, \ldots, a_m)\}$ in the probability space $(\mbW_n^{\mcA'}, \mbbP^{\mcA'})$.
\end{lem}

\begin{defin}\label{definition of saturation}{\bf (Saturation and unsaturation)} {\rm
Let $\bar{x}$ and $\bar{y}$ be sequences of different variables such that 
$\rng(\bar{x}) \cap \rng(\bar{y}) = \es$ and let
$p(\bar{x}, \bar{y})$ and $q(\bar{x})$ be atomic $\sigma$-types such that $q \subseteq p$.
Let also $0 \leq \alpha \leq 1$ and $d = \dim_{\bar{y}}(p)$.
\begin{itemize}
\item[(a)] A finite $\sigma$-structure $\mcA$
is called {\em $(p, q, \alpha)$-saturated} if, whenever $\bar{a} \in A^{|\bar{x}|}$ and $\mcA \models q(\bar{a})$, then
$\big| \{\bar{b} \in A^{|\bar{y}|} : \mcA \models p(\bar{a}, \bar{b})\} \big| \geq \alpha |A|^d$.

\item[(b)] A finite $\sigma$-structure $\mcA$
is called {\em $(p, q, \alpha)$-unsaturated} if, whenever $\bar{a} \in A^{|\bar{x}|}$ and $\mcA \models q(\bar{a})$, then
$\big| \{\bar{b} \in A^{|\bar{y}|} : \mcA \models p(\bar{a}, \bar{b})\} \big| \leq \alpha |A|^d$.
\end{itemize}
If $p'(\bar{x}, \bar{y})$ and $q'(\bar{x})$ are atomic $\sigma'$-types and $q' \subseteq p'$, 
then the notions of {\em $(p', q', \alpha)$-saturated} and {\em $(p', q', \alpha)$-unsaturated} are defined in the same way, but considering
finite $\sigma'$-structures instead.
}\end{defin}

\begin{assump}\label{inductive assumptions}{\bf (Induction hypothesis)} {\rm 
Suppose that $k \in \mbbN^+$, $\varepsilon' > 0$, $\delta' : \mbbN^+ \to \mbbR^{\geq 0}$ and 
$\mbY'_n \subseteq \mbW'_n$, for $n \in \mbbN^+$, are such that the following hold:
\begin{itemize}
\item[(1)] $\lim_{n\to\infty} \delta'(n) = 0$.

\item[(2)] $\mbbP'_n(\mbY'_n) \geq 1 - \delta'(n)$ for all sufficiently large $n$.

\item[(3)] For every complete atomic $\sigma'$-type $p'(\bar{x})$ with $|\bar{x}| \leq k$
there is a number which we denote 
$\msfP'(p'(\bar{x}))$, or just $\msfP'(p')$, such that for all sufficiently large $n$ and all $\bar{a} \in [n]$ 
which realize the identity fragment of $p'$,
\[
\big| \mbbP'_n\big(\{\mcA' \in \mbW'_n : \mcA' \models p'(\bar{a})\}\big) \ -  \ \msfP'(p'(\bar{x})) \big| \ \leq \ \delta'(n).
\]

\item[(4)] For every complete atomic $\sigma'$-type $p'(\bar{x}, \bar{y})$ with $|\bar{x}\bar{y}| \leq k$ and
$0 < \dim_{\bar{y}}(p'(\bar{x}, y)) = |\bar{y}|$,
if $q'(\bar{x}) = p' \uhrc \bar{x}$ and $\msfP'(q') > 0$, then for all sufficiently large $n$, every
$\mcA' \in \mbY'_n$ is $(p', q', \alpha/(1 + \varepsilon'))$-saturated and $(p', q', \alpha(1 + \varepsilon'))$-unsaturated
if $\alpha = \msfP'(p'(\bar{x}, \bar{y})) / \msfP'(q'(\bar{x}))$.

\item[(5)] For every $\chi_{R, i}(\bar{x})$ as in Assumption~\ref{basic assumptions in the proofs}
there is a quantifier-free $\sigma'$-formula 
$\chi_{R, i}^*(\bar{x})$ such that for all sufficiently large $n$ and all $\mcA' \in \mbY'_n$,
$\mcA' \models \forall \bar{x} \big(\chi_{R, i}(\bar{x}) \leftrightarrow \chi^*_{R, i}(\bar{x})\big)$.

\end{itemize}
}\end{assump}

\begin{rem}\label{special cases of inductive assumptions} {\bf (Some special cases)} {\rm
(i) As a technical convenience we allow {\em empty types} (and this does not contradict our definition of an atomic type). 
For example, in Definition~\ref{definition of saturation},
we allow the possibility that $\bar{x}$ is an empty sequence and consequently $q(\bar{x}) = \es$ and
$p(\bar{x}, \bar{y})$ is really just $p(\bar{y})$.\\
(ii) For an empty atomic $\sigma'$-type $p'$ we let $\msfP'(p') = 1$ and in this case we also interpret
the set $\{\mcA' \in \mbW'_n : \mcA' \models p'(\bar{a})\}$ as being equal to $\mbW'_n$. 
Then part~(3) of Assumption~\ref{inductive assumptions} makes sense also for a empty type $p'$.\\
(iii) If $p'(\bar{y})$ is a complete atomic $\sigma'$-type and $\msfP'(p') = 0$, then for all sufficiently large $n$ and all $\mcA' \in \mbY'$,
$p'$ is not realized in $\mcA'$ (i.e. $p'(\mcA') = \es$).
The reason is this: Let $\bar{x}$ denote an empty sequence and let $q'(\bar{x})$ be the empty atomic $\sigma'$-type, so 
$q' \subseteq p'$.
For large enough $n$, every $\mcA' \in \mbW'_n$ is $(p', q', \msfP'(p')(1 + \varepsilon'))$-unsaturated by part~(4)
of Assumption~\ref{inductive assumptions}. If $\msfP'(p') = 0$ this implies that $p'$ has no realization in $\mcA$.
}\end{rem}

\begin{lem}\label{uniformity of probabilities under qf}
Suppose that $p'(\bar{x})$ is a complete atomic $\sigma'$-type and that $p(\bar{x}) \supseteq p'(\bar{x})$
is a (possibly partial) atomic $\sigma$-type.
There is a number which we denote $\msfP(p(\bar{x}) \ | \ p'(\bar{x}))$, or just $\msfP(p \ | \ p')$,
such that for all sufficiently large $n$, all $\bar{a} \in [n]$ and all $\mcA' \in \mbY'_n$
such that $\mcA' \models p'(\bar{a})$,
\[
\mbbP^{\mcA'}\big(\{\mcA \in \mbW^{\mcA'} : \mcA \models p(\bar{a})\} \big) = \msfP(p(\bar{x}) \ | \ p'(\bar{x})).
\]
Moreover, the number $\msfP(p(\bar{x}) \ | \ p'(\bar{x}))$ is a product of numbers of the form 
$\mu(R \ | \ \chi_{R, i})$ or $1 - \mu(R \ | \ \chi_{R, i})$.
\end{lem}

\noindent
{\bf Proof.}
Suppose that $\bar{a}, \bar{b} \in [n]$ and $\mcA', \mcB' \in \mbY'_n$ are such that $\mcA' \models p'(\bar{a})$ and $\mcB' \models p'(\bar{b})$.
Let $R \in \sigma \setminus \sigma'$. 
By part~(5) of Assumption~\ref{inductive assumptions}, for each $1 \leq i \leq \nu_R$,
there is a quantifier free formula $\chi^*_{R, i}$ such that (if $n$ is large enough)
$\chi_{R, i}$ is equivalent to $\chi^*_{R, i}$ in every structure in $\mbY'_n$.
It follows that if $\bar{c}'$ and $\bar{d}'$ are subsequences of $\bar{a}$ and $\bar{b}$, respectively, of length equal to the arity of $R$,
then either $\mcA' \models \chi_{R, i}(\bar{c})$ and $\mcB' \models \chi_{R, i}(\bar{d})$,
or $\mcA' \not\models \chi_{R, i}(\bar{c})$ and $\mcB' \not\models \chi_{R, i}(\bar{d})$.
The conclusion of the lemma now follows from~(a) and~(b) of Remark~\ref{remark on P restricted to A'}.
\hfill $\square$

\begin{lem}\label{the case of empty sigma'} {\bf (The base case)} 
For every $k \in \mbbN^+$ and every $\varepsilon' > 0$,
if $\sigma' = \es$, $\mbbP'_n$ is the uniform probability distribution\footnote{
In fact the uniform probability distribution is the only probability distribution on $\mbW'_n$ since 
$\mbW'_n$ is a singleton if $\sigma' = \es$ (which we assume in this lemma).} 
on $\mbW'_n$ for all $n$
and $\delta' : \mbbN^+ \to \mbbR^{\geq 0}$ is any function such that $\lim_{n\to\infty} \delta'(n) = 0$, 
then there are
$\mbY'_n \subseteq \mbW'_n$, for $n \in \mbbN^+$,  such that
(1)--(4) in Assumption~\ref{inductive assumptions} hold.
Moreover, for every  $\varepsilon'$-noncritical\footnote{
In the sense of Definition~\ref{definition of noncritical formula}.}
$\varphi(\bar{x}) \in CPL(\es)$ with $|\bar{x}| + \mr{qr}(\varphi) \leq k$
there is a quantifier-free formula $\varphi^*(\bar{x})$ such that for all sufficiently large $n$
and all $\mcA' \in \mbY'_n$, $\mcA' \models \forall \bar{x} \big(\varphi(\bar{x}) \leftrightarrow \varphi^*(\bar{x})\big)$.
\end{lem}

\noindent
{\bf Proof.}
Suppose that $\sigma' = \es$ and let $k \in \mbbN^+$ and $\varepsilon' > 0$ be given.
Then, for every $n$, $\mbW'_n$ contains a unique structure which is just the set $[n]$ which has probability 1.
Let $\delta' : \mbbN^+ \to \mbbR^{\geq 0}$ be any function such that $\lim_{n\to\infty} \delta'(n) = 0$.
For every complete atomic $\sigma'$-type $p'(\bar{x})$ let $\msfP'(p'(\bar{x})) = 1$.
Observe that, for every $n$, if $\bar{a} \in [n]$ and $\bar{a}$ realizes the identity fragment of $p'(\bar{x})$,
then $\bar{a}$ realizes $p'(\bar{x})$ in the unique $\mcA'$ of $\mbW'_n$.
Hence, for trivial reasons we have~(3).

For every $n$ let $\mbY'_n$ be the set of all $\mcA' \in \mbW'_n$ such that for every
complete atomic $\sigma'$-type $p'(\bar{x}, \bar{y})$ with $|\bar{x}\bar{y}| \leq k$ and
$0 < \dim_{\bar{y}}(p'(\bar{x}, y)) = |\bar{y}|$,
if $q(\bar{x}) = p \uhrc \bar{x}$, then for all sufficiently large $n$, every
$\mcA' \in \mbY'_n$ is $(p', q', 1/(1 + \varepsilon'))$-saturated and $(p', q', (1 + \varepsilon'))$-unsaturated.
Suppose that $p'(\bar{x}, \bar{y})$ is a complete atomic $\sigma'$-type with $|\bar{x}\bar{y}| \leq k$ and
$0 < \dim_{\bar{y}}(p'(\bar{x}, y)) = |\bar{y}|$.
Let $q'(\bar{x}) = p' \uhrc \bar{x}$ and suppose that $\mcA' \models q'(\bar{a})$ where $\mcA' \in \mbW'_n$.
Then $\mcA' \models p'(\bar{a}, \bar{b})$ for every $\bar{b} \in [n]$ consisting of different elements no one of which 
occurs in $\bar{a}$. There are $n^{|\bar{y}|} - Cn^{|\bar{y}|-1}$ such $\bar{b}$ for some constant $C$.
So if $n^{|\bar{y}|} - Cn^{|\bar{y}|-1} \geq \frac{n^{|\bar{y}|}}{1 + \varepsilon'}$ then
$\mcA'$ is $(p', q', 1/(1 + \varepsilon'))$-saturated.
For trivial reasons, $\mcA'$ is also $(p', q', (1 + \varepsilon'))$-unsaturated.
Hence, we have proved~(4).
The last claim of the lemma follows from 
Proposition~\ref{quantifier elimination in a structure}
the proof of which works out in exactly the same way if $\sigma$ and $\mbY_n$ (in that proof) is replaced by $\sigma'$
and $\mbY'_n$, respectively,
and we assume~(4). In other words, the almost everywhere elimination of quantifiers follows
from the saturation and unsaturation properties stated in~(4).
\hfill $\square$

\begin{lem}\label{an upper bound}
Suppose that $\mbX_n \subseteq \mbW_n$.
Then for all sufficiently large $n$, $\mbbP_n(\mbX_n) \leq \mbbP_n(\mbX_n \cap \mbW^{\mbY'_n}) + \delta'(n)$.
\end{lem}

\noindent
{\bf Proof.}
We have 
\[
\mbP_n(\mbX_n) = \mbbP_n(\mbX_n \cap \mbW^{\mbY'_n}) + \mbbP_n(\mbX_n \setminus \mbW^{\mbY'_n})
\]
and, using
Lemma~\ref{P and P' agree when constraints speak only about W'},
we have
\[
\mbbP_n(\mbW_n \setminus \mbW^{\mbY'_n})
= 1 - \mbbP_n(\mbW^{\mbY'_n}) = 1 - \mbbP'_n(\mbY'_n) \leq \delta'(n).
\]
Hence $\mbbP_n(\mbX_n) \leq \mbbP_n(\mbX_n \cap \mbW^{\mbY'_n}) + \delta'(n)$.
\hfill $\square$

\begin{lem}\label{uniformity of probabilities under qf, second version}
Suppose that $p'(\bar{x})$ is a complete atomic $\sigma'$-type and that $p(\bar{x}) \supseteq p'(\bar{x})$
is an (possibly partial) atomic $\sigma$-type.
Letting $n$ be sufficiently large, then for all $\bar{a} \in [n]$ and letting
$\mbZ'_n$ be the set of all $\mcA' \in \mbY'_n$ such that $\mcA' \models p'(\bar{a})$ we have
\begin{align*}
&\mbbP_n\big(\{\mcA \in \mbW_n : \mcA \models p(\bar{a})\} \ | \ 
\mbW^{\mbY'_n} \cap \{\mcA \in \mbW_n : \mcA \models p'(\bar{a})\}\big) = \\
&\mbbP^{\mbZ'_n}\big(\{\mcA \in \mbW^{\mbZ'_n} : \mcA \models p(\bar{a})\}\big) = \\
&\msfP(p(\bar{x}) \ | \ p'(\bar{x}))
\end{align*}
where $\msfP(p(\bar{x}) \ | \ p'(\bar{x}))$ is like in Lemma~\ref{uniformity of probabilities under qf}.
\end{lem}

\noindent
{\bf Proof.}
For every $\mcA \in \mbW_n$ we have $\mcA \models p'(\bar{a})$ if and only if $\mcA \uhrc \sigma' \models p'(\bar{a})$.
Therefore $\mbW^{\mbY'_n} \cap \{\mcA \in \mbW_n : \mcA \models p'(\bar{a})\} = \mbW^{\mbZ'_n}$.
By Lemma~\ref{P(X given Y) equals P_Y(X cut Y)}
we have 
\begin{align*}
&\mbbP_n\big(\{\mcA \in \mbW_n : \mcA \models p(\bar{a})\} \ | \ 
\mbW^{\mbY'_n} \cap \{\mcA \in \mbW_n : \mcA \models p'(\bar{a})\}\big) = \\
&\mbbP^{\mbZ'_n}\big(\{\mcA \in \mbW^{\mbZ'_n} : \mcA \models p(\bar{a})\}\big).
\end{align*}
Then, using~(\ref{first basic equality about restricted distributions}) and 
Lemma~\ref{uniformity of probabilities under qf}, we get
\begin{align*}
&\mbbP^{\mbZ'_n}\big(\{\mcA \in \mbW^{\mbZ'_n} : \mcA \models p(\bar{a})\}\big) = 
\sum_{\mcA' \in \mbZ'_n} \mbbP^{\mbZ'_n}\big(\mcA \in \mbW^{\mcA'} : \mcA \models p(\bar{a})\big) = \\
&\sum_{\mcA' \in \mbZ'_n} \frac{\mbbP'_n(\mcA')}{\mbbP'_n(\mbZ'_n)} 
\mbbP^{\mcA'}\big(\{\mcA \in \mbW^{\mcA'} : \mcA \models p(\bar{a})\}\big) = \\
&\sum_{\mcA' \in \mbZ'_n} \frac{\mbbP'_n(\mcA')}{\mbbP'_n(\mbZ'_n)}
\msfP(p(\bar{x}) \ | \ p'(\bar{x})) = 
\msfP(p(\bar{x}) \ | \ p'(\bar{x})).
\end{align*}
\hfill $\square$

\begin{lem}\label{P(p) is a product}
Suppose that $p'(\bar{x})$ is a complete atomic $\sigma'$-type and that $p(\bar{x}) \supseteq p'(\bar{x})$
is a (possibly partial) atomic $\sigma$-type.
Then for all sufficiently large $n$ and all $\bar{a} \in [n]$ which realize the identity fragment of $p'(\bar{x})$ (and hence of $p$)
we have
\[
\big| \mbbP_n\big(\{\mcA \in \mbW_n : \mcA \models p(\bar{a})\} \ | \ \mbW^{\mbY'_n} \big) \ - \
\msfP(p(\bar{x} \ | \ p'(\bar{x})) \cdot \msfP'(p'(\bar{x})) \big| \ < \ 3\delta'(n).
\]
\end{lem}

\noindent
{\bf Proof.}
Let $\bar{a} \in [n]$ realize the identity fragment of $p'(\bar{x})$. Furthermore,
\begin{itemize}
\item[] let $\mbX_n$ be the set of all $\mcA \in \mbW_n$ such that $\mcA \models p(\bar{a})$,
\item[] let $\mbX'_n$ be the set of all $\mcA' \in \mbW'_n$ such that $\mcA' \models p'(\bar{a})$, and
\item[] let $\mbZ'_n$ be the set of all $\mcA' \in \mbY'_n$ such that $\mcA' \models p'(\bar{a})$.
\end{itemize}
From parts~(2) and~(3) of Assumption~\ref{inductive assumptions} it easily follows that (for large enough $n$)
\begin{itemize}
\item[] $\mbbP'_n(\mbZ'_n) / \mbbP'_n(\mbY'_n)$ differs from $\mbbP'_n(\mbZ'_n)$ by at most $\delta'(n)$,
\item[] $\mbbP'_n(\mbZ'_n)$ differs from $\mbbP'_n(\mbX'_n)$ by at most $\delta'(n)$ and
\item[] $\mbbP'_n(\mbX'_n)$ differs from $\msfP'(p'(\bar{x}))$ by at most $\delta'(n)$.
\end{itemize}
By Lemma~\ref{P(X given Y) equals P_Y(X cut Y)},
$\mbbP_n(\mbX_n \ | \ \mbW^{\mbY'_n}) = \mbbP^{\mbY'_n}(\mbX \cap \mbW^{\mbY'_n})$.
Then, using~(\ref{first basic equality about restricted distributions}) 
and Lemma~\ref{uniformity of probabilities under qf}, we have
\begin{align*}
&\mbbP^{\mbY'_n}\big(\mbX \cap \mbW^{\mbY'_n}\big) = 
\sum_{\mcA' \in \mbY'_n} \mbbP^{\mbY'_n}\big(\mbX_n \cap \mbW^{\mcA'}\big) =
\sum_{\mcA' \in \mbZ'_n} \mbbP^{\mbY'_n}\big(\mbX_n \cap \mbW^{\mcA'}\big) = \\
&\sum_{\mcA' \in \mbZ'_n} \ \sum_{\mcA \in \mbX_n \cap \mbW^{\mcA'}}\mbbP^{\mbY'_n}(\mcA) = 
\sum_{\mcA' \in \mbZ'_n} \ \sum_{\mcA \in \mbX_n \cap \mbW^{\mcA'}}
\frac{\mbbP'_n(\mcA')}{\mbbP'_n(\mbY'_n)} \mbbP^{\mcA'}(\mcA) = \\
&\sum_{\mcA' \in \mbZ'_n}\frac{\mbbP'_n(\mcA')}{\mbbP'_n(\mbY'_n)}\sum_{\mcA \in \mbX_n \cap \mbW^{\mcA'}}
\mbbP^{\mcA'}(\mcA) = 
\sum_{\mcA' \in \mbZ'_n}
\frac{\mbbP'_n(\mcA')}{\mbbP'_n(\mbY'_n)}\mbbP^{\mcA'}\big(\mbX_n \cap \mbW^{\mcA'}\big) = \\
&\sum_{\mcA' \in \mbZ'_n}
\frac{\mbbP'_n(\mcA')}{\mbbP'_n(\mbY'_n)} \msfP(p(\bar{x}) \ | \ p'(\bar{x})) = 
\msfP(p(\bar{x}) \ | \ p'(\bar{x})) 
\sum_{\mcA' \in \mbZ'_n}\frac{\mbbP'_n(\mcA')}{\mbbP'_n(\mbY'_n)} = \\
&\msfP(p(\bar{x}) \ | \ p'(\bar{x})) \ \frac{\mbbP'_n(\mbZ'_n)}{\mbbP'_n(\mbY'_n)},
\end{align*}
where 
\[
\msfP'(p'(\bar{x})) - 3\delta'(n) \leq \frac{\mbbP'_n(\mbZ'_n)}{\mbbP'_n(\mbY'_n)} \leq \msfP'(p'(\bar{x})) + 3\delta'(n).
\]
\hfill $\square$

\begin{lem}\label{P(p) is a product, version 2}
Suppose that $p'(\bar{x})$ is a complete atomic $\sigma'$-type and that $p(\bar{x}) \supseteq p'(\bar{x})$
is an (possibly partial) atomic $\sigma$-type.
Then for all sufficiently large $n$ and all $\bar{a} \in [n]$ which realize the identity fragment of $p'(\bar{x})$ we have
\[
\big| \mbbP_n\big(\{\mcA \in \mbW_n : \mcA \models p(\bar{a})\}\big) \ - \
\msfP(p(\bar{x}) \ | \ p'(\bar{x})) \cdot \msfP'(p'(\bar{x})) \big| \ < \ 5\delta'(n).
\]
\end{lem}

\noindent
{\bf Proof.}
Let $\bar{a} \in [n]$ realize the identity fragment of $p'(\bar{x})$.
Let $\mbX_n$ be the set of all $\mcA \in \mbW_n$ such that $\mcA \models p(\bar{a})$.
We have 
\[
\mbbP_n\big(\mbX_n\big) = \mbbP_n\big(\mbX_n \ | \ \mbW^{\mbY'_n}\big) \mbbP_n\big(\mbW^{\mbY'_n}\big)
+ 
\mbbP_n\big(\mbX \ | \ \mbW_n \setminus \mbW^{\mbY'_n}\big) \mbbP_n\big(\mbW_n \setminus \mbW^{\mbY'_n}\big).
\]
By the use of Lemma~\ref{P and P' agree when constraints speak only about W'} and by
part~(2) of Assumption~\ref{inductive assumptions},
we also have
\[
\mbbP_n\big(\mbW_n \setminus \mbW^{\mbY'_n}\big) = 1 - \mbbP_n\big(\mbW^{\mbY'_n}\big) =
1 - \mbbP'_n(\mbY'_n) \leq \delta'(n).
\]
It follows that 
$\mbbP_n\big(\mbX \ | \ \mbW_n \setminus \mbW^{\mbY'_n}\big) \mbbP_n\big(\mbW_n \setminus \mbW^{\mbY'_n}\big)
\leq \delta'(n)$.
By Lemma~\ref{P and P' agree when constraints speak only about W'} 
and part~(2) of Assumption~\ref{inductive assumptions},
$\mbbP_n\big(\mbW^{\mbY'_n}\big) = \mbbP'_n\big(\mbY'_n\big) \geq 1 - \delta'(n)$.
It now follows from 
Lemma~\ref{P(p) is a product}
that $\mbbP_n\big(\mbX_n\big)$ differs from 
$\msfP(p(\bar{x} \ | \ p'(\bar{x})) \cdot \msfP'(p'(\bar{x}))$ 
by at most $5\delta'(n)$ (for sufficiently large $n$).
\hfill $\square$

\begin{defin}\label{definition of P(p)}{\rm
For every (possibly partial) $\sigma$-type $p(\bar{x})$ such that 
$p'(\bar{x}) = p \uhrc \sigma'$ is a complete atomic $\sigma'$-type, we define
$
\msfP(p(\bar{x})) = 
\msfP'(p'(\bar{x})) \cdot \msfP(p(\bar{x}) \ | \ p'(\bar{x})).
$
}\end{defin}

\noindent
With this definition we can reformulate Lemma~\ref{P(p) is a product, version 2} as follows:

\begin{cor}\label{probability of p converges to P(p)}
If $p(\bar{x})$ is an (possibly partial) atomic $\sigma$-type such that $p \uhrc \sigma'$  is a complete atomic $\sigma'$-type,
then, for all sufficiently large $n$ and all $\bar{a} \in [n]$ which realize the identity fragment of $p(\bar{x})$ we have
\[
\big| \mbbP_n\big(\{\mcA \in \mbW_n : \mcA \models p(\bar{a})\}\big) \ - \ \msfP(p(\bar{x})) \big| \ < \ 5\delta'(n).
\]
\end{cor}

\begin{lem}\label{an extension and independence}
Suppose that $p(\bar{x}, \bar{y})$ is a complete atomic $\sigma$-type.
Let $p'(\bar{x}, \bar{y}) = p \uhrc \sigma'$, $q(\bar{x}) = p \uhrc \bar{x}$ and
let $p^{\bar{y}}(\bar{x}, \bar{y})$  include $p'(\bar{x}, \bar{y})$ and all formulas in $p$ 
in which at least one variable from $\bar{y}$ occurs.
Then
\[
\msfP(p(\bar{x}, \bar{y}) \ | \ p'(\bar{x}, \bar{y})) = 
\msfP(p^{\bar{y}}(\bar{x}, \bar{y}) \ | \ p'(\bar{x}, \bar{y})) \cdot \msfP(q(\bar{x}) \ | \ p'(\bar{x}, \bar{y})).
\]
\end{lem}

\noindent
{\bf Proof.}
By Lemma~\ref{uniformity of probabilities under qf}, 
for any sufficently large $n$, any $\bar{a}, \bar{b} \in [n]$
and any $\mcA' \in \mbY'_n$ such that $\mcA' \models p'(\bar{a}, \bar{b})$, we have
\begin{align*}
&\msfP(p(\bar{x}, \bar{y}) \ | \ p'(\bar{x}, \bar{y})) = 
\mbbP^{\mcA'}\big(\{\mcA \in \mbW^{\mcA'} : \mcA \models p(\bar{a}, \bar{b})\}\big), \\
&\msfP(p^y(\bar{x}, y) \ | \ p'(\bar{x}, \bar{y})) = 
\mbbP^{\mcA'}\big(\{\mcA \in \mbW^{\mcA'} : \mcA \models p^y(\bar{a}, \bar{b})\}\big) \ \text{ and} \\
&\msfP(q(\bar{x}) \ | \ p'(\bar{x}, \bar{y})) = 
\mbbP^{\mcA'}\big(\{\mcA \in \mbW^{\mcA'} : \mcA \models q(\bar{a})\}\big).
\end{align*}
Note that $p(\bar{x}, \bar{y}) = p'(\bar{x}, \bar{y}) \cup p^{\bar{y}}(\bar{x}, \bar{y}) \cup q(\bar{x})$.
By Lemma~\ref{independence in P restricted to expansions of A', generalized},
the event 
$\{\mcA \in \mbW^{\mcA'} : \mcA \models p^{\bar{y}}(\bar{a}, \bar{b})\}$ is independent, in $(\mbW^{\mcA'}, \mbbP^{\mcA'})$, 
from the event
$\{\mcA \in \mbW^{\mcA'} : \mcA \models q(\bar{a})\}$.
Therefore,
\begin{align*}
&\mbbP^{\mcA'}\big(\{\mcA \in \mbW^{\mcA'} : \mcA \models p(\bar{a}, \bar{b}) \}\big) = \\
&\mbbP^{\mcA'}\big(\{\mcA \in \mbW^{\mcA'} : \mcA \models p^{\bar{y}}(\bar{a}, \bar{b})\}\big) \cdot
\mbbP^{\mcA'}\big(\{\mcA \in \mbW^{\mcA'} : \mcA \models q(\bar{a})\}\big)
\end{align*}
and from this the lemma follows. 
\hfill $\square$

\begin{lem}\label{reduntant information on level W'}
Let $p'(\bar{x}, \bar{y})$ be a complete atomic $\sigma'$-type, $q'(\bar{x}) = p' \uhrc \bar{x}$
and suppose that $q(\bar{x})$ is a complete atomic $\sigma$-type such that $q \supseteq q'$.
Then $\msfP(q(\bar{x}) \ | \ p'(\bar{x}, \bar{y})) = \msfP(q(\bar{x}) \ | \ q'(\bar{x}))$.
\end{lem}

\noindent
{\bf Proof.}
Since $q'(\bar{x}) \subseteq p'(\bar{x}, \bar{y})$ it follows
from Lemma~\ref{uniformity of probabilities under qf} that
for any sufficently large $n$, any $\bar{a}, \bar{b} \in [n]$
and any $\mcA' \in \mbY'_n$ such that $\mcA' \models p'(\bar{a}, \bar{b})$, we have
\begin{align*}
&\msfP(q(\bar{x}) \ | \ p'(\bar{x}, \bar{y})) = 
\mbbP^{\mcA'}\big(\{\mcA \in \mbW^{\mcA'} : \mcA \models q(\bar{a})\}\big) \ \text{ and} \\
&\msfP(q(\bar{x}) \ | \ q'(\bar{x})) = 
\mbbP^{\mcA'}\big(\{\mcA \in \mbW^{\mcA'} : \mcA \models q(\bar{a})\}\big).
\end{align*}
Hence $\msfP(q(\bar{x}) \ | \ p'(\bar{x}, \bar{y})) = \msfP(q(\bar{x}) \ | \ q'(\bar{x}))$.
\hfill $\square$
\\

\noindent
In Lemma~\ref{uniformity of probabilities under qf}
we defined the notation $\msfP(p(\bar{x}) \ | \ p'(\bar{x}))$ when the atomic $\sigma$-type
$p$ has no more variables than the complete atomic $\sigma'$-type $p'$.
From Definition~\ref{definition of P(p)} of $\msfP(p(\bar{x}))$ it follows that
$\msfP(p(\bar{x}) \ | \ p'(\bar{x})) = \msfP(p(\bar{x})) / \msfP'(p'(\bar{x}))$.
Now we extend this notation to pairs of $(p(\bar{x}, \bar{y}), q(\bar{x}))$ where 
$p(\bar{x}, \bar{y})$ is a complete atomic $\sigma$-type and $q(\bar{x}) = p \uhrc \bar{x}$.

\begin{defin}\label{limit of conditional probabilities of extension pairs}{\rm
Suppose that $p(\bar{x}, y)$ is a complete atomic $\sigma$-type and let $q(\bar{x}) = p \uhrc \bar{x}$.
We define
\[
\msfP(p(\bar{x}, \bar{y}) \ | \ q(\bar{x})) = \frac{\msfP(p(\bar{x}, \bar{y}))}{\msfP(q(\bar{x}))}.
\]
In the same way, if $p'(\bar{x}, \bar{y})$ is a complete atomic $\sigma'$-type and $q'(\bar{x}) = p' \uhrc \bar{x}$,
then we define
\[
\msfP'(p'(\bar{x}, \bar{y}) \ | \ q'(\bar{x})) = \frac{\msfP'(p'(\bar{x}, \bar{y}))}{\msfP'(q'(\bar{x}))}.
\]
}\end{defin}

\begin{lem}\label{probability of an extension}
Suppose that $p(\bar{x}, \bar{y})$ is a complete atomic $\sigma$-type, let $q(\bar{x}) = p \uhrc \bar{x}$
and let $p^{\bar{y}}(\bar{x}, \bar{y})$ be defined as in 
Lemma~\ref{an extension and independence}.
Then
$\msfP(p(\bar{x}, \bar{y}) \ | \ q(\bar{x})) = 
\msfP'(p'(\bar{x}, \bar{y}) \ | \ q'(\bar{x})) \cdot \msfP(p^{\bar{y}}(\bar{x}, \bar{y}) \ | \ p'(\bar{x}, \bar{y})).$
\end{lem}

\noindent
{\bf Proof.}
Using Definition~\ref{definition of P(p)} and Lemmas~\ref{an extension and independence}
and~\ref{reduntant information on level W'} we get
\begin{align*}
&\msfP(p(\bar{x}, \bar{y}) \ | \ q(\bar{x})) = 
\frac{\msfP(p(\bar{x}, \bar{y}))}{\msfP(q(\bar{x}))} = 
\frac{\msfP'(p'(\bar{x}, \bar{y})) \cdot \msfP(p(\bar{x}, \bar{y}) \ | \ p'(\bar{x}, \bar{y}))} 
{\msfP'(q'(\bar{x})) \cdot \msfP(q(\bar{x}) \ | \ q'(\bar{x}))} = \\
&\frac{\msfP'(p'(\bar{x}, \bar{y})) \cdot \msfP(p^{\bar{y}}(\bar{x}, \bar{y}) \ | \ p'(\bar{x}, \bar{y})) 
\cdot \msfP(q(\bar{x}) \ | \ p'(\bar{x}, \bar{y}))}
{\msfP'(q'(\bar{x})) \cdot \msfP(q(\bar{x}) \ | \ q'(\bar{x}))} = \\
&\frac{\msfP'(p'(\bar{x}, \bar{y})) \cdot \msfP(p^{\bar{y}}(\bar{x}, \bar{y}) \ | \ p'(\bar{x}, \bar{y})) 
\cdot \msfP(q(\bar{x}) \ | \ q'(\bar{x}))}
{\msfP'(q'(\bar{x})) \cdot \msfP(q(\bar{x}) \ | \ q'(\bar{x}))} = \\
&\msfP'(p'(\bar{x}, \bar{y}) \ | \ q'(\bar{x})) \cdot \msfP(p^{\bar{y}}(\bar{x}, \bar{y}) \ | \ p'(\bar{x}, \bar{y})).
\end{align*}
\hfill $\square$

\begin{lem}\label{the number of realizations in A'}
Suppose that $n$ is large enough that part~(4) of Assumption~\ref{inductive assumptions} holds.
Suppose that $p(\bar{x}, y)$ and $q(\bar{x})$ are complete atomic $\sigma$-types such that $|\bar{x}y| \leq k$,
$\dim_{y}(p) = 1$ and $q \subseteq p$. 
Let $\gamma = \msfP(p(\bar{x}, y) \ | \ q(\bar{x}))$ and $\mcA' \in \mbY'_n$.
Then 
\begin{align*}
\mbbP_n^{\mcA'}\big(\{\mcA \in \mbW^{\mcA'} : &\text{ $\mcA$ is $(p, q, \gamma/(1 + \varepsilon')^2)$-saturated}\\
 &\text{and $(p, q, \gamma(1 + \varepsilon')^2)$-unsaturated} \}\big)
\end{align*}
is at least
$
1 \ - \ 2 n^{|\bar{x}|} e^{-c_{\varepsilon'} \gamma n}
$
where the constant $c_{\varepsilon'} > 0$ depends only on $\varepsilon'$.
\end{lem}

\noindent
{\bf Proof.}
Suppose that $p(\bar{x}, y)$ and $q(\bar{x})$ are complete atomic $\sigma$-types such that
$|\bar{x}y| \leq k$, $\dim_y(p) = 1$ and
$q \subseteq p$.
Let $p' = p \uhrc \sigma$ and $q' = q \uhrc \sigma'$.
Moreover, let $p^y(\bar{x}, y)$ include $p'(\bar{x}, y)$ and all $(\sigma \setminus \sigma')$-formulas
in $p(\bar{x}, y)$ which contain the variable $y$.
Also, let 
\begin{align*}
&\alpha = \msfP'(p'(\bar{x}, y) \ | \ q'(\bar{x})), \\ 
&\beta = \msfP(p^y(\bar{x}, y) \ | \ p'(\bar{x}, y)) \quad \text{ and } \\
&\gamma = \msfP(p(\bar{x}, y) \ | \ q(\bar{x})).
\end{align*}
By Lemma~\ref{probability of an extension}
we have $\gamma = \alpha \beta$.

Let $\mcA' \in \mbY'_n$.
By~(4) of 
Assumption~\ref{inductive assumptions} 
$\mcA'$ is $(p', q', \alpha/(1 + \varepsilon'))$-saturated and $(p', q', \alpha(1 + \alpha))$-unsaturated if $n$ is large enough.
For every $\bar{a} \in [n]^{|\bar{x}|}$ let
\[
B'_{\bar{a}} \ = \ \big\{b \in [n] : \mcA' \models p'(\bar{a}, b)\big\}.
\]
By the mentioned (un)saturation property, if $\mcA' \models q'(\bar{a})$ then
$\alpha n/(1 + \varepsilon') \leq |B'_{\bar{a}}| \leq \alpha n (1 + \varepsilon')$.
For every $\bar{a} \in [n]^{|\bar{x}|}$ and every $\mcA \in \mbW^{\mcA'}$ let
\[
B_{\bar{a}, \mcA} \  = \ \big\{ b \in [n] : \mcA \models p^y(\bar{a}, b)\big\}
\]
and note that $B_{\bar{a}, \mcA} \subseteq B'_{\bar{a}}$ for every $\bar{a}$ and every $\mcA \in \mbW^{\mcA'}$.
Let 
\[
\mbX_{\bar{a}} \ = \ \big\{ \mcA \in \mbW^{\mcA'} : 
\text{ either $\mcA \not\models q(\bar{a})$ or 
$\gamma/(1 + \varepsilon')^2 \leq |B_{\bar{a}, \mcA}| \leq \gamma(1 + \varepsilon')^2$} \big\}.
\]
Observe that if $\mcA \in \mbW^{\mcA'}$, $\mcA \models q(\bar{a})$ and $\mcA \models p^y(\bar{a}, b)$,
then $\mcA \models p(\bar{a}, b)$. 
Hence every $\mcA \in \bigcap_{\bar{a} \in [n]^{|\bar{x}|}} \mbX_{\bar{a}}$ is 
$(p, q, \gamma/(1 + \varepsilon')^2)$-saturated and $(p, q, \gamma(1 + \varepsilon')^2)$-unsaturated.

Fix any $\bar{a}$ such that $\mcA' \models q'(\bar{a})$ (and note that $\mcA \models q(\bar{a})$ implies $\mcA' \models q'(\bar{a})$).
By Lemma~\ref{independence in P restricted to expansions of A', generalized},
for all distinct $b, c \in B'_{\bar{a}}$, the events
\[
\mbE_b = \{\mcA \in \mbW^{\mcA'} : \mcA \models p^y(\bar{a}, b)\} \ \text{ and } \
\mbE_c = \{\mcA \in \mbW^{\mcA'} : \mcA \models p^y(\bar{a}, c)\}
\]
are independent.
Moreover,
by Lemma~\ref{uniformity of probabilities under qf},
for each $b \in B'_{\bar{a}}$, $\mbbP^{\mcA'}_n(\mbE_b) = \beta$.
Let $Z : \mbW^{\mcA'} \to \mbbN$ be the random variable defined by
\[
Z(\mcA) = \big| \{b \in B'_{\bar{a}} : \mcA \models p^y(\bar{a}, b) \}\big|.
\]
Let $\varepsilon = \varepsilon' / (1 + \varepsilon')$ and note that $\varepsilon < \varepsilon'$ and $1 - \varepsilon = 1/(1 + \varepsilon')$.
By Lemma~\ref{independent bernoulli trials},
\[
\mbbP^{\mcA'}\big(\big|Z - \beta | B'_{\bar{a}} | \big| > \varepsilon \beta |B'_{\bar{a}}|\big) \ < \  
2\exp\big(-c_{\varepsilon} \beta|B'_{\bar{a}}|\big)
\]
where $c_{\varepsilon}$ depends only on $\varepsilon$ and hence only on $\varepsilon'$.
Recall that $\alpha \beta = \gamma$ and $\alpha n/(1 + \varepsilon') \leq |B'_{\bar{a}}| \leq \alpha n (1 + \varepsilon')$.
From this it follows that $(1 + \varepsilon')^2 \gamma n  \geq (1 + \varepsilon') \beta |B'_{\bar{a}}|$ and
$\gamma n / (1 + \varepsilon')^2  \leq \beta |B'_{\bar{a}}| / (1 + \varepsilon')$.
Therefore, if $Z > (1 + \varepsilon')^2\gamma n$ or $Z < \gamma n / (1 + \varepsilon')^2$, then
$\big|Z - \beta | B'_{\bar{a}} | \big| > \varepsilon \beta |B'_{\bar{a}}|$.
Hence, if $c_{\varepsilon'} = c_\varepsilon / (1 + \varepsilon')$,
\[
\mbbP^{\mcA'}\big(\mbW^{\mcA'} \setminus \mbX_{\bar{a}}\big) \  < \
2\exp\big(-c_{\varepsilon} \beta|B'_{\bar{a}}|\big) \ \leq \ 
2\exp\big(-c_{\varepsilon'} \gamma n\big).
\]

Since the argument works for all $\bar{a} \in [n]^{|\bar{x}|}$ such that $\mcA' \models q'(\bar{a})$ it follows that
\[
\mbbP^{\mcA'}\bigg(\bigcap_{\bar{a} \in [n]^{|\bar{x}|}} \mbX_{\bar{a}}\bigg) \ \geq \
1 - 2 n^{|\bar{x}|} e^{-c_{\varepsilon'} \gamma n}
\]
and this proves the lemma.
\hfill $\square$
\\

\noindent
The next lemma generalizes the previous one to types $p(\bar{x}, \bar{y})$ where the length of
$\bar{y}$ is greater than one.

\begin{lem}\label{the number of realizations in A' for y-bar}
Suppose that $n$ is large enough that part~(4) of Assumption~\ref{inductive assumptions} holds.
Suppose that $p(\bar{x}, \bar{y})$ and $q(\bar{x})$ are complete atomic $\sigma$-types such that $|\bar{x}\bar{y}| \leq k$,
$\dim_{\bar{y}}(p) = |\bar{y}|$ and $q \subseteq p$.
Let $\gamma = \msfP(p(\bar{x}, \bar{y}) \ | \ q(\bar{x}))$ and $\mcA' \in \mbY'_n$.
Then 
\begin{align*}
\mbbP_n^{\mcA'}\big(\{\mcA \in \mbW^{\mcA'} : &\text{ $\mcA$ is $(p, q, \gamma/(1 + \varepsilon')^{2|\bar{y}|})$-saturated}\\
 &\text{and $(p, q, \gamma(1 + \varepsilon')^{2|\bar{y}|})$-unsaturated} \}\big)
\end{align*}
is at least
$
1 \ - \ 2^{|\bar{y}|} n^{|\bar{x}| + |\bar{y}| - 1} e^{-c_{\varepsilon'} \gamma n}
$
where the constant $c_{\varepsilon'} > 0$ depends only on $\varepsilon'$.
\end{lem}

\noindent
{\bf Proof.}
We prove the lemma by induction on $m = |\bar{y}|$. The base case $m = 1$ is given by
Lemma~\ref{the number of realizations in A'}.
Let $p(\bar{x}, \bar{y})$ and $q(\bar{x})$ be as assumed in the lemma where $\bar{y} = (y_1, \ldots, y_{m+1})$.
Let $p_m(\bar{x}, y_1, \ldots, y_m)$ be the restriction of $p$ to formulas with variables among $\bar{x}, y_1, \ldots, y_m$.
Furthermore, let $\alpha = \msfP(p_m \ | \ q)$, $\beta = \msfP(p \ | \ p_m)$ and
$\gamma = \msfP(p \ | \ q)$.
Observe that by Definition~\ref{limit of conditional probabilities of extension pairs} we have
\[
\gamma = \frac{\msfP(p)}{\msfP(q)} = \frac{\msfP(p)}{\msfP(p_m)} \cdot \frac{\msfP(p_m)}{\msfP(q)} = \beta \alpha.
\]

Let $\mcA' \in \mbY'_n$.
By the induction hypothesis,
the probability (with the distribution $\mbbP^{\mcA'}$) that 
\begin{itemize}
\item[(a)] $\mcA \in \mbW^{\mcA'}$ is 
$(p_m, q, \alpha/(1 + \varepsilon')^{2m})$-saturated and 
$(p_m, q, \alpha(1 + \varepsilon')^{2m})$-unsaturated
\end{itemize}
 is at least
$1 \ - \ 2^m n^{|\bar{x}| + m - 1} e^{-c_{\varepsilon'} \alpha n}$  
where the constant $c_{\varepsilon'}$ depends only on $\varepsilon'$.
By the induction hypothesis again,
the probability that 
\begin{itemize}
\item[(b)] $\mcA \in \mbW^{\mcA'}$ is 
$(p, p_m, \beta/(1 + \varepsilon')^2)$-saturated and 
$(p, p_m, \beta(1 + \varepsilon')^2)$-unsaturated 
\end{itemize}
is at least
$1 \ - \ 2 n^{|\bar{x}| + m} e^{-c_{\varepsilon'} \beta n}$ where $c_{\varepsilon'}$ is the same constant as above
(since it depends only on $\varepsilon'$).
It is straightforward to check that if $\mcA \in \mbW^{\mcA'}$ satisfies both~(a) and~(b) then 
$\mcA$ is $(p, q, \gamma/(1 + \varepsilon')^{2(m+1)})$-saturated and
$(p, q, \gamma (1 + \varepsilon')^{2(m+1)})$-unsaturated.
Since $\gamma = \alpha\beta \leq \min\{\alpha, \beta\}$ it follows that the probability that $\mcA \in \mbW^{\mcA'}$ is 
$(p, q, \gamma/(1 + \varepsilon')^{2(m+1)})$-saturated and
$(p, q, \gamma (1 + \varepsilon')^{2(m+1)})$-unsaturated
is at least $1 \ - \ 2^{m+1} n^{|\bar{x}| + m} e^{-c_{\varepsilon'} \gamma n}$.
\hfill $\square$

\begin{defin}\label{definition of Y-n}{\rm
For every $n$, let $\mbY_n$ be the set of all $\mcA \in \mbW^{\mbY'_n}$ such that
whenever $p(\bar{x}, \bar{y})$ and $q(\bar{x})$ are complete atomic $\sigma$-types with
$|\bar{x}\bar{y}| \leq k$, $\dim_{\bar{y}}(p) = |\bar{y}|$, $q \subseteq p$ and
$\gamma = \msfP(p \ | \ q)$, then 
$\mcA$ is $(p, q, \gamma/(1 + \varepsilon')^{2|\bar{y}|})$-saturated and
$(p, q, \gamma (1 + \varepsilon')^{2|\bar{y}|})$-unsaturated.
}\end{defin}

\noindent
The following corollary follows directly from the definition of $\mbY_n$ and 
Lemma~\ref{the number of realizations in A' for y-bar}.

\begin{cor}\label{all members of Y-n are sufficiently saturated}
Let $p(\bar{x}, \bar{y})$ and $q(\bar{x})$ are complete atomic $\sigma$-types such that
$|\bar{x}\bar{y}| \leq k$, $d = \dim_{\bar{y}}(p) > 0$, $q \subseteq p$ and
$\gamma = \msfP(p \ | \ q)$.
For every $n$, every $\mcA \in \mbY_n$ is
$(p, q, \gamma/(1 + \varepsilon')^{2d})$-saturated and
$(p, q, \gamma (1 + \varepsilon')^{2d})$-unsaturated.
\end{cor}

\begin{lem}\label{Y-n has large probability}
There is a constant $c > 0$ such that for all sufficiently large  $n$,
$\mbbP_n\big(\mbY_n\big) \geq \big(1 - e^{-c n}\big) \big(1 - \delta'(n)\big)$.
\end{lem}

\noindent
{\bf Proof.}
There are, up to changing variables, only finitely many atomic $\sigma$-types $p(\bar{x})$ such that $|\bar{x}| \leq k$.
It follows from 
Lemma~\ref{the number of realizations in A' for y-bar} 
that there is a constant $c > 0$ such that for all large enough $n$ and
all $\mcA' \in \mbY'_n$,
\[
\mbbP_n^{\mcA'}\big(\mbY_n \cap \mbW^{\mcA'}\big) \geq 1 - e^{-c n}.
\]
Note that $\mbbP_n(\mbY_n) = \mbbP_n\big(\mbY_n \ | \ \mbW^{\mbY'_n}\big) \mbbP_n\big(\mbW^{\mbY'_n}\big)$.
By Lemma~\ref{P and P' agree when constraints speak only about W'},
$\mbbP_n(\mbW^{\mbY'_n}) = \mbbP'_n(\mbY'_n)$ and
by Lemma~\ref{P(X given Y) equals P_Y(X cut Y)} we have
$\mbbP_n(\mbY'_n \ | \ \mbW^{\mbY'_n}) = \mbbP^{\mbY'_n}(\mbY_n \cap \mbW^{\mbY'_n})$.
Hence
$\mbbP_n(\mbY_n) = \mbbP^{\mbY'_n}(\mbY_n \cap \mbW^{\mbY'_n}) \mbbP'_n\big(\mbY'_n\big)$.
Then, reasoning similarly as in the proof of Lemma~\ref{P(p) is a product} (using~(\ref{first basic equality about restricted distributions})),  
we get
\begin{align*}
&\mbbP^{\mbY'_n}\big(\mbY_n \cap \mbW^{\mbY'_n}\big) \ = \
\sum_{\mcA' \in \mbY'_n} \mbbP^{\mbY'_n}\big(\mbY_n \cap \mbW^{\mcA'}\big) \ = \
\sum_{\mcA' \in \mbY'_n}\sum_{\mcA \in \mbY_n \cap \mbW^{\mcA'}}\mbbP^{\mbY'_n}(\mcA) = \\
&\sum_{\mcA' \in \mbY'_n}\sum_{\mcA \in \mbY_n \cap \mbW^{\mcA'}}
\frac{\mbbP'_n(\mcA')}{\mbbP'_n(\mbY'_n)} \mbbP^{\mcA'}(\mcA) \ = \
\sum_{\mcA' \in \mbY'_n}\frac{\mbbP'_n(\mcA')}{\mbbP'_n(\mbY'_n)}\sum_{\mcA \in \mbY_n \cap \mbW^{\mcA'}}
\mbbP^{\mcA'}(\mcA) = \\
&\sum_{\mcA' \in \mbY'_n}\frac{\mbbP'_n(\mcA')}{\mbbP'_n(\mbY'_n)}\mbbP^{\mcA'}\big(\mbY_n \cap \mbW^{\mcA'}\big) 
\ \geq \ 
\sum_{\mcA' \in \mbY'_n}\frac{\mbbP'_n(\mcA')}{\mbbP'_n(\mbY'_n)}
\big(1 - e^{-c n}\big) = \\
&\big(1 - e^{-c n}\big).
\end{align*}
Using part~(2) of Assumption~\ref{inductive assumptions} we know get
\[
\mbbP_n(\mbY_n) \ = \ \mbbP^{\mbY'_n}\big(\mbY_n \cap \mbW^{\mbY'_n}\big) \mbbP'_n(\mbY'_n) \ \geq \
\big(1 - e^{-c n}\big) \big(1 - \delta'(n)\big).
\]
\hfill $\square$

\begin{defin}\label{definition of critical number}{\rm
A real number is called  {\em critical} if it is {\em $m$-critical} for some positive integer $m$.
We say that a real number $\alpha$ is {\em $m$-critical} if at least one of the following holds:
\begin{itemize}
\item[(a)] There are a complete atomic $\sigma$-type $q(\bar{x})$, 
distinct complete atomic $\sigma$-types $p_1(\bar{x}, \bar{y})$, $\ldots$, $p_l(\bar{x}, \bar{y})$
and a number $1 \leq l' \leq l$
such that $|\bar{x}\bar{y}| \leq m$, $q \subseteq p_i$ for all $1 \leq i \leq l$ and 

\[
\alpha = \frac{\sum_{i=1}^{l'} \msfP(p_i \ | \ q)}{\sum_{i=1}^l \msfP(p_i \ | \ q)}.
\]

\item[(b)] $\alpha = l'/l$ where $0 \leq l' \leq l$ are integers and $l$ is, for any choice of distinct variables
$x_1, \ldots, x_m$, less or equal to 
the number of pairs $(p(x_1, \ldots, x_{m'}), q(x_1, \ldots, x_d))$ where $d < m' \leq m$,
$p$ and $q$ are complete atomic $\sigma$-types such that $q \subseteq p$ and $\dim_{(x_d, \ldots, x_{m'})}(p) = 0$.
\end{itemize}
}\end{defin}

\noindent
From the definition it follows that (for every $m \in \mbbN$) there are only finitely many $m$-critical numbers.
It also follows (from part~(b)) that, for every $m$, $0$ and $1$ are $m$-critical.
In part~(a) above we allow $\bar{x}$ to be empty in which case the type $q(\bar{x})$ is omitted and $\msfP(p_i \ | \ q)$
is replaced by $\msfP(p_i)$.

\begin{defin}\label{definition of noncritical formula}{\rm
Let $\varphi(\bar{x}) \in CPL(\sigma)$ and let $l = |\bar{x}| + \text{qr}(\varphi)$.\\
(i) We call $\varphi(\bar{x})$ {\em noncritical} if the following holds:
\begin{itemize}
\item[] If 
\begin{align*}
&\Big( r + \| \psi(\bar{z}, \bar{y}) \ |  \ \theta(\bar{z}, \bar{y}) \|_{\bar{y}} \ \geq \ 
\| \psi^*(\bar{z}, \bar{y}) \ |  \ \theta^*(\bar{z}, \bar{y}) \|_{\bar{y}} \Big) \quad \text{ or} \\
&\Big( \| \psi(\bar{z}, \bar{y}) \ |  \ \theta(\bar{z}, \bar{y}) \|_{\bar{y}} \ \geq \ 
\| \psi^*(\bar{z}, \bar{y}) \ |  \ \theta^*(\bar{z}, \bar{y}) \|_{\bar{y}}  + r \Big)
\end{align*}
is a subformula of $\varphi(\bar{x})$ (where $\psi, \theta, \psi^*$ and $\theta^*$ denote formulas in $CPL(\sigma)$ and
$\bar{z}$ and $\bar{y}$ may have variables in common with $\bar{x}$)
then, for all $l$-critical numbers $\alpha$ and $\beta$, $r \neq \alpha - \beta$.
\end{itemize}
(ii) Let $\varepsilon > 0$. We say that $\varphi(\bar{x})$ is {\em $\varepsilon$-noncritical} if 
\begin{itemize}
\item $\varphi(\bar{x})$ is noncritical and

\item whenever $r$ appears in a subformula as in part~(i) and $\alpha$ and $\beta$ are $l$-critical numbers, 
then the following implications hold:
\begin{itemize}
\item[] If $r + \alpha > \beta$ then $r + \alpha/(1 + 2\varepsilon)^2 > \beta (1 + 2\varepsilon)^2$, and \\
if $\alpha > \beta + r$ then $\alpha/(1 + 2\varepsilon)^2 > \beta (1 + 2\varepsilon)^2 + r$.
\end{itemize}
\end{itemize}
}\end{defin}

\noindent
Since, for every $l \in \mbbN$, there are only finitely many $l$-critical numbers it follows that for every
noncritical $\varphi(\bar{x}) \in CPL(\sigma)$, if one just chooses $\varepsilon > 0$ sufficiently small, then
$\varphi(\bar{x})$ is $\varepsilon$-noncritical.
Definition~\ref{definition of Y-n} and Lemma~\ref{Y-n has large probability} motivate the next definition.

\begin{defin}\label{definition of epsilon}{\rm
Let $\varepsilon > 0$ be such that $1 + \varepsilon = (1 + \varepsilon')^{2k}$.
}\end{defin}

\noindent
It follows from Definition~\ref{definition of epsilon}  and
Lemma~\ref{all members of Y-n are sufficiently saturated} that
if $p(\bar{x}, \bar{y})$ and $q(\bar{x})$ are complete atomic $\sigma$-types such that
$|\bar{x}\bar{y}| \leq k$, $d = \dim_{\bar{y}}(p) > 0$, $q \subseteq p$, $\msfP(q) > 0$, and
$\gamma = \msfP(p \ | \ q)$, then
for every $n$, every $\mcA \in \mbY_n$ is
$(p, q, \gamma/(1 + \varepsilon))$-saturated and
$(p, q, \gamma (1 + \varepsilon))$-unsaturated.
By an analogous argument as in Remark~\ref{special cases of inductive assumptions}~(iii), it now follows 
that if $p(\bar{x})$ is a complete atomic $\sigma$-type such that $|\bar{x}| \leq k$ and 
$\msfP(p) = 0$, then for all sufficiently large $n$, $p$ is not realized
in any member of $\mbY_n$.

In the proof of the proposition below we will sometimes abuse notation by treating an atomic type 
$p(\bar{x})$ as the formula obtained by taking the
conjunction of all formulas in $p(\bar{x})$.
So when writing, for example, `$\bigvee_{i=1}^m\bigvee_{j=1}^{m_i} p_{i, j}(\bar{x}, y)$' 
in the proof below we view $p_{i, j}(\bar{x}, y)$ in this expression as the conjunction of all formulas in the
complete atomic type $p_{i, j}(\bar{x}, y)$.

\begin{prop}\label{quantifier elimination in a structure} {\bf (Elimination of quantifiers)}
Suppose that $\varphi(\bar{x}) \in CPL(\sigma)$ is $\varepsilon$-noncritical and $|\bar{x}| + \text{qr}(\varphi) \leq k$.
Then there is a quantifier-free formula $\varphi^*(\bar{x})$ such that 
for all sufficiently large $n$ and every $\mcA \in \mbY_n$,
$\mcA \models \forall \bar{x} (\varphi(\bar{x}) \leftrightarrow \varphi^*(\bar{x}))$.
\end{prop}

\noindent
{\bf Proof.}
Let an $\varepsilon$-noncritical $\varphi(\bar{x}) \in CPL(\sigma)$ be given with $|\bar{x}| + \text{qr}(\varphi) \leq k$.
We will assume that $\bar{x}$ is nonempty (i.e. that $\varphi$ has free variables). 
In Remark~\ref{the case when bar-x is empty} it is indicated which changes we need to make in
the simpler case when $\varphi$ has no free variable.
The proof proceeds by induction on quantifier-rank.
Suppose that $\text{qr}(\varphi) > 0$ since otherwise we can just let $\varphi^*$ be $\varphi$ and then we are done.
If for all sufficiently large $n$, for all $\mcA \in \mbY_n$ and for all $\bar{a} \in [n]^{|\bar{x}|}$ we have 
$\mcA \not\models \varphi(\bar{a})$ then we can let $\varphi^*(\bar{x})$ be the formula $x_1 \neq x_1$ and then
$\mcA \models \forall \bar{x} (\varphi(\bar{x}) \leftrightarrow \varphi^*(\bar{x}))$ for all sufficiently large $n$
and all $\mcA \in \mbY_n$.
So from now on we assume that, for arbitrarily large $n$, there are $\mcA \in \mbY_n$ and $\bar{a}$ such that
$\mcA \models \varphi(\bar{a})$.

Suppose that $\varphi(\bar{x})$ is $\exists y \psi(\bar{x}, y)$ for some $\psi(\bar{x}, y)$.
Then we have $|\bar{x}y| + \text{qr}(\psi) \leq k$ and $\text{qr}(\psi) < \text{qr}(\varphi)$ so, 
by the induction hypothesis, we may assume that
$\psi(\bar{x}, y)$ is quantifier-free.
By assumption there are $n$, $\mcA \in \mbY_n$, $\bar{a}$ and $b$ such that $\mcA \models \psi(\bar{a}, b)$.
Then there are $m \geq 1$, different complete atomic $\sigma$-types $q_i(\bar{x})$, $i = 1, \ldots, m$, and, for each $i$,
$m_i \geq 1$ and different complete atomic $\sigma$-types $p_{i, j}(\bar{x}, y)$, $j = 1, \ldots, m_i$, such that
$q_i \subseteq p_{i, j}$ for all $j$ and
$\psi(\bar{x}, y)$ is equivalent to $\bigvee_{i=1}^m\bigvee_{j=1}^{m_i} p_{i, j}(\bar{x}, y)$.
If, for some $i$, $\msfP(q_i(\bar{x})) = 0$, then
$q_i$ is not realized in any $\mcA \in \mbY_n$ (for large enough $n$) and can be removed.
So we may assume that all $\msfP(q_i) > 0$ for all $i$.
If, for some $i$ and $j$, $\msfP(p_{i, j} \ | \ q) = 0$ then $\msfP(p_{i, j}) = 0$ so
$p_{i, j}$ is not realized in any $\mcA \in \mbY_n$ for large enough $n$.
So we may also assume that $\msfP(p_{i, j} \ | \ q_i) > 0$ for all $i$ and $j$.
If $\dim_y(p_{i, j}) = 1$ then, by the definitions of $\mbY_n$ and $\varepsilon$, it follows that for all
sufficiently large $n$ and all $\mcA \in \mbY_n$, if $\mcA \models q_i(\bar{a})$ then $\mcA \models \exists y p_{i, j}(\bar{a}, y)$.
If $\dim_y(p_{i, j}) = 0$ then, for all $n$ and all $\mcA \in \mbW_n$, if $\mcA \models q_i(\bar{a})$
then $\mcA \models p_{i, j}(\bar{a}, b)$ for some $b \in \rng(\bar{a})$.
It follows that for all sufficiently large $n$ and all $\mcA \in \mbY_n$,
$\mcA \models \forall \bar{x} 
\big( \exists y \psi(\bar{x}, y) \leftrightarrow \bigvee_{i=1}^m q_i(\bar{x})\big)$.

Now we consider the case when $\varphi(\bar{x})$ has the form
\begin{align}\label{the formula}
&\Big( r + \| \psi(\bar{x}, \bar{y}) \ |  \ \theta(\bar{x}, \bar{y}) \|_{\bar{y}} \ \geq \ 
\| \psi^*(\bar{x}, \bar{y}) \ |  \ \theta^*(\bar{x}, \bar{y}) \|_{\bar{y}} \Big)
\ \ \text{ or} \\
&\Big( \| \psi(\bar{x}, \bar{y}) \ |  \ \theta(\bar{x}, \bar{y}) \|_{\bar{y}} \ \geq \ 
\| \psi^*(\bar{x}, \bar{y}) \ |  \ \theta^*(\bar{x}, \bar{y}) \|_{\bar{y}} + r \Big). \label{the second possibility}
\end{align}
Since the second case~(\ref{the second possibility}) is treated by straightforward variations of the arguments for taking care of the first 
case~(\ref{the formula})
we only consider the first case~(\ref{the formula}).
Observe that $|\bar{x}\bar{y}| + \text{qr}(\psi) \leq k$ 
(because $\text{qr}(\varphi) = |\bar{y}| + \max\{\text{qr}(\psi), \text{qr}(\theta), \text{qr}(\psi^*), \text{qr}(\theta^*)\}$)
and similarly for $\theta$, $\psi^*$ and $\theta^*$.
Since all the formulas $\psi$, $\theta$, $\psi^*$ and $\theta^*$ have smaller quantifier-rank than $\varphi$ we may,
by the induction hypothesis,
assume that $\psi(\bar{x}, \bar{y}), \theta(\bar{x}, \bar{y}), \psi^*(\bar{x}, \bar{y})$ and $\theta^*(\bar{x}, \bar{y})$
are quantifier-free formulas.

If $\theta(\bar{x}, \bar{y})$ or $\theta^*(\bar{x}, \bar{y})$ is unsatisfiable, then, by the provided semantics,
we have $\mcA \not\models \varphi(\bar{a})$ for every $\sigma$-structure $\mcA$ and every sequence of elements
$\bar{a}$ from the domain of $\mcA$. In this case $\varphi(\bar{x})$ is equivalent to any contradictory quantifier-free
formula with free variables among $\bar{x}$, for example the formula $x_1 \neq x_1$.
So from now on we assume that $\theta(\bar{x}, \bar{y})$ and $\theta^*(\bar{x}, \bar{y})$ are satisfiable.

Until further notice, assume also that $\psi(\bar{x}, \bar{y}) \wedge \theta(\bar{x}, \bar{y})$ and 
$\psi^*(\bar{x}, \bar{y}) \wedge \theta^*(\bar{x}, \bar{y})$ are satisfiable.
Then there are distinct complete atomic $\sigma$-types $q_i(\bar{x})$, $p_{i, j}(\bar{x}, \bar{y})$, for $i = 1, \ldots, m$ and $j = 1, \ldots, m_i$,
and distinct complete atomic $\sigma$-types $t_i(\bar{x})$, $s_{i, j}(\bar{x}, \bar{y})$, for $i = 1, \ldots, l$ and $j = 1, \ldots, l_i$, such that
the following conditions hold:
\begin{itemize}
\item $q_i(\bar{x}) \subseteq p_{i,j}(\bar{x}, \bar{y})$ for all $i = 1, \ldots, m$ and all $j = 1, \ldots, m_i$.

\item $\psi(\bar{x}, \bar{y}) \wedge \theta(\bar{x}, \bar{y})$ is equivalent to 
$\bigvee_{i=1}^m\bigvee_{j=1}^{m_i}p_{i,j}(\bar{x}, \bar{y})$.

\item $t_i(\bar{x}) \subseteq s_{i,j}(\bar{x}, \bar{y})$ for all $i = 1, \ldots, l$ and all $j = 1, \ldots, l_i$.

\item $\theta(\bar{x}, \bar{y})$ is equivalent to 
$\bigvee_{i=1}^l\bigvee_{j=1}^{l_i}s_{i,j}(\bar{x}, \bar{y})$.
\end{itemize}
Since $\bigvee_{i=1}^l\bigvee_{j=1}^{l_i}s_{i,j}(\bar{x}, \bar{y})$ is a consequence 
of $\bigvee_{i=1}^m\bigvee_{j=1}^{m_i}p_{i,j}(\bar{x}, \bar{y})$ it follows that 
$m \leq l$ and $m_i \leq l_i$ for all $i \leq m$.
Moreover, for every $i \leq m$ there is $i'$ such that $q_i = t_{i'}$, and for all $i \leq m$ and all $j \leq m_i$ there are $i', j'$ such that $p_{i,j} = s_{i',j'}$.
Therefore we may assume in addition (by reordering if necessary) that
\begin{equation}\label{some identities between the types}
\text{$q_i = t_i$ for all $i \leq m$ and $p_{i, j} = s_{i, j}$ for all $i \leq m$ and all $j \leq m_i$.}
\end{equation}

\noindent
For the same reasons as in the previous case we may assume that all of $\msfP(q_i)$, $\msfP(p_{i, j})$, $\msfP(t_i)$ and $\msfP(s_{i, j})$
are positive for all $i$ and $j$.
Next we define
\begin{align*}
&d_{i, j} = \dim_{\bar{y}}(p_{i, j})  \ \text{ for all $i = 1, \ldots, m$ and $j = 1, \ldots, m_i$}, \\
&e_{i, j} = \dim_{\bar{y}}(s_{i, j})  \ \text{ for all $i = 1, \ldots, l$ and $j = 1, \ldots, l_i$}, \\
&d_i = \max\{d_{i, 1}, \ldots, d_{i, m_i}\} \ \text{ for all $i = 1, \ldots, m$}, \\
&e_i = \max\{e_{i, 1}, \ldots, e_{i, l_i}\} \ \text{ for all $i = 1, \ldots, l$}, \\
&\alpha_{i,j} = \msfP(p_{i,j}(\bar{x}, \bar{y}) \ | \ q_i(\bar{x})) \ \text{ for all $i = 1, \ldots, m$ and $j = 1, \ldots, m_i$}, \\
&\alpha_i = \text{ the sum of all $\alpha_{i, j}$ such that $d_{i, j} = d_i$}, \\
&\beta_{i,j} = \msfP(s_{i,j}(\bar{x}, \bar{y}) \ | \ t_i(\bar{x})) \ \text{ for all $i = 1, \ldots, l$ and $j = 1, \ldots, l_i$}, \\
&\beta_i = \text{ the sum of all $\beta_{i, j}$ such that $e_{i, j} = e_i$.}
\end{align*}
It follows that for all $i = 1, \ldots, m$ we have $d_i \leq e_i$ and $\alpha_i \leq \beta_i$.

\begin{defin}\label{definition of gamma}{\rm
For all $i = 1, \ldots, l$ we define a number $\gamma_i$ as follows:
\begin{enumerate}
\item If $i \leq m$ and $d_i = e_i > 0$ then we define $\gamma_i = \alpha_i / \beta_i$.

\item If $i \leq m$ and $d_i = e_i = 0$ then we define $\gamma_i = m_i / l_i$.

\item If $i \leq m$ and $d_i < e_i$ then we define $\gamma_i = 0$.

\item If $m < i \leq l$ then we define $\gamma_i = 0$.
\end{enumerate}
}\end{defin}

\noindent
Now we can reason in exactly the same way with regard to the formulas $\psi^*(\bar{x}, \bar{y})$ and $\theta^*(\bar{x}, \bar{y})$.
So there are numbers $m^*, l^*, m_i^*$ and $l_i^*$ and complete atomic $\sigma$-types
$q_i^*(\bar{x})$ for $i = 1, \ldots, m^*$, $p_{i, j}^*(\bar{x}, \bar{y})$ for $i \leq m^*$ and $j = 1, \ldots, m_i^*$, 
$t_i^*(\bar{x})$ for $i = 1, \ldots, l^*$ and $s_{i, j}^*(\bar{x}, \bar{y})$ for $i \leq l^*$ and $j = 1, \ldots, l_i^*$
such that all which has been said about $\psi$, $\theta$, $q_i$, $p_{i, j}$, $t_i$ and $s_{i,j}$ holds if these formulas and
types
are replaced by $\psi^*$, $\theta^*$, $q_i^*$, $p_{i, j}^*$ etcetera, and the numbers $m$, $l$, $m_i$, $l_i$ are
replaced by $m^*, l^*, m_i^*$ and $l_i^*$. 
Moreover, we define numbers $d_{i, j}^*$, $e_{i, j}^*$, $d_i^*$, $e_i^*$, 
$\alpha_{i, j}^*$, $\alpha_i^*$, $\beta_{i, j}^*$, $\beta_i^*$ and $\gamma_i^*$
in the same way as above, using the types $q_i^*$, $p_{i, j}^*$, $t_i^*$ and $s_{i, j}^*$ 
instead of $q_i$, $p_{i, j}$, $t_i$ and $s_{i, j}$.

So far we have assumed that $\psi(\bar{x}, \bar{y}) \wedge \theta(\bar{x}, \bar{y})$ and 
$\psi^*(\bar{x}, \bar{y}) \wedge \theta^*(\bar{x}, \bar{y})$ are satisfiable.
If $\psi(\bar{x}, \bar{y}) \wedge \theta(\bar{x}, \bar{y})$ is not satisfiable, then
we let $m = 0$ and we view the disjunction 
$\bigvee_{i=1}^m\bigvee_{j=1}^{m_i}p_{i,j}(\bar{x}, \bar{y})$
as ``empty'' and hence always false.
In this case we always have $i > m$ so it follows that $\gamma_i = 0$ for all $i = 1, \ldots, l$.
Similar conventions apply if $\psi^*(\bar{x}, \bar{y}) \wedge \theta^*(\bar{x}, \bar{y})$ is not satisfiable.
With these conventions the case when any one of the mentioned formulas is unsatisfiable is taken care of by
the rest of the proof.

\begin{lem}\label{upper and lower bounds on p-i-j divided by s-i-j}
Let $i \in \{1, \ldots, m\}$.\\
(a) For all suffiently large $n$ and all $\mcA \in \mbY_n$,
\[
\frac{\gamma_i}{(1 + 2\varepsilon)^2} \ \leq \ \frac{\Big| \bigcup_{j=1}^{m_i} p_{i,j}(\bar{a}, \mcA)\Big|}
{\Big| \bigcup_{j=1}^{l_i} s_{i,j}(\bar{a}, \mcA)\Big|}.
\]
(b) If $d_i = e_i$ then, for all sufficiently large $n$ and all $\mcA \in \mbY_n$,
\[
\frac{\Big| \bigcup_{j=1}^{m_i} p_{i,j}(\bar{a}, \mcA)\Big|}
{\Big| \bigcup_{j=1}^{l_i} s_{i,j}(\bar{a}, \mcA)\Big|} \ \leq \ (1 + 2\varepsilon)^2 \gamma_i.
\]
(c) If $d_i < e_i$ then, for all sufficiently large $n$ and all $\mcA \in \mbY_n$,
\begin{equation*}
\frac{\Big| \bigcup_{j=1}^{m_i} p_{i,j}(\bar{a}, \mcA)\Big|}{\Big| \bigcup_{j=1}^{l_i} s_{i,j}(\bar{a}, \mcA)\Big|} 
\ \leq \ \frac{C}{n}
\end{equation*}
where the constant $C > 0$ depends only on the types $p_{i, j}$ and $s_{i, j}$.\\
(d) Parts (a), (b) and (c) hold if $m, m_i, l_i, \gamma_i$, $d_i$, $e_i$, $p_{i, j}$ and $s_{i, j}$ are replaced
by $m^*$, $m_i^*$, $l_i^*$, $\gamma_i^*$, $d_i^*$, $e_i^*$, $p_{i, j}^*$ and $s_{i, j}^*$, respectively.
\end{lem}

\noindent
{\bf Proof.}
We split the argument into cases corresponding to the three first cases of Definition~\ref{definition of gamma}.
Let $\mcA \in \mbY_n$.

{\em First suppose that $d_i = e_i > 0$} and hence $\gamma_i = \alpha_i / \beta_i$.
Since $\mcA$ is assumed to be $(p_{i,j}, q_i, (1+\varepsilon)\alpha_{i,j})$-unsaturated if $d_{i,j} > 0$
it follows that 
\[
\big| p_{i, j}(\bar{a}, \mcA) \big| \leq (1 + \varepsilon) \alpha_{i, j} n^{d_{i,j}} \ \text{ if $d_{i, j} > 0$.}
\]
If $d_{i,j} = 0$ then $\big|p_{i, j}(\bar{a}, \mcA)\big| = 1$ and each member of the unique tuple realizing $p_{i, j}(\bar{a}, \bar{y})$
belongs to $\bar{a}$.
It follows that 
\begin{equation}\label{upper bound of p}
\bigg|\bigcup_{j=1}^{m_i} p_{i,j}(\bar{a}, \mcA)\bigg| 
\ \leq \ (1 + 2\varepsilon)\alpha_i n^{d_i} \text{ for all sufficiently large $n$.}
\end{equation}
By similiar reasoning (and since we assume $d_i = e_i$) we get
\begin{equation}\label{upper bound of s}
\bigg|\bigcup_{j=1}^{l_i} s_{i,j}(\bar{a}, \mcA)\bigg| \  \leq \ (1 + 2\varepsilon)\beta_i n^{d_i}.
\end{equation}

Since $\mcA$ is assumed to be $(p_{i,j}, q_i, \alpha_{i,j}/(1+\varepsilon))$-saturated if $d_{i, j} > 0$
and $(s_{i,j}, t_i, \beta_{i,j}/(1+\varepsilon))$-saturated if $e_{i, j} > 0$
it follows that $\big|p_{i, j}(\bar{a}, \mcA)\big| \geq \alpha_{i, j}n^{d_{i, j}}/(1 + \varepsilon)$ if $d_{i, j} > 0$ and
 $\big|s_{i, j}(\bar{a}, \mcA)\big| \geq \beta_{i, j}n^{e_{i, j}}/(1 + \varepsilon)$ if $e_{i, j} > 0$.
This (and $d_i = e_i$) implies that
\begin{equation}\label{lower bounds of p and s}
\bigg|\bigcup_{j=1}^{m_i} p_{i,j}(\bar{a}, \mcA)\bigg| \  \geq \ \frac{\alpha_i n^{d_i}}{1+\varepsilon} \quad \text{and} \quad
\bigg|\bigcup_{j=1}^{l_i} s_{i,j}(\bar{a}, \mcA)\bigg| \ \geq \ \frac{\beta_i n^{d_i}}{1+\varepsilon}.
\end{equation}
From~(\ref{upper bound of p}), (\ref{upper bound of s}) and~(\ref{lower bounds of p and s}) we get
\begin{equation}\label{upper and lower bounds of p/s}
\frac{\gamma_i}{(1 + 2\varepsilon)^2} \ = \ \frac{\alpha_i}{(1 + 2\varepsilon)^2 \beta_i} \ \leq \
\frac{\Big| \bigcup_{j=1}^{m_i} p_{i,j}(\bar{a}, \mcA)\Big|}
{\Big| \bigcup_{j=1}^{l_i} s_{i,j}(\bar{a}, \mcA)\Big|} \ \leq \ (1 + 2\varepsilon)^2\frac{\alpha_i}{\beta_i}
\ = \ (1 + 2\varepsilon)^2 \gamma_i .
\end{equation}

{\em Now suppose that $d_i = e_i = 0$.} Then $\gamma_i = m_i / l_i$.
Also, each $p_{i, j}(\bar{a}, \bar{y})$ and each $s_{i, j}(\bar{a}, \bar{y})$  has a unique realization in $\mcA$.
Since we assume that $p_{i, j} \neq p_{i, j'}$ if $j \neq j'$ and $s_{i, j} \neq s_{i, j'}$ if $j \neq j'$ we get
\[
\gamma_i = 
\frac{m_i}{l_i} = \frac{\Big| \bigcup_{j=1}^{m_i} p_{i,j}(\bar{a}, \mcA)\Big|}
{\Big| \bigcup_{j=1}^{l_i} s_{i,j}(\bar{a}, \mcA)\Big|},
\]
and now the inequalities of~(a) and~(b) follow trivially.

{\em Next, suppose that $d_i < e_i$.}
Then $\gamma_i = 0$.
By similar reasoning as before,
\[
0 \  < \ \frac{\beta_i n^{e_i}}{1 + 2\varepsilon} \ \leq \ \bigg| \bigcup_{j=1}^{l_i} s_{i, j}(\bar{a}, \mcA) \bigg|
\ \text{ for sufficiently large $n$}.
\]
It follows that 
\[
\frac{\gamma_i}{(1 + 2\varepsilon)^2} \ = \ 0 \leq \frac{\Big| \bigcup_{j=1}^{m_i} p_{i,j}(\bar{a}, \mcA)\Big|}
{\Big| \bigcup_{j=1}^{l_i} s_{i,j}(\bar{a}, \mcA)\Big|}.
\]
Since $e_i > 0$ we can argue as we did to get~(\ref{lower bounds of p and s}), so we have
\[
\bigg|\bigcup_{j=1}^{l_i} s_{i,j}(\bar{a}, \mcA)\bigg| \ \geq \ \frac{\beta_i n^{e_i}}{1+\varepsilon}.
\]
Depending on whether $d_i > 0$ or $d_i = 0$ we get, by arguing as in previous cases,
\[
\bigg|\bigcup_{j=1}^{m_i} p_{i,j}(\bar{a}, \mcA)\bigg| 
\ \leq \ (1 + 2\varepsilon)\alpha_i n^{d_i} 
\quad  \text{ or } \quad
\bigg|\bigcup_{j=1}^{m_i} p_{i,j}(\bar{a}, \mcA)\bigg|  = m_i.
\]
Since $d_i < e_i$ we get, in either case,
\begin{equation*}
\frac{\Big| \bigcup_{j=1}^{m_i} p_{i,j}(\bar{a}, \mcA)\Big|}{\Big| \bigcup_{j=1}^{l_i} s_{i,j}(\bar{a}, \mcA)\Big|} 
\ \leq \ \frac{C}{n} \ \text{ for sufficiently large $n$}
\end{equation*}
where $C > 0$ is a constant that depends only on the types $p_{i, j}$ and $s_{i, j}$.

The proof of part~(d) is, of course, exactly the same (besides the relevant replacements of symbols).
\hfill $\square$

\begin{defin}\label{definition of I}{\rm
Let $I$ be the set of all $i \in \{1, \ldots, l\}$ such that there exists some $i' \in \{1, \ldots, l^*\}$
such that $t_i = t_{i'}^*$ and $r + \gamma_i \geq \gamma_{i'}^*$.
}\end{defin}

\begin{rem}\label{remark on computing I} {\bf (The computational problem of finding $I$)} {\rm
The number $\alpha_{i, j}$ is obtained from numbers given by 
assumptions~\ref{basic assumptions in the proofs}
and~\ref{inductive assumptions}
and applying a number of arithmetic operations which is linear in $|p_{i, j}|$.
It follows that the number of arithmetic operations needed to compute $\alpha_i$ is linear in $\sum_{j=1}^{m_i} |p_{i, j}|$,
where by an arithmetic operation I mean addition, multiplication or division.
The case is similar for $\beta_i$, $\alpha^*_i$ and $\beta^*_i$ .
The number of comparisons of literals needed to check if $t_i = t_{i'}^*$ is $|t_i|$ if we assume that
we use some uniform way of listing the literals in complete atomic $\sigma$-types.
So to decide if $i \in I$ we need to perform a number of arithmetic operations, comparisons of literals and comparisons of numbers
which is linear in 
\[
\sum_{j=1}^{m_i} |p_{i, j}| + \sum_{j=1}^{l_i} |s_{i, j}| + \sum_{j=1}^{m^*_i} |p^*_{i, j}| + \sum_{j=1}^{l^*_i} |s^*_{i, j}|.
\]
Consequently the number of arithmetic operations, comparisons of literals and comparisons of numbers that are needed to create $I$
is linear in 
\[
\sum_{i=1}^m\sum_{j=1}^{m_i} |p_{i, j}| + \sum_{i=1}^{l}\sum_{j=1}^{l_i} |s_{i, j}| + 
\sum_{i=1}^{m^*}\sum_{j=1}^{m^*_i} |p^*_{i, j}| + \sum_{i=1}^{l^*}\sum_{j=1}^{l^*_i} |s^*_{i, j}|.
\]

}\end{rem}

\noindent
Lemmas~\ref{the case when I is nonempty}
and~\ref{the case when I is empty}
below show that $\varphi(\bar{x})$ 
is equivalent, in every $\mcA \in \mbW_n$ for all large enough $n$,
to a quantifier-free formula which depends only on $\varphi(\bar{x})$ and the 
lifted Bayesian network $\mfG$.
As noted after
Definition~\ref{definition of critical number},
0 is a $\zeta$-critical number for every $\zeta$, so $r > 0$ (since $\varphi$ is noncritical).
Observe that it follows from
Definitions~\ref{definition of critical number}
and~\ref{definition of gamma}
that $\gamma_i$ and $\gamma_i^*$ are 
$(|\bar{x}| + \mr{qr}(\varphi))$-critical numbers for all $i$.

\begin{lem}\label{the case when I is nonempty}
Suppose that $I \neq \es$. 
Then for all sufficiently large $n$, all $\mcA \in \mbY_n$ and all $\bar{a} \in [n]^{|\bar{x}|}$,
\[
\mcA \models
\Big( r + \| \psi(\bar{a}, \bar{y}) \ |  \ \theta(\bar{a}, \bar{y}) \|_{\bar{y}} 
\ \geq \ 
\| \psi^*(\bar{a}, \bar{y}) \ |  \ \theta^*(\bar{a}, \bar{y}) \|_{\bar{y}} \Big)
\]
if and only if $\mcA \models \bigvee_{i \in I} t_i(\bar{a})$.
\end{lem}

\noindent
{\bf Proof.}
Suppose that 
\[
\mcA \models
\Big( r + \| \psi(\bar{a}, \bar{y}) \ |  \ \theta(\bar{a}, \bar{y}) \|_{\bar{y}} 
\ \geq \ 
\| \psi^*(\bar{a}, \bar{y}) \ |  \ \theta^*(\bar{a}, \bar{y}) \|_{\bar{y}} \Big).
\]
Then both $\| \psi(\bar{a}, \bar{y}) \ |  \ \theta(\bar{a}, \bar{y}) \|_{\bar{y}}$ and 
$\| \psi^*(\bar{a}, \bar{y}) \ |  \ \theta^*(\bar{a}, \bar{y}) \|_{\bar{y}}$
are defined in $\mcA$ and
\[
r + \frac{\big|\psi(\bar{a}, \mcA) \cap \theta(\bar{a}, \mcA)\big|}{\big|\theta(\bar{a}, \mcA)\big|} \  \geq \
\frac{\big|\psi^*(\bar{a}, \mcA) \cap \theta^*(\bar{a}, \mcA)\big|}{\big|\theta^*(\bar{a}, \mcA)\big|},
\]
so
\begin{equation}\label{the inequality stated in terms of types}
r + \frac{\Big| \bigcup_{\iota =1}^m \bigcup_{j=1}^{m_\iota} p_{\iota,j}(\bar{a}, \mcA)\Big|}
{\Big| \bigcup_{\iota =1}^l \bigcup_{j=1}^{l_\iota} s_{\iota,j}(\bar{a}, \mcA)\Big|} \ \geq \
\frac{\Big| \bigcup_{\iota =1}^{m^*} \bigcup_{j=1}^{m_\iota^*} p_{\iota,j}^*(\bar{a}, \mcA)\Big|}
{\Big| \bigcup_{\iota=1}^{l^*} \bigcup_{j=1}^{l_\iota^*} s_{\iota,j}^*(\bar{a}, \mcA)\Big|}.
\end{equation}
Now $\bar{a}$ realizes exactly one of $t_1(\bar{x}), \ldots, t_1(\bar{x})$ 
and exactly one of $t_1^*(\bar{x}), \ldots, t_1^*(\bar{x})$ so
there are $1 \leq i \leq l$ and $1 \leq i' \leq l^*$ such that $\mcA \models t_i(\bar{a}) \wedge t_{i'}^*(\bar{a})$
and hence $t_i = t_{i'}^*$. 
If $r + \gamma_i \geq \gamma_{i'}^*$ then $i \in I$ and hence $\mcA \models \bigvee_{i \in I} t_i(\bar{a})$ so we are done.
Hence it remains to prove that $r + \gamma_i \geq \gamma_{i'}^*$.
We divide the argument into cases.

{\bf Case 1}: {\em Suppose that $i > m$ and $i' > m^*$.}
Then, by the definition of $\gamma_i$ and $\gamma_i^*$ (Definition~\ref{definition of gamma}),
we have $\gamma_i = \gamma_{i'}^* = 0$ so we get
$r + \gamma_i \geq \gamma_{i'}^*$.

{\bf Case 2}: {\em Suppose that $i \leq m$ and $i' > m^*$.}
Then $\gamma_{i'}^* = 0$ and as $\gamma_i$ is always nonnegative we get $r + \gamma_i \geq \gamma_{i'}^*$.

{\bf Case 3}: {\em Suppose that $i > m$ and $i' \leq m^*$.}
By Lemma~\ref{upper and lower bounds on p-i-j divided by s-i-j}~(a), assuming that $n$ is sufficiently large,
\[
\frac{\gamma_{i'}^*}{(1 + 2\varepsilon)^2} \ \leq \ \frac{\Big| \bigcup_{j=1}^{m_{i'}^*} p_{i',j}^*(\bar{a}, \mcA)\Big|}
{\Big| \bigcup_{j=1}^{l_{i'}^*} s_{i',j}^*(\bar{a}, \mcA)\Big|}.
\]
Since $i > m$ we have $p_{\iota, j}(\bar{a}, \mcA) = \es$ for every $1 \leq \iota \leq m$ and every $1 \leq j \leq m_\iota$, so
\[
 \frac{\Big| \bigcup_{\iota =1}^m \bigcup_{j=1}^{m_\iota} p_{\iota,j}(\bar{a}, \mcA)\Big|}
{\Big| \bigcup_{\iota =1}^l \bigcup_{j=1}^{l_\iota} s_{\iota,j}(\bar{a}, \mcA)\Big|} \ = \ 0.
\]
This together with~(\ref{the inequality stated in terms of types}) implies that
\begin{equation}\label{an equality in the third case}
r \  \geq \ \frac{\Big| \bigcup_{\iota =1}^{m^*} \bigcup_{j=1}^{m_\iota^*} p_{\iota,j}^*(\bar{a}, \mcA)\Big|}
{\Big| \bigcup_{\iota =1}^{l^*} \bigcup_{j=1}^{l_\iota^*} s_{\iota,j}^*(\bar{a}, \mcA)\Big|} 
\ = \ \frac{\Big| \bigcup_{j=1}^{m^*_{i'}} p_{i',j}^*(\bar{a}, \mcA)\Big|}
{\Big| \bigcup_{j=1}^{l^*_{i'}} s_{i',j}^*(\bar{a}, \mcA)\Big|} 
\ \geq \
\frac{\gamma_{i'}^*}{(1 + 2\varepsilon)^2}.
\end{equation}
If $r < \gamma_{i'}^*$ then, since $\varphi(\bar{x})$ is $\varepsilon$-noncritical, we get
$r < \gamma_{i'}^*/(1 + 2\varepsilon)^2$ which contradicts~(\ref{an equality in the third case}).
Hence $r \geq \gamma_{i'}^*$ and since $\gamma_i = 0$ (because $i > m$) we get $r + \gamma_i \geq \gamma_{i'}^*$.

{\bf Case 4}: {\em Suppose that $i \leq m$ and $i' \leq m^*$.}
Then~(\ref{the inequality stated in terms of types}) reduces to
\begin{equation}\label{r + p over s at least p* over s*}
r + \frac{\Big| \bigcup_{j=1}^{m_i} p_{i,j}(\bar{a}, \mcA)\Big|}
{\Big| \bigcup_{j=1}^{l_i} s_{i,j}(\bar{a}, \mcA)\Big|} \ \geq \
\frac{\Big| \bigcup_{j=1}^{m^*_{i'}} p_{i',j}^*(\bar{a}, \mcA)\Big|}
{\Big| \bigcup_{j=1}^{l^*_{i'}} s_{i',j}^*(\bar{a}, \mcA)\Big|}.
\end{equation}
Towards a contradiction, suppose that $r + \gamma_i < \gamma_{i'}^*$.
Since $\varphi(\bar{x})$ is assumed to be $\varepsilon$-noncritical we get
\begin{equation}\label{r + gamma-i < gamma-i-star with epsilon}
r + (1 + 2\varepsilon)^2 \gamma_i < \frac{\gamma_{i'}^*}{(1 + 2\varepsilon)^2}.
\end{equation}
Recall, from the definition of $d_i$ and $e_i$, that $d_i \leq e_i$.
We now consider two subcases and in each subcase we will derive a contradiction to~(\ref{r + p over s at least p* over s*}).

{\em Subcase 4(a): Suppose that $d_i = e_i$.}
By parts~(a) and~(b) of Lemma~\ref{upper and lower bounds on p-i-j divided by s-i-j},
\begin{align}\label{proportions in case 1}
\frac{\Big| \bigcup_{j=1}^{m_i} p_{i,j}(\bar{a}, \mcA)\Big|}
{\Big| \bigcup_{j=1}^{l_i} s_{i,j}(\bar{a}, \mcA)\Big|} \ \leq \ (1 + 2\varepsilon)^2 \gamma_i \quad \text{and} \quad
\frac{\Big| \bigcup_{j=1}^{m_{i'}^*} p_{i',j}^*(\bar{a}, \mcA)\Big|}
{\Big| \bigcup_{j=1}^{l_{i'}^*} s_{i',j}^*(\bar{a}, \mcA)\Big|} \ \geq \ 
\frac{\gamma_{i'}^*}{(1 + 2\varepsilon)^2}. 
\end{align}
From~(\ref{r + gamma-i < gamma-i-star with epsilon}) and~(\ref{proportions in case 1}) we get
\begin{align*}
r \ + \ \frac{\Big| \bigcup_{j=1}^{m_i} p_{i,j}(\bar{a}, \mcA)\Big|}{\Big| \bigcup_{j=1}^{l_i} s_{i,j}(\bar{a}, \mcA)\Big|}
\ \leq \
r + (1 + 2\varepsilon)^2 \gamma_i \ < \ \frac{\gamma_{i'}^*}{(1 + 2\varepsilon)^2}
\ \leq \ 
\frac{\Big| \bigcup_{j=1}^{m_{i'}^*} p_{i',j}^*(\bar{a}, \mcA)\Big|}
{\Big| \bigcup_{j=1}^{l_{i'}^*} s_{i',j}^*(\bar{a}, \mcA)\Big|}
\end{align*}
which contradicts~(\ref{r + p over s at least p* over s*}).

{\em Subcase 4(b): Suppose that $d_i < e_i$.}
Then Lemma~\ref{upper and lower bounds on p-i-j divided by s-i-j}~(c) gives
\begin{equation}\label{upper bound in case 1c}
\frac{\Big| \bigcup_{j=1}^{m_i} p_{i,j}(\bar{a}, \mcA)\Big|}{\Big| \bigcup_{j=1}^{l_i} s_{i,j}(\bar{a}, \mcA)\Big|} 
\ \leq \ \frac{C}{n}
\end{equation}
where the constant $C > 0$ depends only on the involved types.
Lemma~\ref{upper and lower bounds on p-i-j divided by s-i-j}~(a) gives
\[
\frac{\Big| \bigcup_{j=1}^{m_{i'}^*} p_{i',j}^*(\bar{a}, \mcA)\Big|}
{\Big| \bigcup_{j=1}^{l_{i'}^*} s_{i',j}^*(\bar{a}, \mcA)\Big|} \ \geq \ \frac{\gamma_{i'}^*}{(1 + 2\varepsilon)^2}.
\]
Since $d_i < e_i$ implies that $\gamma_i = 0$ it follows from~(\ref{r + gamma-i < gamma-i-star with epsilon}) that
$r < \gamma_{i'}^*/(1 + 2\varepsilon)^2$.
Note that the right hand term in~(\ref{upper bound in case 1c}) tends to 0 as $n$ tends to infinity.
So for all sufficiently large $n$ we have
\[
r \ + \ \frac{\Big| \bigcup_{j=1}^{m_i} p_{i,j}(\bar{a}, \mcA)\Big|}{\Big| \bigcup_{j=1}^{l_i} s_{i,j}(\bar{a}, \mcA)\Big|} 
\ \leq \ r + \frac{C}{n} \ < \ \frac{\gamma_{i'}^*}{(1 + 2\varepsilon)^2} 
\ \leq \ \frac{\Big| \bigcup_{j=1}^{m_{i'}^*} p_{i',j}^*(\bar{a}, \mcA)\Big|}
{\Big| \bigcup_{j=1}^{l_{i'}^*} s_{i',j}^*(\bar{a}, \mcA)\Big|}
\]
and this contradicts~(\ref{r + p over s at least p* over s*}).

{\em Now suppose that $\mcA \models \bigvee_{i \in I} t_i(\bar{a})$}, so $\mcA \models t_i(\bar{a})$ for some $i \in I$.
By Definition~\ref{definition of I} of $I$ there is $i' \in \{1, \ldots, l^*\}$ such that $t_i = t_{i'}^*$
and $r + \gamma_i \geq \gamma_{i'}^*$.
Since $\varphi(\bar{x})$ is an $\varepsilon$-noncritical formula it is, in particular, noncritical which implies that
$r + \gamma_i \neq \gamma_{i'}^*$ and hence
$r + \gamma_i > \gamma_{i'}^*$.
Since $\varphi(\bar{x})$ is $\varepsilon$-noncritical it follows that
\begin{equation}\label{the inequality involving r and gamma in the inverse direction}
r + \gamma_i / (1 + 2\varepsilon)^2  > \gamma_{i'}^* (1 + 2\varepsilon)^2.
\end{equation}
It suffices to prove that 
\begin{equation}\label{the inequality stated in terms of types again}
r + \frac{\Big| \bigcup_{\iota =1}^m \bigcup_{j=1}^{m_\iota} p_{\iota,j}(\bar{a}, \mcA)\Big|}
{\Big| \bigcup_{\iota =1}^l \bigcup_{j=1}^{l_\iota} s_{\iota,j}(\bar{a}, \mcA)\Big|} \ \geq \
\frac{\Big| \bigcup_{\iota =1}^{m^*} \bigcup_{j=1}^{m_\iota^*} p_{\iota,j}^*(\bar{a}, \mcA)\Big|}
{\Big| \bigcup_{\iota =1}^{l^*} \bigcup_{j=1}^{l_\iota^*} s_{\iota,j}^*(\bar{a}, \mcA)\Big|}.
\end{equation}
Again we divide the proof into cases.

{\bf Case 1}: {\em Suppose that $i' > m^*$.}
Then the term to the right of `$\geq$' in~(\ref{the inequality stated in terms of types again}) is zero, 
so~(\ref{the inequality stated in terms of types again})  holds.

{\bf Case 2}: {\em Suppose that $i > m$ and $i' \leq m^*$.}
Then the term immediately to the left of `$\geq$' in~(\ref{the inequality stated in terms of types again})
is zero, so we need to prove that
\begin{equation}\label{the reduction of the main inequality in case 2 of converse direction}
r \ \geq \ \frac{\Big| \bigcup_{\iota =1}^{m^*} \bigcup_{j=1}^{m_\iota^*} p_{\iota,j}^*(\bar{a}, \mcA)\Big|}
{\Big| \bigcup_{\iota =1}^{l^*} \bigcup_{j=1}^{l_\iota^*} s_{\iota,j}^*(\bar{a}, \mcA)\Big|}.
\end{equation}
From $i > m$ we get $\gamma_i = 0$ so~(\ref{the inequality involving r and gamma in the inverse direction})
reduces to
\begin{equation}\label{reduction of equality with r and gamma in case 3 of converse direction}
r > \gamma_{i'}^*(1 + 2\varepsilon)^2.
\end{equation}
Recall that from the definition it follows that $d_{i'}^* \leq e_{i'}^*$.

{\em Subcase 2(a): Suppose that $d_{i'}^* = e_{i'}^*$.}
Then using parts~(d) and (b) of Lemma~\ref{upper and lower bounds on p-i-j divided by s-i-j} we get
\[
\gamma_{i'}^*(1 + 2\varepsilon)^2 \ \geq \  
\frac{\Big| \bigcup_{j=1}^{m_{i'}^*} p_{i',j}^*(\bar{a}, \mcA)\Big|}
{\Big| \bigcup_{j=1}^{l_{i'}^*} s_{i',j}^*(\bar{a}, \mcA)\Big|}
\ = \ 
 \frac{\Big| \bigcup_{\iota =1}^{m^*} \bigcup_{j=1}^{m_\iota^*} p_{\iota,j}^*(\bar{a}, \mcA)\Big|}
{\Big| \bigcup_{\iota =1}^{l^*} \bigcup_{j=1}^{l_\iota^*} s_{\iota,j}^*(\bar{a}, \mcA)\Big|}.
\]
which together with~(\ref{reduction of equality with r and gamma in case 3 of converse direction})
gives~(\ref{the reduction of the main inequality in case 2 of converse direction}).

{\em Subcase 2(b): Suppose that $d_{i'}^* < e_{i'}^*$.}
Then parts~(d) and~(c) of Lemma~\ref{upper and lower bounds on p-i-j divided by s-i-j}
imply that 
\[
\frac{\Big| \bigcup_{\iota =1}^{m^*} \bigcup_{j=1}^{m_\iota^*} p_{\iota,j}^*(\bar{a}, \mcA)\Big|}
{\Big| \bigcup_{\iota =1}^{l^*} \bigcup_{j=1}^{l_\iota^*} s_{\iota,j}^*(\bar{a}, \mcA)\Big|}
\ = \ 
\frac{\Big| \bigcup_{j=1}^{m_{i'}^*} p_{i',j}^*(\bar{a}, \mcA)\Big|}
{\Big| \bigcup_{j=1}^{l_{i'}^*} s_{i',j}^*(\bar{a}, \mcA)\Big|}
\ \leq \
\frac{C}{n}
\]
for some constant $C > 0$ depending only on the involved types. 
Since $r > 0$ it follows that~(\ref{the reduction of the main inequality in case 2 of converse direction})
holds for all sufficiently large $n$.

{\bf Case 3}: {\em Suppose that $i \leq m$ and $i' \leq m^*$.}
Now~(\ref{the inequality stated in terms of types again}) is equivalent to  
\begin{equation}\label{r + p over s at least p* over s* again}
r + \frac{\Big| \bigcup_{j=1}^{m_i} p_{i,j}(\bar{a}, \mcA)\Big|}
{\Big| \bigcup_{j=1}^{l_i} s_{i,j}(\bar{a}, \mcA)\Big|} \ \geq \
\frac{\Big| \bigcup_{j=1}^{m^*_{i'}} p_{i',j}^*(\bar{a}, \mcA)\Big|}
{\Big| \bigcup_{j=1}^{l^*_{i'}} s_{i',j}^*(\bar{a}, \mcA)\Big|}.
\end{equation}
So it remains to prove~(\ref{r + p over s at least p* over s* again}).
By Lemma~\ref{upper and lower bounds on p-i-j divided by s-i-j}
and~(\ref{the inequality involving r and gamma in the inverse direction})
we have 
\begin{equation}\label{the almost final equation}
r \ + \ \frac{\Big| \bigcup_{j=1}^{m_i} p_{i,j}(\bar{a}, \mcA)\Big|}{\Big| \bigcup_{j=1}^{l_i} s_{i,j}(\bar{a}, \mcA)\Big|} 
\ \geq \
r + \frac{\gamma_i}{(1 + 2\varepsilon)^2} \ > \ (1 + 2\varepsilon)^2 \gamma_{i'}^*.
\end{equation}
If $d_{i'}^* = e_{i'}^*$, then, by Lemma~\ref{upper and lower bounds on p-i-j divided by s-i-j},
\[
(1 + 2\varepsilon)^2 \gamma_{i'}^* \ \geq \
\frac{\Big| \bigcup_{j=1}^{m^*_{i'}} p_{i',j}^*(\bar{a}, \mcA)\Big|}
{\Big| \bigcup_{j=1}^{l^*_{i'}} s_{i',j}^*(\bar{a}, \mcA)\Big|}
\]
which together with~(\ref{the almost final equation}) gives~(\ref{r + p over s at least p* over s* again}).

Now suppose that $d_{i'}^* < e_{i'}^*$.
Then
$\gamma_{i'}^* = 0$ and~(\ref{the almost final equation}) reduces to
\begin{equation}\label{the almost final equation in final case}
r \ + \ \frac{\Big| \bigcup_{j=1}^{m_i} p_{i,j}(\bar{a}, \mcA)\Big|}{\Big| \bigcup_{j=1}^{l_i} s_{i,j}(\bar{a}, \mcA)\Big|} 
\ \geq \
r + \frac{\gamma_i}{(1 + 2\varepsilon)^2} \ > \ 0
\end{equation}
Lemma~\ref{upper and lower bounds on p-i-j divided by s-i-j} gives
\[
\frac{\Big| \bigcup_{j=1}^{m^*_{i'}} p_{i',j}^*(\bar{a}, \mcA)\Big|}
{\Big| \bigcup_{j=1}^{l^*_{i'}} s_{i',j}^*(\bar{a}, \mcA)\Big|}
\ \leq \ 
\frac{C}{n} 
\]
This together with~(\ref{the almost final equation in final case}) gives~(\ref{r + p over s at least p* over s* again})
for all sufficiently large $n$.
This completes the proof of
Lemma~\ref{the case when I is nonempty}.
\hfill $\square$

\begin{lem}\label{the case when I is empty}
Suppose that $I = \es$.
Then for all sufficiently large $n$, all $\mcA \in \mbY_n$ and all $\bar{a} \in [n]^{|\bar{x}|}$,
\[
\mcA \not\models
\Big( r + \| \psi(\bar{a}, \bar{y}) \ |  \ \theta(\bar{a}, \bar{y}) \|_{\bar{y}} 
\ \geq \ 
\| \psi^*(\bar{a}, \bar{y}) \ |  \ \theta^*(\bar{a}, \bar{y}) \|_{\bar{y}} \Big).
\]
Hence the formula~(\ref{the formula}) is equivalent, in every such $\mcA$, to any contradictory quantifier-free formula.
\end{lem}

\noindent
{\bf Proof.}
Suppose that $I = \es$.
Suppose towards a contradiction that there are arbitrarily large $n$, $\mcA \in \mbY_n$ and  $\bar{a} \in [n]^{|\bar{x}|}$
such that
\[
\mcA \models
\Big( r + \| \psi(\bar{a}, \bar{y}) \ |  \ \theta(\bar{a}, \bar{y}) \|_{\bar{y}} 
\ \geq \ 
\| \psi^*(\bar{a}, \bar{y}) \ |  \ \theta^*(\bar{a}, \bar{y}) \|_{\bar{y}} \Big).
\]
Then we argue just as we did in the beginning of the proof of
Lemma~\ref{the case when I is nonempty} to get~(\ref{the inequality stated in terms of types})
and find 
$1 \leq i \leq l$ and $1 \leq i' \leq l^*$ such that $t_i = t_{i'}^*$.
Since $I = \es$ we must have $i \notin I$ and therefore $r + \gamma_i < \gamma_{i'}^*$.
Now we can continue to argue exactly as in the proof of Lemma~\ref{the case when I is nonempty}
to get a contradiction in each one of the cases 1--4 in that proof.
\hfill $\square$

\begin{rem}\label{the case when bar-x is empty}{\bf (The case when $\bar{x}$ is empty)} {\rm
Suppose now that $\bar{x}$ is empty, so the formula~(\ref{the formula}) becomes
\begin{equation}\label{when bar-x is empty}
\Big( r + \| \psi(\bar{y}) \ |  \ \theta(\bar{y}) \|_{\bar{y}} \ \geq \ 
\| \psi^*(\bar{y}) \ |  \ \theta^*(\bar{y}) \|_{\bar{y}} \Big),
\end{equation}
where we can assume that $\psi$, $\theta$, $\psi^*$ and $\theta^*$ are quantifier-free.
Then there are distinct types $p_i(\bar{y})$, $i = 1, \ldots, m$ and distinct types $s_i(\bar{y})$, $i = 1, \ldots, l$.
We can now define numbers $\gamma$ and $\gamma^*$ similarly as each $\gamma_i$ (and $\gamma_i^*$) was defined above.
We now get an analogoue of 
Lemma~\ref{upper and lower bounds on p-i-j divided by s-i-j}
which gives the same kind of upper and lower bounds of $\big| \bigcup_{i=1}^m p_i(\mcA) \big| \Big/ \big| \bigcup_{i=1}^l s_i(\mcA)\big|$
in terms of $\gamma$.
If $r + \gamma \geq \gamma^*$ then, by the noncriticality of~(\ref{when bar-x is empty}), we get
$r + \gamma > \gamma^*$ and by the $\varepsilon$-noncriticality of the same formula we get
$r + \gamma/(1 + 2\varepsilon)^2 >  \gamma^* (1 + 2\varepsilon)^2$.
Now we can argue similarly as in the ``converse direction'' in the proof of
Lemma~\ref{the case when I is nonempty}
and conclude that~(\ref{when bar-x is empty}) is true in all $\mcA \in \mbY_n$ for all sufficiently large $n$;
hence~(\ref{when bar-x is empty}) is equivalent to $\top$ in all such $\mcA$.
Now suppose that $r + \gamma < \gamma^*$ and suppose, towards a contradiction, that
there are arbitrarily large $n$ and $\mcA \in \mbY_n$ in which~(\ref{when bar-x is empty}) holds.
Then we can argue as in the first part of the proof of
Lemma~\ref{the case when I is nonempty}
and get a contradiction. Hence, for all sufficiently large $n$,~(\ref{when bar-x is empty}) is false in all $\mcA \in \mbY_n$;
consequently,~(\ref{when bar-x is empty}) is equivalent to $\neg\top$ in all such $\mcA$.
(The case when $\varphi$ has the form $\exists y \psi(\bar{y})$ is easier and analogous to the 
argument in the beginning of the proof of 
Proposition~\ref{quantifier elimination in a structure}
so this part is left to the reader.)
}\end{rem}

\noindent
Now the proof of Proposition~\ref{quantifier elimination in a structure}
is completed.
\hfill $\square$

\begin{defin}\label{definition of delta}{\rm
Define a function $\delta : \mbbN^+ \to \mbbR^{\geq 0}$ by
$\delta(n) = 5 \cdot \max\{\delta'(n), e^{-cn}\}$ where $c > 0$ is like in Lemma~\ref{Y-n has large probability}.
}\end{defin}

\begin{prop}\label{completing the inductive step with regard to the BN} {\bf (Completion of the induction step)}
Let $\mbY_n \subseteq \mbW_n$, $\varepsilon > 0$ and $\delta(n)$ be as in definitions~\ref{definition of Y-n}, \ref{definition of epsilon}
and~\ref{definition of delta}, respectively.
Then:
\begin{itemize}
\item[(1)] $\lim_{n\to\infty} \delta(n) = 0$.

\item[(2)] $\mbbP_n(\mbY_n) \geq 1 - \delta(n)$ for all sufficiently large $n$.

\item[(3)] For every complete atomic $\sigma$-type $p(\bar{x})$ with $|\bar{x}| \leq k$
there is a number which we denote 
$\msfP(p(\bar{x}))$ such that for all sufficiently large $n$ and all $\bar{a} \in [n]$ 
which realize the identity fragment of $p$,
\[
\big| \mbbP_n\big(\{\mcA \in \mbW_n : \mcA \models p(\bar{a})\}\big) \ -  \ \msfP(p(\bar{x})) \big| \ \leq \ \delta(n).
\]

\item[(4)] For every complete atomic $\sigma$-type $p(\bar{x}, \bar{y})$ with $|\bar{x}\bar{y}| \leq k$ and
$0 < \dim_{\bar{y}}(p(\bar{x}, y)) = |\bar{y}|$,
if $q(\bar{x}) = p \uhrc \bar{x}$ and $\msfP(q) > 0$, then for all sufficiently large $n$, every
$\mcA \in \mbY_n$ is $(p, q, \alpha/(1 + \varepsilon))$-saturated and $(p, q, \alpha(1 + \varepsilon))$-unsaturated
if $\alpha = \msfP(p(\bar{x}, \bar{y})) \ | \ \msfP(q(\bar{x}))$.

\item[(5)] For every $\varepsilon$-noncritical $\varphi(\bar{x}) \in CPL(\sigma)$ with $|\bar{x}| + \text{qr}(\varphi) \leq k$,
there is a quantifier-free $\sigma$-formula 
$\varphi^*(\bar{x})$ such that for all sufficiently large $n$ and all $\mcA \in \mbY_n$,
$\mcA \models \forall \bar{x} \big(\varphi(\bar{x}) \leftrightarrow \varphi^*(\bar{x})\big)$.
\end{itemize}
\end{prop}

\noindent
{\bf Proof.}
Parts~(1) and (2) follows from the definition of $\delta(n)$, Assumption~\ref{inductive assumptions}
and Lemma~\ref{Y-n has large probability}.
Part~(3) follows from Corollary~\ref{probability of p converges to P(p)}.
Part~(4) follows from 
Corollary~\ref{all members of Y-n are sufficiently saturated} 
and the definition of $\varepsilon$.
Part~(5) follows from Proposition~\ref{quantifier elimination in a structure}.
\hfill $\square$

\begin{cor}\label{corollary to completion of the inductive step}
Let $\varepsilon > 0$ be as in Definition~\ref{definition of epsilon}.\\
(a) If $\varphi(\bar{x}) \in CPL(\sigma)$ is $\varepsilon$-noncritical and  $|\bar{x}| + \mr{qr}(\varphi) \leq k$,
then there are $c > 0$ and  $0 \leq d \leq 1$ which is a sum of numbers of the form
$\msfP(p)$, where $p$ is a complete atomic $\sigma$-type, such that for every $m \in \mbbN^+$ and every $\bar{a} \in [m]^{|\bar{x}|}$
such that $\mcA \models \varphi(\bar{a})$ for some $\mcA \in \mbW_m$,
$\big| \mbbP_n(\varphi(\bar{a})) - d \big| \ \leq \ C\delta(n)$ for all sufficiently large $n$ where the constant $C$ depends only on $\varphi$.
Moreover, $d$ is $l$-critical for some $l$.\\
(b) If $\varphi \in CPL(\sigma)$ has no free variable, is $\varepsilon$-noncritical and $\mr{qr}(\varphi) \leq k$,
then either $\mbbP_n(\varphi) \leq \delta(n)$ for all sufficiently large $n$,
or $\mbbP_n(\varphi) \geq 1 - \delta(n)$ for all sufficiently large $n$.\\
(c) Suppose that for every $R \in \sigma \setminus \sigma'$,
if $\bar{x}$ is the sequence of free variables of $\chi_{R, i}$
then $|\bar{x}| + \mr{qr}(\chi_{R, i}) \leq k$.
Let $\mbbP_n^*$ be defined as $\mbbP_n$ except that we replace $\chi_{R, i}$ by $\chi_{R, i}^*$ in
Definition~\ref{definition of conditional probability distribution} where $\chi_{R, i}^*$ is a quantifier-free formula which his equivalent
to $\chi_{R, i}$ in every structure in $\mbY_n$ for all large enough $n$.
If $\varphi(\bar{x}) \in CPL(\sigma)$ is $\varepsilon$-noncritical, $|\bar{x}| + \mr{qr}(\varphi) \leq k$ and
$\mcA \models \varphi(\bar{a})$ for some $\mcA \in \mbW_m$ and some $m$, then
\[
\big| \mbbP_n^*(\varphi(\bar{a})) - \mbbP_n(\varphi(\bar{a})) \big| \  \leq \ \delta(n)
\ \text{ for all sufficiently large $n$.}
\]
\end{cor}

\noindent
{\bf Proof.}
(a) Suppose that $\varphi(\bar{x}) \in CPL(\sigma)$ is $\varepsilon$-noncritical and $|\bar{x}| + \mr{qr}(\varphi)  \leq k$.
By part~(5) of Proposition~\ref{completing the inductive step with regard to the BN}
$\varphi(\bar{x})$ is equivalent, in every $\mcA \in \mbY_n$ (for large enough $n$), to a quantifier-free formula $\varphi^*(\bar{x})$.
Then $\varphi^*(\bar{x})$ is equivalent to a disjunction of complete atomic $\sigma$-types 
$\bigvee_{i=1}^l p_i(\bar{x})$ where we assume that $p_i \neq p_j$ if $i \neq j$.
Suppose that $\mcA \models \varphi(\bar{a})$ for some $\mcA \in \mbW_m$ and some $m$.
Let $I$ be the set of indices $i$ such that $\mcA \models p_i(\bar{a})$ for some $\mcA \in \mbW_n$ and some $n$.
By assumption, $I \neq \es$. Let $d = \sum_{i \in I}\msfP(p_i)$.
By part~(3) of Proposition~\ref{completing the inductive step with regard to the BN},
we have 
$\big| \mbbP_n(\varphi^*(\bar{a})) - d \big| \ \leq \ |I|\delta(n)$ for all sufficiently large $n$,
and now~(a) follows from part~(2) of
Proposition~\ref{completing the inductive step with regard to the BN}.
It now follows from Definition~\ref{definition of critical number}
that $d$ is $l$-critical for some $l$.

(b) Suppose that $\varphi \in CPL(\sigma)$ has no free variable, is $\varepsilon$-noncritical and $\mr{qr}(\varphi) \leq k$.
By Proposition~\ref{completing the inductive step with regard to the BN}~(5),
there is a quantifier-free sentence $\varphi^*$ such that for all sufficiently large $n$ and all $\mcA \in \mbY_n$,
$\mcA \models \varphi \leftrightarrow \varphi^*$. 
Then $\varphi^*$ must be equivalent to $\bot$ or $\top$.
The conclusion of part~(b) now follows from parts~(1) and~(2) of
Proposition~\ref{completing the inductive step with regard to the BN}.

(c) Since $\chi_{R, i}$ is equivalent to $\chi_{R, i}^*$ in every $\mcA \in \mbY_n$ it follows from the definitions
of $\mbbP_n$ and $\mbbP_n^*$ that if
$\mcA \in \mbY_n$ then $\mbbP_n^*(\mcA) = \mbbP_n(\mcA)$.
It follows that if $\mbX_n \subseteq \mbY_n$ then $\mbbP_n^*(\mbX_n) = \mbbP_n(\mbX_n)$
and in particular
$\mbbP_n^*(\mbY_n) = \mbbP_n(\mbY_n)$.
Since $\mbbP_n^*(\mbW_n \setminus \mbY_n) = 1 - \mbbP_n^*(\mbY_n)$, and similarly for $\mbbP_n$, 
it follows that $\mbbP_n^*(\mbW_n \setminus \mbY_n) = \mbbP_n(\mbW_n \setminus \mbY_n)$.
From part~(2) of Proposition~\ref{completing the inductive step with regard to the BN}
we get $\mbbP_n^*(\mbW_n \setminus \mbY_n) = \mbbP_n(\mbW_n \setminus \mbY_n) \leq \delta(n)$.
Let $\mbX_n = \{\mcA \in \mbW_n : \mcA \models \varphi(\bar{a})\}$.
Then
\begin{align*}
&\mbbP_n^*(\mbX_n) \ \leq \ \mbbP_n^*(\mbX_n \ | \ \mbY_n)\mbbP_n^*(\mbY_n) + \delta(n) \ = \ 
\mbbP_n^*(\mbX_n \cap \mbY_n) + \delta(n)  = \\
&\mbbP_n(\mbX_n \cap \mbY_n) + \delta(n)
\ \leq \
\mbbP_n(\mbX_n) + \delta(n),
\end{align*}
and by similar reasoning
$
\mbbP_n(\mbX_n) \ \leq \ \mbbP_n^*(\mbX_n) + \delta(n).
$
\hfill $\square$

\section{Concluding remarks}

\noindent
The results of this article considers one particular formal logic and one type of lifted graphical model. 
Also, given these two things, 
choices have been made for example regarding exactly how to define a probability distribution on the set of structures with a
common finite domain. From the point of view machine learning and artificial intelligence, as well as mathematical curiosity, one could ask
a number of questions, of which I suggest a few below.

In finite model theory, theoretical computer science and linguistics a number of extensions of first-order logic have been considered \cite{Lib}. 
For example, a generic way of extending first-order logic is by adding one or more so-called generalized quantifiers \cite{Kai, KW}.
In machine learning, data mining and artificial intelligence a number of different (lifted) graphical models, including the popular
{\em Markov networks} \cite{DL, KMG} have been considered.
For which combinations of formal logical language and lifted graphical model do we get ``almost sure elimination of quantifiers''
and/or ``logical limit laws''? Do we get more expressive formalisms by using aggregation functions than if we use aggregation rules, or vice versa?
How do different combinations of formal language and graphical model relate to each other?
In what sense is a combination (formal language 1, graphical model 1) ``better'' than a combination
(formal language 2, graphical model 2)? What are reasonable candidates for the relation ``A is better/stronger than B''?
Some thoughts in this direction appear in
the last part of~\cite{DFKM}.

One can consider conditional probabilities which are not constant, but depend on the size of the set of elements (or tuples) satisfying
the condition in question. As a special case we have probabilities that depend on the size of the whole domain, as in
previous work on logical zero-one laws in random graphs \cite{SS, Spen}.)

What if the probability of a tuple $\bar{a}$ satisfying a relation is dependent on whether another tuple $\bar{b}$ satisfies the same relation
(as in \cite{KPR, Lyn} for example)?

A situation that seems natural in the context of artificial intelligence is to have an underlying fixed structure and on top of it relations that 
are ``governed'' by some probabilistic graphical model.
The underlying fixed structure could be represented by a $\tau$-structure $\mcA$ for some signature $\tau$.
For another signature $\sigma$ (disjoint from $\tau$) we could consider the set of expansions of $\mcA$ to $(\tau \cup \sigma)$-structures
where the probabilities of these extensions are governed by some probabilistic model {\em and} the underlying structure $\mcA$.
To formalize this using the set up of this article, one can modify $\mbW^\es_n$ in Definition~\ref{empty BN}
to contain exactly one $\tau$-structure with domain $[n]$ and $\mbW_n$ will be the set of all $(\tau \cup \sigma)$-structures
that expand the uniquen structure in $\mbW^\es_n$.
The definition of the probability distribution $\mbbP_n$ 
on $\mbW_n$ can now depend not only on the lifted Bayesian network $\mfG$ but also on the unique structure in $\mbW^\es_n$.
It seems obvious that, in order to get similar results as in this article, 
one needs to assume some sort of uniformity regarding the unique structure in $\mbW^\es_n$ for cofinitely many $n$.

\end{document}